%


\ifx\UsualIsLoaded\undefined
\let\UsualIsLoaded=\relax		

\font\fourteenrm=cmr12  scaled \magstep1
\font\fourteenbf=cmbx12 scaled \magstep1
\font\fourteentt=cmtt12 scaled \magstep1
\font\fourteensl=cmsl12 scaled \magstep1
\font\fourteensy=cmsy10 scaled \magstep2
\font\fourteeni=cmmi12  scaled \magstep1
\font\fourteenit=cmti12 scaled \magstep1
\font\fourteensc=cmcsc10 scaled \magstep2
\font\fourteenbit=cmssi12 scaled \magstep1
\font\fourteenbbb=msbm10 scaled \magstep2

\font\twelverm=cmr12
\font\twelvebf=cmbx12
\font\twelvett=cmtt12
\font\twelvesl=cmsl12
\font\twelvesy=cmsy10 scaled \magstep1
\font\twelvei=cmmi12
\font\twelveit=cmti12
\font\twelvesc=cmcsc10 scaled \magstep1
\font\twelvebit=cmssi12
\font\twelvebbb=msbm10 scaled \magstep1
\font\tenrm=cmr10
\font\tenbf=cmb10 
\font\tentt=cmtt10 
\font\tensl=cmsl10 
\font\tensy=cmsy10 
\font\teni=cmmi10 
\font\tenit=cmti10 
\font\tensc=cmcsc10 
\font\tenbit=cmssi10
\font\tenbbb=msbm10
\font\ninei=cmmi9 
\font\ninerm=cmr9

\font\ninesy=cmsy9 
\font\ninebbb=msbm10 at 9 pt
\font\eighti=cmmi8
\font\eightrm=cmr8

\font\eightsy=cmsy8
\font\eightbbb=msbm10 at 8 pt
\font\seveni=cmmi7 
\font\sevenrm=cmr7
\font\sevensy=cmsy7
 
%
%

%

\newfam\bbbfam

\def\tenpoint{
\def\rm{\fam0\tenrm}
\textfont0=\tenrm \scriptfont0=\eightrm \scriptscriptfont0=\sevenrm
\textfont1=\teni \scriptfont1=\eighti \scriptscriptfont1=\seveni
\textfont2=\tensy \scriptfont2=\eightsy \scriptscriptfont2=\sevensy
\textfont3=\tenex \scriptfont3=\tenex \scriptscriptfont3=\tenex

\textfont\itfam=\tenit
\def\it{\fam\itfam\tenit}

\textfont\slfam=\tensl
\def\sl{\fam\slfam\tensl}

\textfont\bffam=\tenbf
\def\bf{\fam\bffam\tenbf}

\textfont\ttfam=\tentt
\def\tt{\fam\ttfam\tentt}

\def\sc{\tensc}

\def\bit{\tenbit}

\def\bbb{\fam\bbbfam\twelvebbb}
\textfont\bbbfam=\tenbbb
\scriptfont\bbbfam=\eightbbb
\scriptscriptfont\bbbfam=\sevenbf
 
\normalbaselineskip=12pt

\setbox\strutbox=\hbox{\vrule height10pt depth4pt width0pt}%
\normalbaselines\rm}


\def\twelvepoint{
\def\rm{\fam0\twelverm}
\textfont0=\twelverm \scriptfont0=\ninerm \scriptscriptfont0=\sevenrm
\textfont1=\twelvei \scriptfont1=\ninei \scriptscriptfont1=\seveni
\textfont2=\twelvesy \scriptfont2=\ninesy \scriptscriptfont2=\sevensy
\textfont3=\tenex \scriptfont3=\tenex \scriptscriptfont3=\tenex

\textfont\itfam=\twelveit
\def\it{\fam\itfam\twelveit}

\textfont\slfam=\twelvesl
\def\sl{\fam\slfam\twelvesl}

\textfont\bffam=\twelvebf
\def\bf{\fam\bffam\twelvebf}

\textfont\ttfam=\twelvett
\def\tt{\fam\ttfam\twelvett}

\def\sc{\twelvesc}
\def\bit{\twelvebit}

\def\bbb{\fam\bbbfam\twelvebbb}
\textfont\bbbfam=\twelvebbb
\scriptfont\bbbfam=\ninebbb 
\scriptscriptfont\bbbfam=\sevenbf
 
\normalbaselineskip=14pt

\setbox\strutbox=\hbox{\vrule height10pt depth4pt width0pt}%
\normalbaselines\rm}


\def\fourteenpoint{
\def\rm{\fam0\fourteenrm}
\textfont0=\fourteenrm \scriptfont0=\twelverm \scriptscriptfont0=\tenrm
\textfont1=\fourteeni \scriptfont1=\twelvei \scriptscriptfont1=\teni
\textfont2=\fourteensy \scriptfont2=\twelvesy \scriptscriptfont2=\tensy
\textfont3=\tenex \scriptfont3=\tenex \scriptscriptfont3=\tenex

\textfont\itfam=\fourteenit
\def\it{\fam\itfam\fourteenit}

\textfont\slfam=\fourteensl
\def\sl{\fam\slfam\fourteensl}

\textfont\bffam=\fourteenbf
\def\bf{\fam\bffam\fourteenbf}

\textfont\ttfam=\fourteentt
\def\tt{\fam\ttfam\fourteentt}

\def\sc{\fourteensc}
\def\bit{\fourteenbit}

\def\bbb{\fam\bbbfam\twelvebbb}
\textfont\bbbfam=\fourteenbbb
\scriptfont\bbbfam=\tenrm 
\scriptscriptfont\bbbfam=\eightrm

\normalbaselineskip=16pt

\setbox\strutbox=\hbox{\vrule height10pt depth4pt width0pt}%
\normalbaselines\rm}

%
\twelvepoint

\abovedisplayskip 14pt plus 3pt minus 10pt%
\belowdisplayskip 14pt plus 3pt minus 10pt%
\abovedisplayshortskip 0pt plus 3pt%
\belowdisplayshortskip 8pt plus 3pt minus 5pt%
\parskip 3pt plus 1.5pt
\hsize=6.5in
\vsize=8.9in

%

\fi			


\def\e{{\rm e}}
\def\i{{\rm i}}
\def\union{\cup}		  
\def\boundary{\partial}		  
\def\interior{\mathaccent"0017}		  


\def\reals{{\bbb R}}		  
\def\R{{\reals}}

\def\C{{\bbb C}}			  
\def\re{{\rm Re\,}}		  
\def\im{{\rm Im\,}}		  
\def\arg{{\rm Arg}}		  
\def\conj#1{\overline{#1}}	  


\def\maps{\rightarrow}		  





\def\st{\thinspace \big| \thinspace }
\def\smbull{{\leavevmode\raise.25ex\hbox{$\scriptstyle\bullet$}}}
\def\smcirc{{\leavevmode\raise.25ex\hbox{$\scriptstyle\circ$}}}







\def\({\left (}
\def\){\right )}

\def\raggedcenter{\leftskip=0pt plus 1fill \rightskip=0pt plus 1fill}

\def\title#1{%
	\vbox to 0pt{}
	\vskip 12pt plus 4pt
	\vbox{\noindent\raggedcenter\let\\=\break
	      \fourteenpoint\bf #1}
	\vskip 14pt plus 10pt}


\newcount\secno			  
\global\secno=0			  
\newcount\thmno
\global\thmno=0
\outer\def\section#1\par{\vskip 0pt plus.25\vsize\penalty-250
	\vskip0pt plus -.25\vsize\bigskip\bigskip\bigskip\vskip\parskip
	\global\thmno=0
	\ifnum\secno<0
		{\centerline {\twelvebf #1}}
	\else
		\advance\secno by 1
		{\centerline {\twelvebf \the\secno: #1}}
	\fi
	\message{<\the\secno: #1>}
	\nobreak\medskip\noindent}
%
\def\proclaim #1. #2\par{\bigbreak
	\noindent{\sc#1.\enspace}{\sl#2\par}
	\ifdim\lastskip<\medskipamount 
		\removelastskip\penalty55\medskip
	\fi}
%
%
\outer\def\theorem #1 #2. #3\par{\medbreak
	\advance\thmno by 1
	\ifnum\secno<1 \edef\thmid{\number\thmno}
		\else  \edef\thmid{{\number\secno}.{\number\thmno}}
	\fi
	\edef#2{{#1 \thmid}}
	\proclaim {#1 \thmid}. {#3}\par\medskip
	}
%
%
\outer\def\figurename#1{\advance\thmno by 1
	\ifnum\secno<1 \edef\thmid{\number\thmno}
		\else  \edef\thmid{{\number\secno}.{\number\thmno}}
	\fi
	\edef#1{{Figure \thmid}}
	}

\def\proof{\smallbreak	\noindent{\it Proof.\enspace}}
\def\QED{\hfill\hbox{\rlap{$\sqcup$}\raise .2ex\hbox{$\sqcap$}}}

\def\definition{\smallbreak\noindent{\sc Definition.\ }}

\def\remark{\smallbreak\noindent{\sc Remark.\ }}
\def\subsection#1{\bigbreak\noindent{\bf #1.\ }}
\input psfig
\let\macrosLoaded=\relax                      
%
%
\def\shorttitle{Dynamics of Certain Non-Conformal Maps of $\C$}
\def\authors{Bielefeld, Sutherland, Tangerman, and Veerman}
%
\global\headline={%
    \ifnum\pageno > 1 {%
	  \ifodd\pageno
	     {\tenpoint \hfil{\sc \shorttitle}  \hfil}
          \else
	     {\tenpoint \hfil {\sc \authors} \hfil}
	  \fi}
    \else
		\hss\vbox{}\hss
    \fi
}

%
\def\cite#1{[#1]}

%
\def\CLocus#1{{\cal C}_{#1}}
\def\ConnLocus{\CLocus{\alpha}}
\def\MSet{\CLocus{1}}
\newif\ifPreviewSyms
\PreviewSymsfalse
\ifPreviewSyms
 \def\OneRep{\circ}
 \def\TwoRep{\lower.4ex\hbox{${\mathop{\kern 0pt\smcirc}\limits^\circ}$}}
 \def\Saddle{+}
 \def\OneAttr{\bullet}
 \def\TwoAttr{\lower.4ex\hbox{${\mathop{\kern 0pt\smbull}\limits^\bullet}$}}
\else
 \def\OneRep{\lower.2ex\psfig{figure=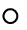,height=1.5ex,silent=}}
 \def\TwoRep{\lower.2ex\psfig{figure=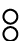,height=1.5ex,silent=}}
 \def\Saddle{\lower.2ex\psfig{figure=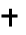,height=1.5ex,silent=}}
 \def\OneAttr{\lower.2ex\psfig{figure=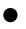,height=1.5ex,silent=}}
 \def\TwoAttr{\lower.2ex\psfig{figure=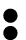,height=1.5ex,silent=}}
\fi

%
\def\advfignumber{%
        \global\advance\thmno by 1
	\ifnum\secno<1 \edef\fignumber{\number\thmno}
		\else  \edef\fignumber{{\number\secno}.{\number\thmno}}
	\fi
}
\def\shortcaption#1{%
    \advfignumber
    \smallskip
    \centerline{\tenpoint\noindent{\sc Figure \fignumber:}\quad #1}\medskip}
\def\longcaption#1{%
        \advfignumber
	\smallskip
	\centerline{\vbox{\hsize=.8\hsize\tenpoint\noindent
		{\sc Figure \fignumber:}\quad #1}}\medskip}
\def\widecaption#1{%
        \advfignumber
	\smallskip
	\centerline{\vbox{\tenpoint\noindent
		{\sc Figure \fignumber:}\quad #1}}\medskip}
\def\namefigure#1{%
    \global\edef#1{Figure~\fignumber}}
\def\namenextfigure#1{%
    \advfignumber
    \namefigure{#1}
    \global\advance\thmno by -1
}

\long\gdef\ignore#1{}
\def\emptyset{/\kern -0.67em{\scriptstyle\bigcirc}}

\chardef \other = 12

\def\makeactive#1{\catcode`#1 = \active\ignorespaces}
{
	\makeactive\^^M %
	\gdef\obeywhitespace{%
		\makeactive\^^M %
		\let^^M = \newline %
		\aftergroup\removebox %
		\obeyspaces %
	}%
}
\def\newline{\par\indent}
\def\removebox{\setbox0=\lastbox}

\def\{}

\def\BoxIt#1#2{
	\vbox{\hrule
	\hbox{\vrule\kern#2\vbox{\kern#2#1\kern#2}\kern#2\vrule}
    \hrule}}


%
\def\IMSmarkvadjust{0 pt}
\def\IMSmarkhadjust{0 pt}
\def\IMSmarkhpadding{0 pt}
\def\IMSpubltext{Published in modified form:}
\def\SBIMSMark#1#2#3{
 \font\SBF=cmss10 at 10 true pt
 \font\SBI=cmssi10 at 10 true pt
 \setbox0=\hbox{\SBF \hbox to \IMSmarkhpadding{\relax}
                Stony Brook IMS Preprint \##1}
 \setbox2=\hbox to \wd0{\hfil \SBI #2}
 \setbox4=\hbox to \wd0{\hfil \SBI #3}
 \setbox6=\hbox to \wd0{\hss
             \vbox{\hsize=\wd0 \parskip=0pt \baselineskip=10 true pt
                   \copy0 \break%
                   \copy2 \break%
                   \copy4 \break}}
 \dimen0=\ht6   \advance\dimen0 by \vsize \advance\dimen0 by 8 true pt
                \advance\dimen0 by -\pagetotal
	        \advance\dimen0 by \IMSmarkvadjust
 \dimen2=\hsize \advance\dimen2 by .25 true in
	        \advance\dimen2 by \IMSmarkhadjust

%
%
  \openin2=publishd.tex
  \ifeof2\setbox0=\hbox to 0pt{}
  \else 
     \setbox0=\hbox to 3.1 true in{
                \vbox to \ht6{\hsize=3 true in \parskip=0pt  \noindent  
                {\SBI \IMSpubltext}\hfil\break
                {\it Experimental Math.}~{\bf 2} (1993), 281--300 
                \vfill}}
  \fi
  \closein2
  \ht0=0pt \dp0=0pt
 \ht6=0pt \dp6=0pt
 \setbox8=\vbox to \dimen0{\vfill \hbox to \dimen2{\copy0 \hss \copy6}}
 \ht8=0pt \dp8=0pt \wd8=0pt
 \copy8
 \message{*** Stony Brook IMS Preprint #1, #2. #3 ***}
}

\SBIMSMark{1991/18}{November 1991}{revised August 1993}
\newif\ifSingleSecs

%
%
	\ifx\macrosLoaded\undefined 
        
%


\ifx\UsualIsLoaded\undefined
\let\UsualIsLoaded=\relax		

\font\fourteenrm=cmr12  scaled \magstep1
\font\fourteenbf=cmbx12 scaled \magstep1
\font\fourteentt=cmtt12 scaled \magstep1
\font\fourteensl=cmsl12 scaled \magstep1
\font\fourteensy=cmsy10 scaled \magstep2
\font\fourteeni=cmmi12  scaled \magstep1
\font\fourteenit=cmti12 scaled \magstep1
\font\fourteensc=cmcsc10 scaled \magstep2
\font\fourteenbit=cmssi12 scaled \magstep1
\font\fourteenbbb=msbm10 scaled \magstep2

\font\twelverm=cmr12
\font\twelvebf=cmbx12
\font\twelvett=cmtt12
\font\twelvesl=cmsl12
\font\twelvesy=cmsy10 scaled \magstep1
\font\twelvei=cmmi12
\font\twelveit=cmti12
\font\twelvesc=cmcsc10 scaled \magstep1
\font\twelvebit=cmssi12
\font\twelvebbb=msbm10 scaled \magstep1
\font\tenrm=cmr10
\font\tenbf=cmb10 
\font\tentt=cmtt10 
\font\tensl=cmsl10 
\font\tensy=cmsy10 
\font\teni=cmmi10 
\font\tenit=cmti10 
\font\tensc=cmcsc10 
\font\tenbit=cmssi10
\font\tenbbb=msbm10
\font\ninei=cmmi9 
\font\ninerm=cmr9

\font\ninesy=cmsy9 
\font\ninebbb=msbm10 at 9 pt
\font\eighti=cmmi8
\font\eightrm=cmr8

\font\eightsy=cmsy8
\font\eightbbb=msbm10 at 8 pt
\font\seveni=cmmi7 
\font\sevenrm=cmr7
\font\sevensy=cmsy7
 
%
%

%

\newfam\bbbfam

\def\tenpoint{
\def\rm{\fam0\tenrm}
\textfont0=\tenrm \scriptfont0=\eightrm \scriptscriptfont0=\sevenrm
\textfont1=\teni \scriptfont1=\eighti \scriptscriptfont1=\seveni
\textfont2=\tensy \scriptfont2=\eightsy \scriptscriptfont2=\sevensy
\textfont3=\tenex \scriptfont3=\tenex \scriptscriptfont3=\tenex

\textfont\itfam=\tenit
\def\it{\fam\itfam\tenit}

\textfont\slfam=\tensl
\def\sl{\fam\slfam\tensl}

\textfont\bffam=\tenbf
\def\bf{\fam\bffam\tenbf}

\textfont\ttfam=\tentt
\def\tt{\fam\ttfam\tentt}

\def\sc{\tensc}

\def\bit{\tenbit}

\def\bbb{\fam\bbbfam\twelvebbb}
\textfont\bbbfam=\tenbbb
\scriptfont\bbbfam=\eightbbb
\scriptscriptfont\bbbfam=\sevenbf
 
\normalbaselineskip=12pt

\setbox\strutbox=\hbox{\vrule height10pt depth4pt width0pt}%
\normalbaselines\rm}


\def\twelvepoint{
\def\rm{\fam0\twelverm}
\textfont0=\twelverm \scriptfont0=\ninerm \scriptscriptfont0=\sevenrm
\textfont1=\twelvei \scriptfont1=\ninei \scriptscriptfont1=\seveni
\textfont2=\twelvesy \scriptfont2=\ninesy \scriptscriptfont2=\sevensy
\textfont3=\tenex \scriptfont3=\tenex \scriptscriptfont3=\tenex

\textfont\itfam=\twelveit
\def\it{\fam\itfam\twelveit}

\textfont\slfam=\twelvesl
\def\sl{\fam\slfam\twelvesl}

\textfont\bffam=\twelvebf
\def\bf{\fam\bffam\twelvebf}

\textfont\ttfam=\twelvett
\def\tt{\fam\ttfam\twelvett}

\def\sc{\twelvesc}
\def\bit{\twelvebit}

\def\bbb{\fam\bbbfam\twelvebbb}
\textfont\bbbfam=\twelvebbb
\scriptfont\bbbfam=\ninebbb 
\scriptscriptfont\bbbfam=\sevenbf
 
\normalbaselineskip=14pt

\setbox\strutbox=\hbox{\vrule height10pt depth4pt width0pt}%
\normalbaselines\rm}


\def\fourteenpoint{
\def\rm{\fam0\fourteenrm}
\textfont0=\fourteenrm \scriptfont0=\twelverm \scriptscriptfont0=\tenrm
\textfont1=\fourteeni \scriptfont1=\twelvei \scriptscriptfont1=\teni
\textfont2=\fourteensy \scriptfont2=\twelvesy \scriptscriptfont2=\tensy
\textfont3=\tenex \scriptfont3=\tenex \scriptscriptfont3=\tenex

\textfont\itfam=\fourteenit
\def\it{\fam\itfam\fourteenit}

\textfont\slfam=\fourteensl
\def\sl{\fam\slfam\fourteensl}

\textfont\bffam=\fourteenbf
\def\bf{\fam\bffam\fourteenbf}

\textfont\ttfam=\fourteentt
\def\tt{\fam\ttfam\fourteentt}

\def\sc{\fourteensc}
\def\bit{\fourteenbit}

\def\bbb{\fam\bbbfam\twelvebbb}
\textfont\bbbfam=\fourteenbbb
\scriptfont\bbbfam=\tenrm 
\scriptscriptfont\bbbfam=\eightrm

\normalbaselineskip=16pt

\setbox\strutbox=\hbox{\vrule height10pt depth4pt width0pt}%
\normalbaselines\rm}

%
\twelvepoint

\abovedisplayskip 14pt plus 3pt minus 10pt%
\belowdisplayskip 14pt plus 3pt minus 10pt%
\abovedisplayshortskip 0pt plus 3pt%
\belowdisplayshortskip 8pt plus 3pt minus 5pt%
\parskip 3pt plus 1.5pt
\hsize=6.5in
\vsize=8.9in

%

\fi			


\def\e{{\rm e}}
\def\i{{\rm i}}
\def\union{\cup}		  
\def\boundary{\partial}		  
\def\interior{\mathaccent"0017}		  


\def\reals{{\bbb R}}		  
\def\R{{\reals}}

\def\C{{\bbb C}}			  
\def\re{{\rm Re\,}}		  
\def\im{{\rm Im\,}}		  
\def\arg{{\rm Arg}}		  
\def\conj#1{\overline{#1}}	  


\def\maps{\rightarrow}		  





\def\st{\thinspace \big| \thinspace }
\def\smbull{{\leavevmode\raise.25ex\hbox{$\scriptstyle\bullet$}}}
\def\smcirc{{\leavevmode\raise.25ex\hbox{$\scriptstyle\circ$}}}







\def\({\left (}
\def\){\right )}

\def\raggedcenter{\leftskip=0pt plus 1fill \rightskip=0pt plus 1fill}

\def\title#1{%
	\vbox to 0pt{}
	\vskip 12pt plus 4pt
	\vbox{\noindent\raggedcenter\let\\=\break
	      \fourteenpoint\bf #1}
	\vskip 14pt plus 10pt}


\newcount\secno			  
\global\secno=0			  
\newcount\thmno
\global\thmno=0
\outer\def\section#1\par{\vskip 0pt plus.25\vsize\penalty-250
	\vskip0pt plus -.25\vsize\bigskip\bigskip\bigskip\vskip\parskip
	\global\thmno=0
	\ifnum\secno<0
		{\centerline {\twelvebf #1}}
	\else
		\advance\secno by 1
		{\centerline {\twelvebf \the\secno: #1}}
	\fi
	\message{<\the\secno: #1>}
	\nobreak\medskip\noindent}
%
\def\proclaim #1. #2\par{\bigbreak
	\noindent{\sc#1.\enspace}{\sl#2\par}
	\ifdim\lastskip<\medskipamount 
		\removelastskip\penalty55\medskip
	\fi}
%
%
\outer\def\theorem #1 #2. #3\par{\medbreak
	\advance\thmno by 1
	\ifnum\secno<1 \edef\thmid{\number\thmno}
		\else  \edef\thmid{{\number\secno}.{\number\thmno}}
	\fi
	\edef#2{{#1 \thmid}}
	\proclaim {#1 \thmid}. {#3}\par\medskip
	}
%
%
\outer\def\figurename#1{\advance\thmno by 1
	\ifnum\secno<1 \edef\thmid{\number\thmno}
		\else  \edef\thmid{{\number\secno}.{\number\thmno}}
	\fi
	\edef#1{{Figure \thmid}}
	}

\def\proof{\smallbreak	\noindent{\it Proof.\enspace}}
\def\QED{\hfill\hbox{\rlap{$\sqcup$}\raise .2ex\hbox{$\sqcap$}}}

\def\definition{\smallbreak\noindent{\sc Definition.\ }}

\def\remark{\smallbreak\noindent{\sc Remark.\ }}
\def\subsection#1{\bigbreak\noindent{\bf #1.\ }}
\input psfig
\let\macrosLoaded=\relax                      
%
%
\def\shorttitle{Dynamics of Certain Non-Conformal Maps of $\C$}
\def\authors{Bielefeld, Sutherland, Tangerman, and Veerman}
%
\global\headline={%
    \ifnum\pageno > 1 {%
	  \ifodd\pageno
	     {\tenpoint \hfil{\sc \shorttitle}  \hfil}
          \else
	     {\tenpoint \hfil {\sc \authors} \hfil}
	  \fi}
    \else
		\hss\vbox{}\hss
    \fi
}

%
\def\cite#1{[#1]}

%
\def\CLocus#1{{\cal C}_{#1}}
\def\ConnLocus{\CLocus{\alpha}}
\def\MSet{\CLocus{1}}
\newif\ifPreviewSyms
\PreviewSymsfalse
\ifPreviewSyms
 \def\OneRep{\circ}
 \def\TwoRep{\lower.4ex\hbox{${\mathop{\kern 0pt\smcirc}\limits^\circ}$}}
 \def\Saddle{+}
 \def\OneAttr{\bullet}
 \def\TwoAttr{\lower.4ex\hbox{${\mathop{\kern 0pt\smbull}\limits^\bullet}$}}
\else
 \def\OneRep{\lower.2ex\psfig{figure=figures/OneRep.ps,height=1.5ex,silent=}}
 \def\TwoRep{\lower.2ex\psfig{figure=figures/TwoRep.ps,height=1.5ex,silent=}}
 \def\Saddle{\lower.2ex\psfig{figure=figures/Saddle.ps,height=1.5ex,silent=}}
 \def\OneAttr{\lower.2ex\psfig{figure=figures/OneAttr.ps,height=1.5ex,silent=}}
 \def\TwoAttr{\lower.2ex\psfig{figure=figures/TwoAttr.ps,height=1.5ex,silent=}}
\fi

%
\def\advfignumber{%
        \global\advance\thmno by 1
	\ifnum\secno<1 \edef\fignumber{\number\thmno}
		\else  \edef\fignumber{{\number\secno}.{\number\thmno}}
	\fi
}
\def\shortcaption#1{%
    \advfignumber
    \smallskip
    \centerline{\tenpoint\noindent{\sc Figure \fignumber:}\quad #1}\medskip}
\def\longcaption#1{%
        \advfignumber
	\smallskip
	\centerline{\vbox{\hsize=.8\hsize\tenpoint\noindent
		{\sc Figure \fignumber:}\quad #1}}\medskip}
\def\widecaption#1{%
        \advfignumber
	\smallskip
	\centerline{\vbox{\tenpoint\noindent
		{\sc Figure \fignumber:}\quad #1}}\medskip}
\def\namefigure#1{%
    \global\edef#1{Figure~\fignumber}}
\def\namenextfigure#1{%
    \advfignumber
    \namefigure{#1}
    \global\advance\thmno by -1
}

\long\gdef\ignore#1{}
\def\emptyset{/\kern -0.67em{\scriptstyle\bigcirc}}

\chardef \other = 12

\def\makeactive#1{\catcode`#1 = \active\ignorespaces}
{
	\makeactive\^^M %
	\gdef\obeywhitespace{%
		\makeactive\^^M %
		\let^^M = \newline %
		\aftergroup\removebox %
		\obeyspaces %
	}%
}
\def\newline{\par\indent}
\def\removebox{\setbox0=\lastbox}

\def\{}

\def\BoxIt#1#2{
	\vbox{\hrule
	\hbox{\vrule\kern#2\vbox{\kern#2#1\kern#2}\kern#2\vrule}
    \hrule}}


\fi
\secno=-1

\title{Dynamics of Certain Non-Conformal \\
       Degree Two Maps of the Plane}

\centerline {Ben Bielefeld}
\centerline {Scott Sutherland}
\centerline {Folkert Tangerman}
\centerline {Peter Veerman}

\section{Introduction}

We consider a family of maps which are similar to quadratic maps in
that they are degree two branched covers of the Riemann sphere, but
not in general conformal.  In particular, for $\alpha$ real and
fixed, we study maps $f_c$ which are given in polar coordinates by 
	$$f_c(r \e^{\i\theta}) = r^{2\alpha}\e^{2\i\theta} + c.$$
For $\alpha=1$, this is the usual quadratic family ($z \mapsto z^2+c$),
which has been extensively studied  and is fairly well understood.  For
$\alpha$ different from one, $f_c$ is only quasi-conformal, and very
different behavior can occur, although there are many strong
similarities to the conformal case.  It is our goal to determine which
results for the quadratic family  can be generalized to maps which are
topologically similar (and when $\alpha$ is close to 1, close to
quadratic), and where such results break down.

In the quadratic family, the orbit of the critical point completely
determines the dynamics.  This is not the case for the maps $f_c$: for
example, we have found periodic attractors which do not attract the
critical point.  For certain parameter values, the dynamics is
dominated by two dimensional real behavior: periodic saddle points,
invariant circles, and so on. 

Another striking difference with the quadratic family is the existence of
smooth Julia sets.  In the conformal case, the only smooth Julia sets are
the segment $[-2,2]$ (for the map $z\mapsto z^2-2$) and the unit circle (for
the map $z\mapsto z^2$). The corresponding Julia sets for $f_c$ are also
smooth, but there are more: we use structural stability techniques to show
that for any $\alpha <1$, the Julia set is $C^k$-smooth for all $c$-values
sufficiently near $0$.

We also study the connectedness locus (the analogue of the Mandelbrot set),
and the bifurcations which occur in the $c$-plane.  Numerical evidence
strongly suggests that the connectedness locus is always connected, and
never locally connected for $\alpha \ne 1$.  Furthermore, the bifurcations
which occur as the parameter $c$ varies are considerably more complicated
than those in the conformal case, although there are many similarities.  We
discuss these issues at some length in sections~4 and~5.

	\ifSingleSecs \vfil\eject \fi

	\ifx\macrosLoaded\undefined 
        
\fi
\secno=0  

\section {Definitions and Elementary Results}

For $\alpha > 0$, consider the map $Q_{\alpha}$:
$$\eqalign {
  Q_{\alpha}\(r \e^{\i \theta}\) = r^{\alpha}\e^{\i\theta} 
	\quad&\quad \hbox{in polar coordinates}\cr
  Q_{\alpha}(z)=z^{{\alpha+1}\over{2}} {\conj{z}}^{{\alpha-1}\over{2}}
	&\quad \hbox{in $(z,\conj{z})$ coordinates, for appropriate branches
		    of the powers}
}$$
Note that the family $\{ Q_{\alpha} \}$ is a one parameter group: 
$Q_{\alpha} \circ Q_{\beta} = Q_{\alpha \beta}.$

\theorem Proposition \QisQC. 
Each $Q_{\alpha}$ is a quasi-conformal homeomorphism of the Riemann sphere
of constant dilatation $\max(\alpha, 1/\alpha)$.

For the definition of quasi-conformal, refer to \cite{L}.

\proof
Consider $Q = Q_{\alpha}$. A straightforward computation yields that its
dilitation is 
$$
 {\rm dil}_Q(z) =  {{|\partial_zQ|  + |\partial_{\bar{z}}Q|}\over 
	      {|\partial_z Q| - |\partial_{\bar{z}}Q|}}		
 	  =  {{\alpha + 1 + |\alpha - 1|}\over
	      {\alpha + 1 - |\alpha - 1|}}
          =  \max(\alpha, 1/\alpha).
$$
\vskip-1.7\baselineskip\QED

Denote by $P_c$ the quadratic map on $\C$: 
\ $\displaystyle{P_c(z) = z^2 + c}$,
and let $f_{\alpha, c}$ represent the composition $P_c \circ Q_\alpha$. Then
$$
 f_{\alpha, c} = \cases{
     |z|^{2\alpha-2} z^2 + c		 & or\cr
     z^{\alpha+1}{\conj z}^{\alpha-1} +c & in $(z,\conj{z})$ coordinates,\cr
     r^{2\alpha} \e^{\i 2 \theta} + c    & in polar coordinates.\cr}
$$

For any $\alpha > 0$ and any $c \in {\C}$, the map $f_{\alpha, c}$ is a
branched cover of $\C$ with a single branch point, the origin, where the map
is ramified of degree two. It extends to the Riemann sphere with a branch
point at $\infty$ of degree two. Throughout this paper we always assume
that $\alpha > 1/2$. This guarantees that the dynamics near infinity
is always the same: the point $\infty$ is attracting. Moreover,
when $\alpha > 1/2$, each $f_{\alpha, c}$ is at least once differentiable 
everywhere.

For any map $f_{\alpha, c}$, define its {\bit filled-in Julia set}
$K(\alpha, c)$ as the set of points whose orbits under $f_{\alpha, c}$ do
not accumulate at $\infty$.  We define the {\bit Julia set} $J(\alpha, c)$
as the the set of points which have no neighborhood in which the iterates of
$f_{\alpha, c}$ form an equicontinuous family in the spherical metric.
Because $f_{\alpha, c}$ is an open map, the Julia set can be split up into
two completely invariant (i.e.\ forward and backward invariant) subsets
$\boundary K(\alpha, c)$ and $\Delta(\alpha, c) =
J(\alpha,c)-\boundary K(\alpha, c)$.  
In the case when $\alpha=1$ the set $\Delta(\alpha, c)$ is empty, but in
general it is nonempty.  For instance, the set $\Delta(\alpha, c)$ may
contain stable manifolds of periodic saddle points.

\midinsert{
  \hbox to \hsize{
    \vsize=.28\hsize
    \hfil
    \vbox to \vsize{\vfil
        \psfig{figure=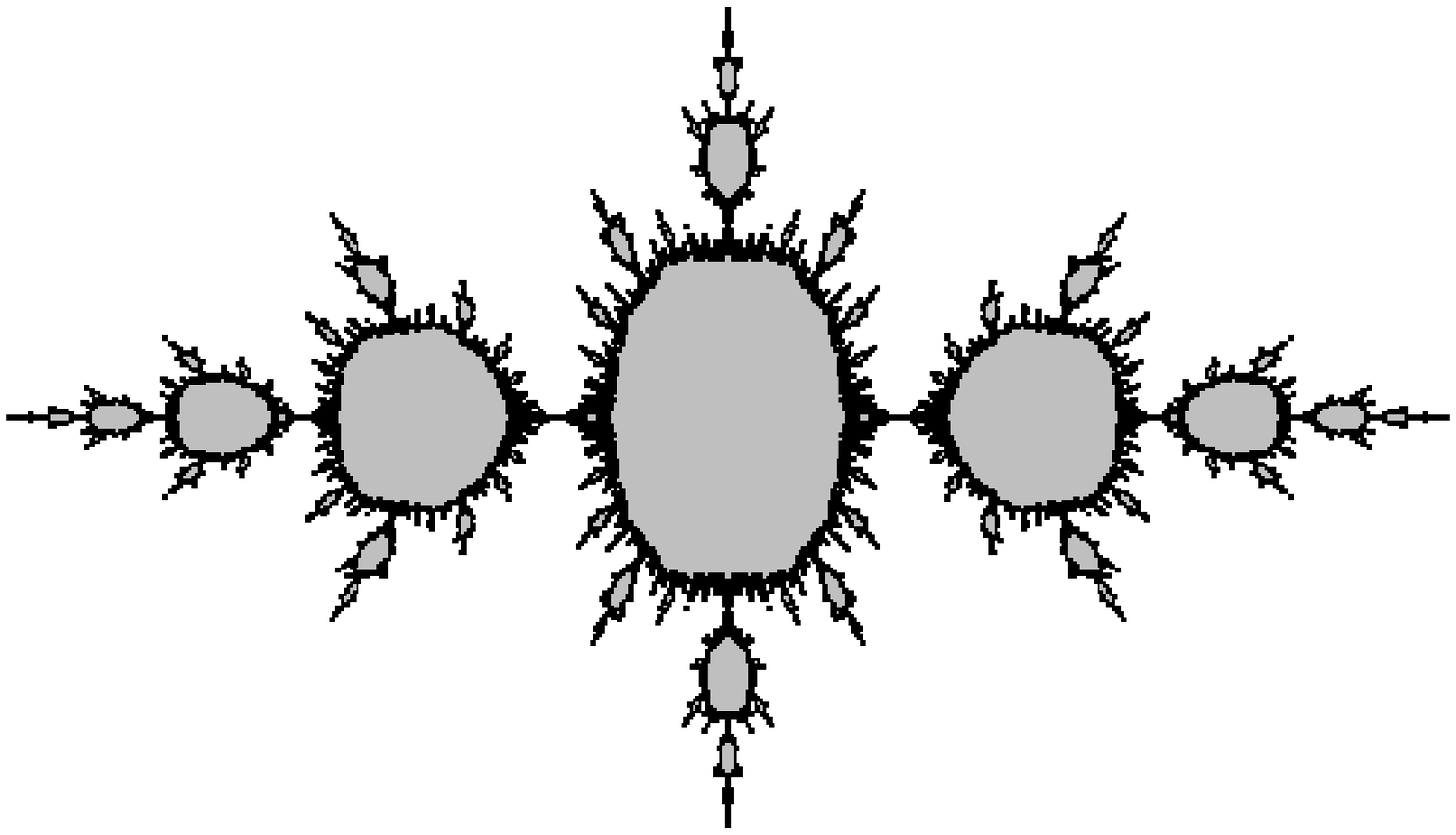,width=.45\hsize}
	\vfil}
    \hfil
    \vbox to \vsize{\vfil
        \psfig{figure=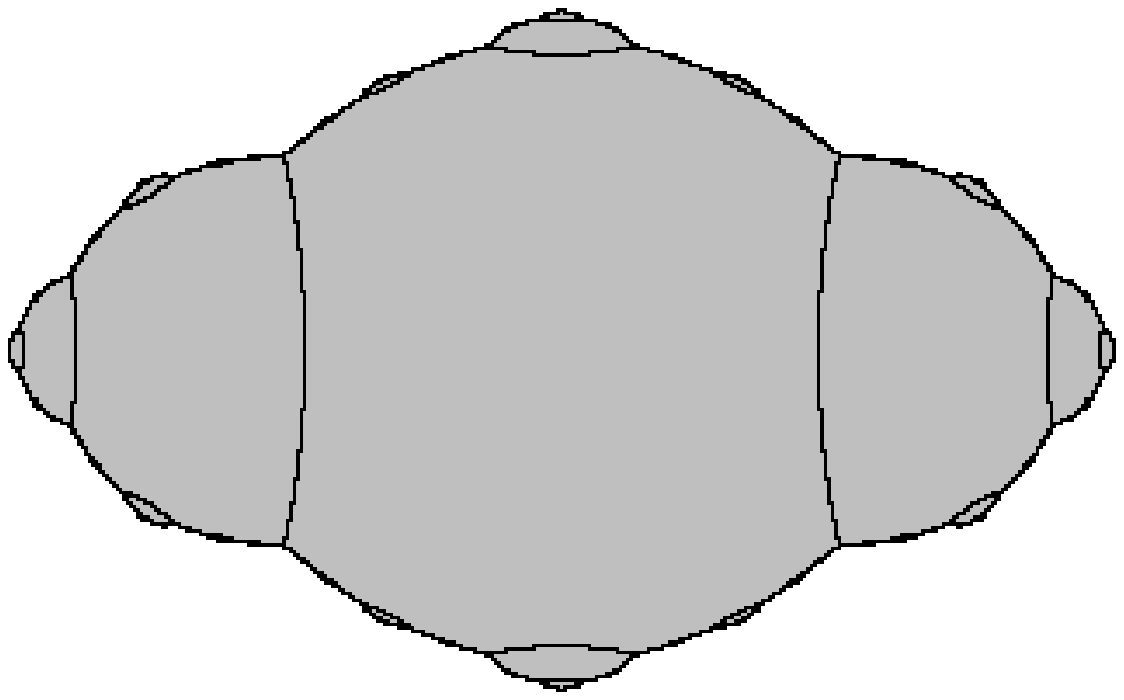,width=.45\hsize}
	\vfil}
    \hfil
  }
  \shortcaption{The filled Julia sets $K(0.75,-0.78)$ and $K(1.5,-0.8)$}
}\endinsert

\theorem Proposition \EasyJKFacts. \hfil
\item{\smbull} $K(\alpha,c)$ and $J(\alpha,c)$ are closed and completely
	   	invariant. 
\item{\smbull} $K(\alpha,c)$ and $J(\alpha,c)$ are connected if and only if
		$0 \in K(\alpha,c)$.
\item{$\smbull$} If K$(\alpha, c)$ is connected, then f is conjugate to the
		map $z \mapsto z^2,\; |z| >1$ on the complement of
		$K(\alpha, c)$.

\proof The proof is essentially the same as the proof for quadratic
polynomials. Refer to \cite{DH}, \cite{Bl}, or \cite{M}.\QED

\theorem Proposition \KPathSimplyConnected.
Every path component of $K(\alpha, c)$ is simply connected.

\proof
If $\gamma$ is a Jordan curve contained in $K(\alpha, c)$, then its iterates
are bounded. Consider the component $D$ of ${\C} - \gamma$ which does not
contain infinity. Since $f_{\alpha, c}\!\!: {\C} \maps {\C}$ is an open map, a
point in $D$ cannot map to a boundary point of $f^{n}(D)$ under
$f^{n}$. Thus $\boundary f^{n}(D) \subset f^{n}(\boundary D)
= f^{n}(\gamma)$, and $f^{n}(D)$ will then be bounded. Therefore
$D$ is contained in $K(\alpha, c)$. 
\QED

Define for fixed $\alpha$ the {\bit connectedness locus} $\ConnLocus$ 
of the family $\{ f_{\alpha, c} \}_{c\in{\C}}$ as 
$$	\ConnLocus = \{ c \st K(\alpha,c) \;\hbox{is connected} \}$$
The set $\MSet$ is known as the Mandelbrot set. An interesting issue is
the dependence of $\ConnLocus$ on the parameter.  An isolated saddle-node
bifurcation  which results in an attractor which attracts the critical
point could ruin the continuity in the Hausdorff topology. We have not
observed such a bifurcation. At this point we formulate the following
conjecture:    

\theorem Conjecture \ConnLocusVariesContinuous.
The connectedness locus $\ConnLocus$ varies continuously with $\alpha$ in
the Hausdorff topology.

There is another set of points $\cal D_\alpha$ in the parameter space that
is interesting,
$$
{\cal D}_\alpha=\{ c \st K(\alpha,c) \;\hbox{is not totally disconnected} \}.
$$
In the conformal case, $K(\alpha,c)$ is not connected if and only if it is
totally disconnected.  In section~2, we show that for large $c$ the set
$K(\alpha,c)$ is totally disconnected.  In the case where $\alpha<1$, there
are $c$ values for which $K(\alpha,c)$ is not connected and not totally
disconnected (see section 4).   It may be that for $\alpha\geq 1$ that 
$\ConnLocus=\cal D_\alpha$.   This is about all we know about $\cal
D_\alpha$.  It would be interesting to find a computer algorithm to draw
this set. 

\medskip
Besides the Mandelbrot set ($\MSet$), the two extreme examples can be
fairly well understood.
 
\theorem Proposition \COneHalf.
The connectedness locus $\CLocus{1/2}$ is a union of half lines,
containing the origin. 

\proof 
Let $f_c$ denote the map $f_{1/2, c}$. Then for any $k > 0$,
$$	f_{k c}(k z) = k f_c(z).$$
Consider the orbit of the critical point. It is easily seen by induction
that $f_{k c}^{n+1}(0) = k f_c^{n+1}(0)$. Therefore, the property that the
orbit of the critical point be bounded is independent of $k$.
\QED

\theorem Proposition \CInfinity.
As $\alpha \rightarrow \infty$, $\ConnLocus$ converges in the Hausdorff
topology to the unit disk.

\proof
Note that for $|c| > 1$ and $\alpha$ large enough, $f_{\alpha, c}^2(0)$ is
close to infinity. Consequently, any Hausdorff limit is contained in the
closed unit disk.  On the other hand, when $|c| < 1$, and $\epsilon$ is
small, then for $\alpha$ large enough, the orbit of the critical point is
contained in the disk of radius  $|c| + \epsilon$. Therefore, any open disk
contained in the closed unit disk is contained in  $\ConnLocus$ for
$\alpha$ large enough.
\QED 

\medskip
In the following sections we describe some theorems, conjectures, and
numerical results.

	\ifSingleSecs \vfil\eject \fi

	\ifx\macrosLoaded\undefined 
        
\fi
\secno=1  

\section{When Disconnected Filled-in Julia Sets are Cantor Sets}

In the holomorphic case ($\alpha$ = 1), filled-in Julia sets are
totally disconnected when they are disconnected. When $\alpha \neq  1$, this
need no longer be true. One can find values of the parameter for which there
are periodic attractors, while the critical point tends to $\infty$. These
examples have only been found when $\alpha < 1$ (see section 4). On the
other hand, when $\alpha$ is fixed and $|c|$ is large enough, this cannot
occur. 

\theorem Theorem \KTotallyDisconn.
Let $\alpha$ be fixed, and $|c|-|2c|^{1/2\alpha}\geq 1$,
then $K(\alpha,c)$ is totally disconnected,
 $K(\alpha,c) = J(\alpha,c)$ and f$_{\alpha,c}|_{J(\alpha,c)}$ is uniformly
expanding.

\proof
The idea of our proof of this theorem is straightforward. First, we show in
the lemma below that there is a disk containing the critical point which
iterates to $\infty$.   The next proposition shows that the map on the
filled-in Julia set is uniformly expanding; the theorem follows immediately.
\QED

\theorem Lemma \ZBounded. 
If $|c| \geq 2^{{{1}\over{2\alpha-1}}}$, then
$$
\(|c|-|2c|^{1/2\alpha}\)^{1/2\alpha} \leq |z| \leq |c|
$$
for any  $z\in K(\alpha,c)$.

\proof
When $|c| \geq 2^{{{1}\over{2\alpha-1}}}$, we have $|c|^{2\alpha} \geq 2 |c|$.
Consider a point $z$ with $|z| > |c|$.
Then
$$
\eqalign{
	|f_{\alpha,c}(z)| &\geq |z|^{2\alpha} - |c|			\cr
			  &= (|z/c|^{2\alpha})|c|^{2\alpha} - |c|	\cr
			  &>  |z/c| (2 |c|) - |c|			\cr
			  &=  2|z| - |c| 				\cr
			  &>   |z|.
}$$
If the orbit of $z$ remains bounded, then by continuity of $f$
there is a limit point $z_\infty$ of the orbit such that
$|f_{\alpha,c}(z_\infty)| = |z_\infty|$, yielding a contradiction.
Therefore, the orbit of $z$ goes to infinity, and so $z\not\in K(\alpha,c)$.

On the other hand suppose $|z| < \(|c|-|2c|^{1/2\alpha}\)^{1/2\alpha}$.
We show that the second iterate of $z$ is outside the disk of radius $|c|$,
and hence  by the above argument, the orbit of $z$ goes to infinity.
We have
$$\eqalign{
  |f_{\alpha,c}(f_{\alpha,c}(z))|	
   	&\geq |f_{\alpha,c}(z)|^{2\alpha} - |c|				\cr
	&\geq \left| |c|-|z|^{2\alpha} \right|^{2\alpha} - |c|		\cr
	&>    (|2c|^{1/2\alpha})^{2\alpha} - |c|			\cr
	&=    |c|.
}$$
\vskip-1.7\baselineskip\QED

\theorem Corollary \LargeCImpliesEscape. \hfil
\item{\smbull}  If $|c|>2^{{{1}\over{2\alpha-1}}}$ then $0\not\in K(\alpha,c)$
	    and thus $c\not\in \ConnLocus$.
\item{\smbull} If $c\in \ConnLocus$ and $z\in K(\alpha,c)$, then
	    $|c|\leq 2^{{{1}\over{2\alpha-1}}}$ and thus 
	    $|z|\leq 2^{{{1}\over{2\alpha-1}}}$.

\theorem Proposition \CLargeImpliesExpanding.
If $|c|-|2c|^{1/2\alpha}\geq 1$, then $f_{\alpha,\,c}$ expands the Euclidean
metric on $K_{\alpha,c}$. 

\proof
Let $f$ denote $f_{\alpha,c}$, let $z$ be a point in $K_{\alpha,c}$, let
$A$ denote $D_{z}f$, and let $v$ be a non-zero tangent vector in $T_{z}\C$.
We must show that $\langle{Av,Av}\rangle > \langle{v,v}\rangle$, or
equivalently $\langle{A^*Av,v}\rangle > \langle{v,v}\rangle$.  Since $A^*A$
has an orthonormal basis of eigenvectors with positive eigenvalues,
it suffices to show that the minimum eigenvalue $\lambda_{\rm min}$ of
$A^*A$ is greater than 1.

For general $f$ we have $\lambda_{\rm min}=(|f_z|-|f_{\bar z}|)^2$,  and in
our case $\lambda_{\rm min}=(\alpha + 1 - |\alpha -1|)^2|z|^{4\alpha - 2}$.
By \ZBounded, we have $|z| \geq (|c|-|2c|^{1/2\alpha})^{1/2\alpha} \geq 1$,
since $z\in K_{\alpha,c}$. When $\alpha \geq 1$,  
$$	\lambda_{\rm min} = 4 |z|^{4 \alpha - 2}\geq 4,$$
and when $1/2<\alpha \leq 1$,
$$	\lambda_{\rm min} = 4 \alpha^2 |z|^{4\alpha - 2}
	\geq 4\alpha^2 > 1.
$$ 
\vskip-1.7\baselineskip\QED

	\ifSingleSecs \vfil\eject \fi

	\ifx\macrosLoaded\undefined 
   
   \ifx\LargeCImpliesEscape\undefined\def\LargeCImpliesEscape{Corollary~2.3}\fi
   \ifx\ZBounded\undefined\def\ZBounded{Lemma~2.2}\fi
   \pageno=6
\fi
\secno=2  

\section{Smooth Julia Sets}

In the holomorphic case there are only two smooth Julia sets. When $c = 0$,
the Julia set is the unit circle, and in this case the Julia set is
hyperbolic. When the critical value $c$ is real and is one of the preimages
of a repelling fixed point, then the Julia set is the closed interval between
$-|c|$ and $|c|$.  This value for $c$ is at the ``tip of the Mandelbrot set",
and in this case the dynamics on the Julia set is sub-hyperbolic (i.e.\ the
dynamics on this Julia set is expanding with respect to a metric which is
smoothly equivalent to the Euclidean metric except at a finite number of
points.)

For fixed values of $\alpha$, one finds readily the parameter value for which
the critical value is a preimage of a repelling fixed point. This fixed point
is real and has coordinate $-c$. The fixed point equation is 
	$$ |c|^{2 \alpha} + c = -c $$
Consequently, $ c =  - 2^{1/(2 \alpha - 1)}$. 
Let us denote the corresponding Julia set by $J_\alpha$. Numerical
observations suggest that when $\alpha$ is between 1/2 and 2, $J_\alpha$ is
indeed an interval, and that when $\alpha$ is greater than 2, $J_\alpha$ is
not contained in the real line.

\theorem Theorem \WhenJIsInterval.
When $\alpha$ is between $0.5$ and $1.7$, the Julia set $J_\alpha$ consists
of an interval. The dynamics on $J_\alpha$ is sub-hyperbolic.

The idea of the proof is straightforward: to find a metric which is
contracted by the inverse (branches) of $f = f_{\alpha,c_\alpha}$. Consider
the metric 
$$
  \rho_\alpha(z) |dz|= {{|dz|}\over{{|c^2-z^2|^{{{2\alpha-1}\over{2\alpha}}}}}}
$$
The restriction of this metric to the interval [$-c,c$] was considered by
Jiang \cite{J1}.

\theorem Proposition \FExpandsRho.
When $0.5 \leq \alpha \leq 1.7$, then $f$ expands the metric $\rho_\alpha$ on
the ball of radius  $2^{{{1}\over{2\alpha-1}}}$.

\proof 
We want to show that $f^*(\rho_\alpha) > \rho_\alpha$. We let
$p={{{2\alpha-1}\over{2\alpha}}}$.  We note that  $2 c =- |c|^{2\alpha}$ and
$f(z)=z^{\alpha+1}{\conj z}^{\alpha - 1} + c$.
$$\openup 2\jot\eqalignno{
  f^*(\rho_\alpha)(z) 
	&= {{|f_z dz + f_{\bar z} d{\bar z}|}\over
	    {\biggl| c^2 - (|z|^{2\alpha-2} z^2 + c)^2\biggr|^p}} 	\cr
	&= {{|(\alpha+1){\bar z}^{\alpha-2} z^{\alpha} {\bar z} dz
	      + (\alpha-1){\bar z}^{\alpha -2} z^{\alpha} z d{\bar z}|}\over
	    {\biggl| 2c|z|^{2\alpha-2}z^2+z^4|z|^{4\alpha-4}\biggr|^p}}	\cr
	&= {{|(\alpha+1)dz  + (\alpha-1){{z}\over{\bar z}} d{\bar z}|}\over
	    {|z|^{1-2\alpha}\biggl| |z|^{4\alpha-4}z^4
			      -|c|^{2\alpha}|z|^{2\alpha-2}z^2\biggr|^p}}\cr
	&= {{|(\alpha+1) dz + (\alpha - 1){{z}\over{\bar z}} d{\bar z}|}\over
	    {\biggl| |z|^{2\alpha-2}z^2 -|c|^{2\alpha}\biggr|^p}}
}
$$
We now wish to show that the ``expansion'' ratio $f^*(\rho_\alpha)/\rho_\alpha$
at a point $z$ in the disk of radius $2^{{{1}\over{2\alpha-1}}}$ is bounded
from below by one.  We have
$$
{{f^*(\rho_\alpha)}\over{\rho_\alpha}}=
\left|(\alpha+1) + (\alpha - 1){{z}\over{\bar z}} {{d{\bar z}}\over{dz}}\right|
\left({{|c^2 - z^2|}\over
       {\left| |c|^{2\alpha}-|z|^{2\alpha-2}z^2\right|}}  \right)^p
$$
Now let $ z=2^{{1}\over{2\alpha-1}}x\e^{\i\theta}$.  Since $z$ is in the
closed disk of radius $2^{{{1}\over{2\alpha-1}}}=|c|$, we have 
$0 \leq x \leq 1$. Denote by $\e^{\i\phi}$ the quantity $d{\bar z}/dz$. We
can express the expansion ratio as the product of two terms. The first is 
$$
	|(\alpha+1) + (\alpha - 1)\e^{\i(2\theta+\phi)}|
$$
which is bounded below by $2\alpha$ when $\alpha\leq 1$ and by
2 when $\alpha\geq 1$.  The second term is 
$$
  \({{|2^{{{2}\over{2\alpha-1}}}(1-\e^{2\i\theta}x^2)|}\over
        {|2^{{{2\alpha}\over{2\alpha-1}}}(1-\e^{2\i\theta}x^{2\alpha})|}} \)^p
  = 2^{{{1-\alpha}\over{\alpha}}} \({{|1-\e^{2\i\theta}x^2|}\over
                                {|1-\e^{2\i\theta}x^{2\alpha}|}} \)^p.
$$
Notice that
$$
  {{|1-\e^{2\i\theta}x^2|}\over{|1-\e^{2\i\theta}x^{2\alpha}|}} > \cases{
	1/\alpha	& when $\alpha > 1$,\cr
      \noalign{\vskip 2pt}       
	1/2		& when $\alpha \le 1.$\cr}
$$
Hence for $\alpha>1$, the expansion factor is greater than 
$\displaystyle{
	2^{1/\alpha}(1/\alpha)^{{2\alpha -1}\over{2\alpha}}
}$,
which is a decreasing function of $\alpha$ and is bigger than 1 for all
$\alpha \in [1,1.7]$. 
For $\alpha<1$, the expansion factor is greater than 
$\displaystyle{\alpha\,2^{{3 - 2\alpha}\over{2\alpha}}
}$,
which is also a decreasing function of $\alpha$, and is  greater than 1 when 
$\alpha=1$. 
We conclude that when $1/2 \leq \alpha \leq 1.7$, the ratio
$f^*(\rho_\alpha)/\rho_\alpha$ is uniformly greater than one for all points
$z$ within the disk of radius $2^{{{1}\over{2\alpha-1}}}$.  
\QED

With this fact in hand, we can now prove \WhenJIsInterval.

\proof
By \LargeCImpliesEscape, the filled-in Julia set is contained in the closed
disk $D$ of radius $2^{{{1}\over{2\alpha-1}}}$.  From the proof of
\ZBounded, it follows that the inverses of $f$ map this disk into itself.
Let $S_0$ and $S_1$ be the two components of the inverse image of the disk.

\midinsert{
	\centerline{\psfig{figure=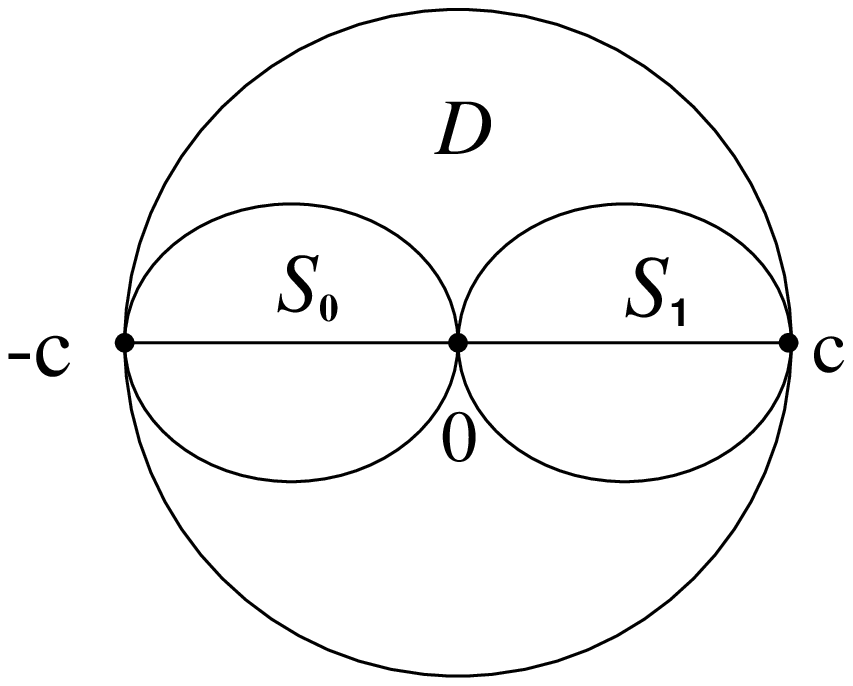,height=1in}}
	\shortcaption{The disk $D$ and its inverse images $S_0$ and $S_1$} 
}\endinsert

The two inverse branches $\psi_i\!:D\to S_i$ are homeomorphisms; by the
previous proposition, they are uniformly contracting.  Thus,  the diameter
of the sets 
$\psi_{\epsilon_1}\circ\psi_{\epsilon_2}\circ\cdots\circ\psi_{\epsilon_n}(D)$
go to $0$ geometrically. Hence, there is exactly one point $x_\epsilon$
with the $n^{\rm th}$ iterate of $x_\epsilon$ in $S_{\epsilon_n}$, where 
$\epsilon=(\epsilon_1,\epsilon_2,\ldots)$ (preimages of the critical value have two
such representations $\epsilon$).  For each $\epsilon$ there is a point on
the real segment with the itinerary $\epsilon$, and so there are no other
points in the filled-in Julia set. 
\QED

\theorem Conjecture \WhenJisIntervalConject.
For all $\alpha$ in (1/2, 2), the Julia set $J_\alpha$ is the interval
$\left[-2^{{{1}\over{2\alpha-1}}},2^{{{1}\over{2\alpha-1}}} \right]$, and
the dynamics on $J_\alpha$ is sub-hyperbolic.

\subsection{Structural Stability}
We now investigate structurally stable properties for $\alpha$ fixed and
$c$ near zero. Consider $f_c = f_{\alpha, c}$. Take $c$ to be zero. The unit
circle $S^1$ is smooth, $f_0$-invariant and repelling. In fact $T_{S^1}\C$
splits as a direct sum of invariant bundles $N$, the radial direction, and
$TS^1$: 
	$$T_{S^1}\C = TS^1 \oplus N.$$
When $v \in TS^1$, $|Df_0(v)| = 2 |v|$, while when $v \in N$, $|Df_0(v)| =
2\alpha |v|$. Let  $m = {{\ln(2 \alpha)}\over{\ln 2}}$, then the dynamics
near $S^1$ is $m$-normally hyperbolic in the sense of  Hirsch-Pugh-Shub
\cite{HPS}.

On $\C - \{0\}$, we have the foliation by concentric circles and the
foliation by radial lines. These foliations are invariant (i.e.\ every
component of $f_0^{-1}$ of a leaf is contained in a leaf), smooth, and
intersect transversly. We consider the stability properties of these
foliations. 

\definition
Let A be an annulus. A foliation on $A$ is {\bit circular} if each of the
boundary components of $A$ are leaves and if every leaf is homeomorphic to a
circle. A foliation on $A$ is {\bit transverse} if every leaf is homeomorphic
to a closed interval and intersects each of the boundary components of $A$ in
a single point. We say that a circular (transverse) foliation on $A$ is $C^k$
when each leaf is $C^k$-diffeomorphic to a round circle (interval) and when
nearby leaves are $C^k$-close. 

Let $A$ and $B$ be domains in the plane. Let $f:B \rightarrow A$ be a
smooth nonsingular map. Then any foliation on $A$ lifts to a foliation on
$B$. We say that a $(C^k)$ foliation on $A$ is {\bit compatible with the
dynamics} if it and its lift to $B$ form a $(C^k)$ foliation of $A \union B$.

Consider a concentric annulus $A$ containing $S^1$. Then $f_0^{-1}(A)$ is
strictly contained in $A$. Choose a circular foliation on $A_0=A-f_0^{-1}(A)$
which is $C^k$-close to the foliation by round circles (in particular,
transverse to the radial foliation).  We can obtain a foliation on $A - S^1$
by repeatedly pulling back by $f_0^{-1}$; adding $S^1$ gives an
$f_0$-invariant foliation $A$ whose leaves are Jordan curves.

One easily shows that every leaf of this foliation is a graph of radial
function ($r(\theta), \theta$); these graphs are uniformly $C^k$ for
all $k \le {{\ln(2 \alpha)}\over{\ln 2}}$. The leaves on
$A_n=A-f_0^{-n-1}(A)$ converge to the round circle in the $C^k$ topology. 

\pageinsert {
 \hbox to \hsize{
	\psfig{figure=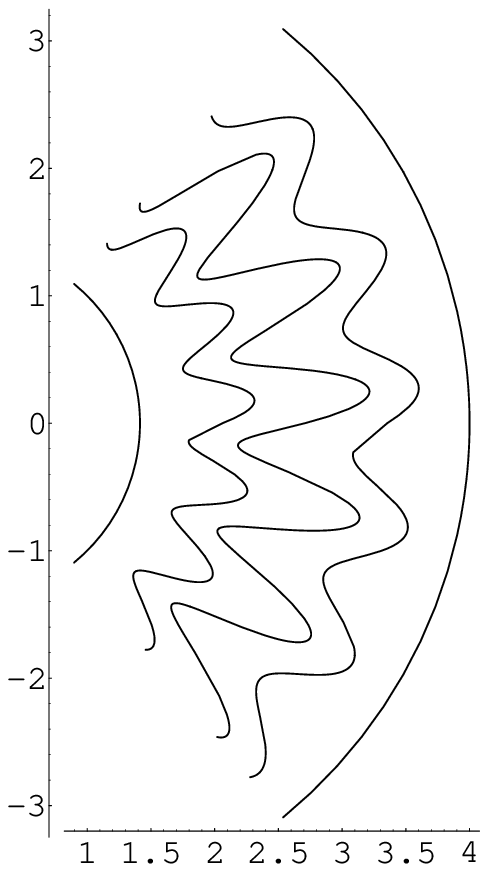,height=1.75in}\hfil
	\psfig{figure=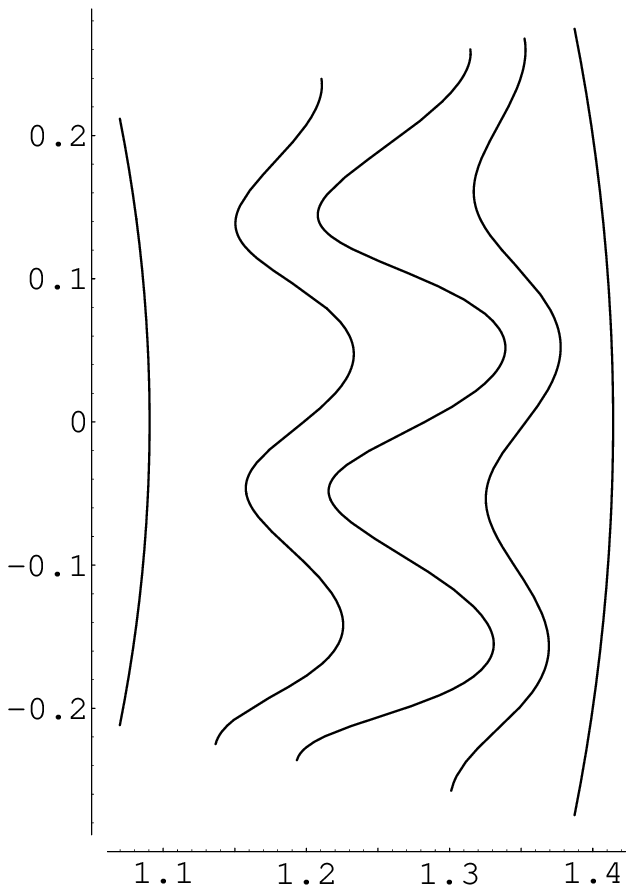,height=1.75in}\hfil
	\psfig{figure=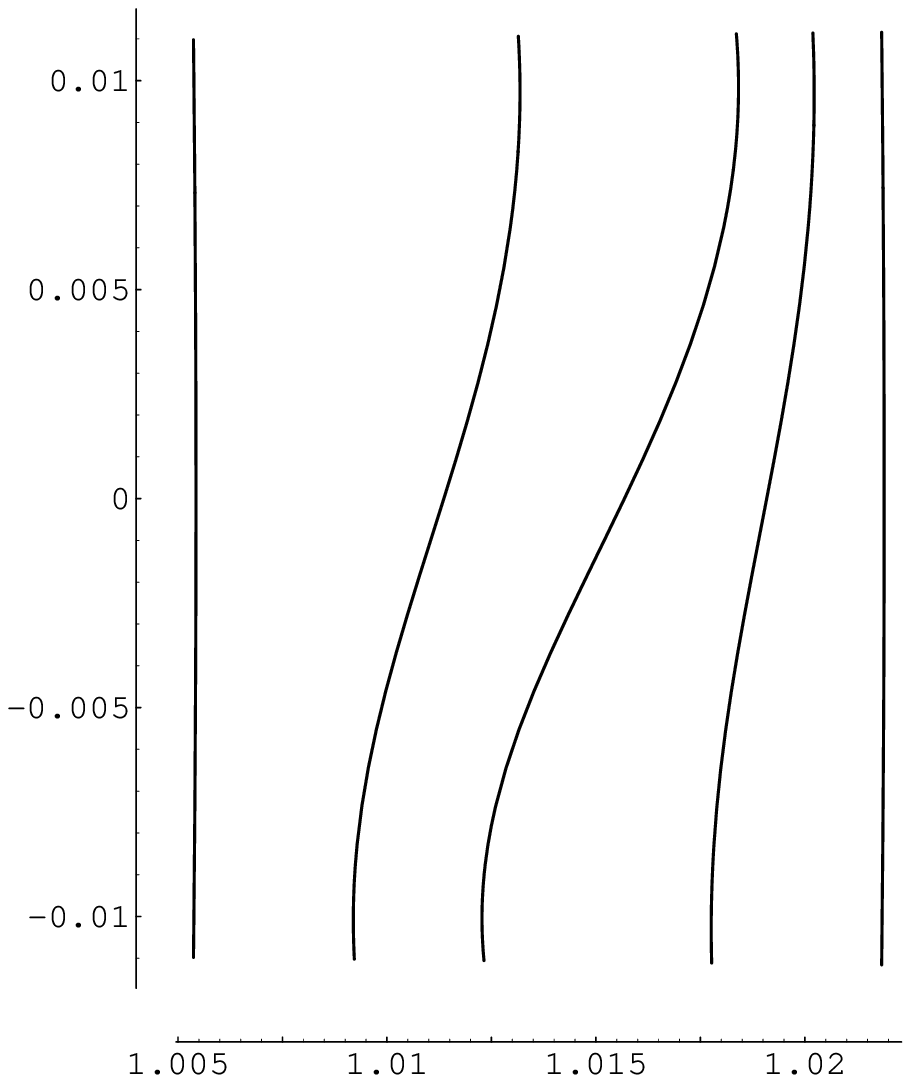,height=1.75in}\hfil
	\psfig{figure=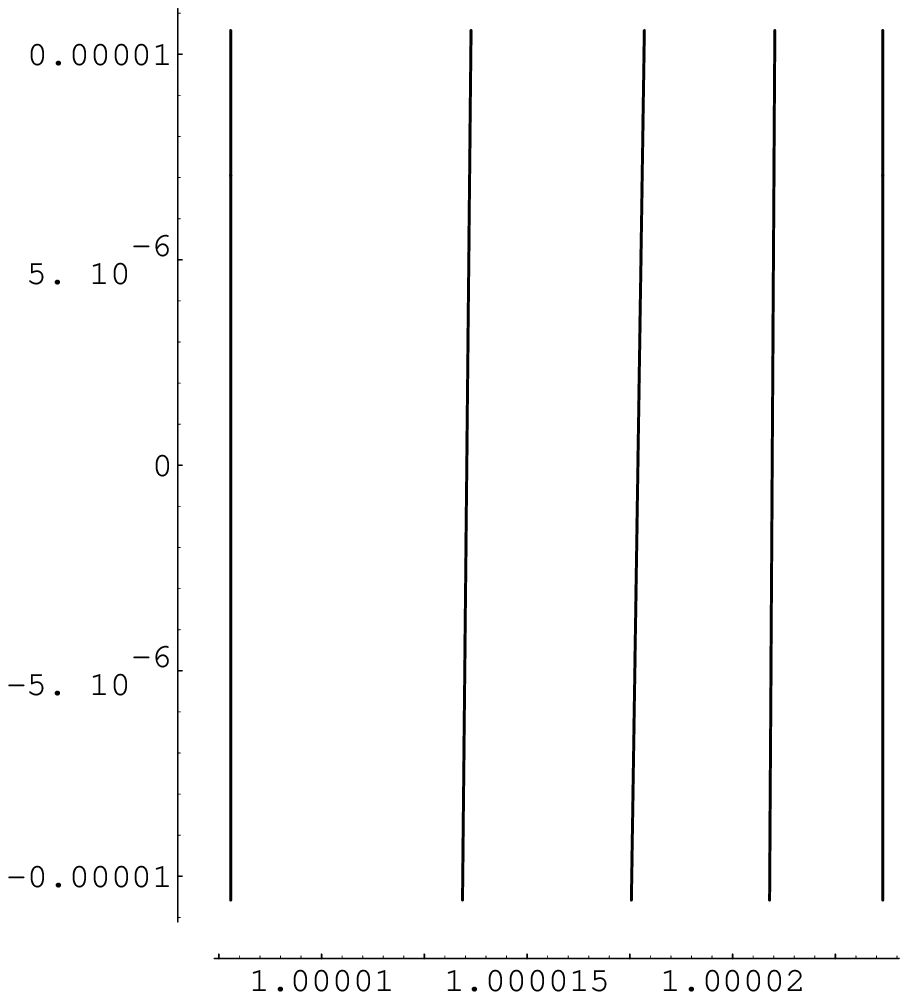,height=1.75in}
 }
 \longcaption{
     Part of several leaves of a circular foliation on the
     outer component of $A_0$, and the pullbacks to $A_1$, $A_3$, and $A_8$
     when $\alpha=2$. These leaves converge to the round circle in the $C^2$
     topology. }

  \vfil
 \hbox to \hsize{
	\psfig{figure=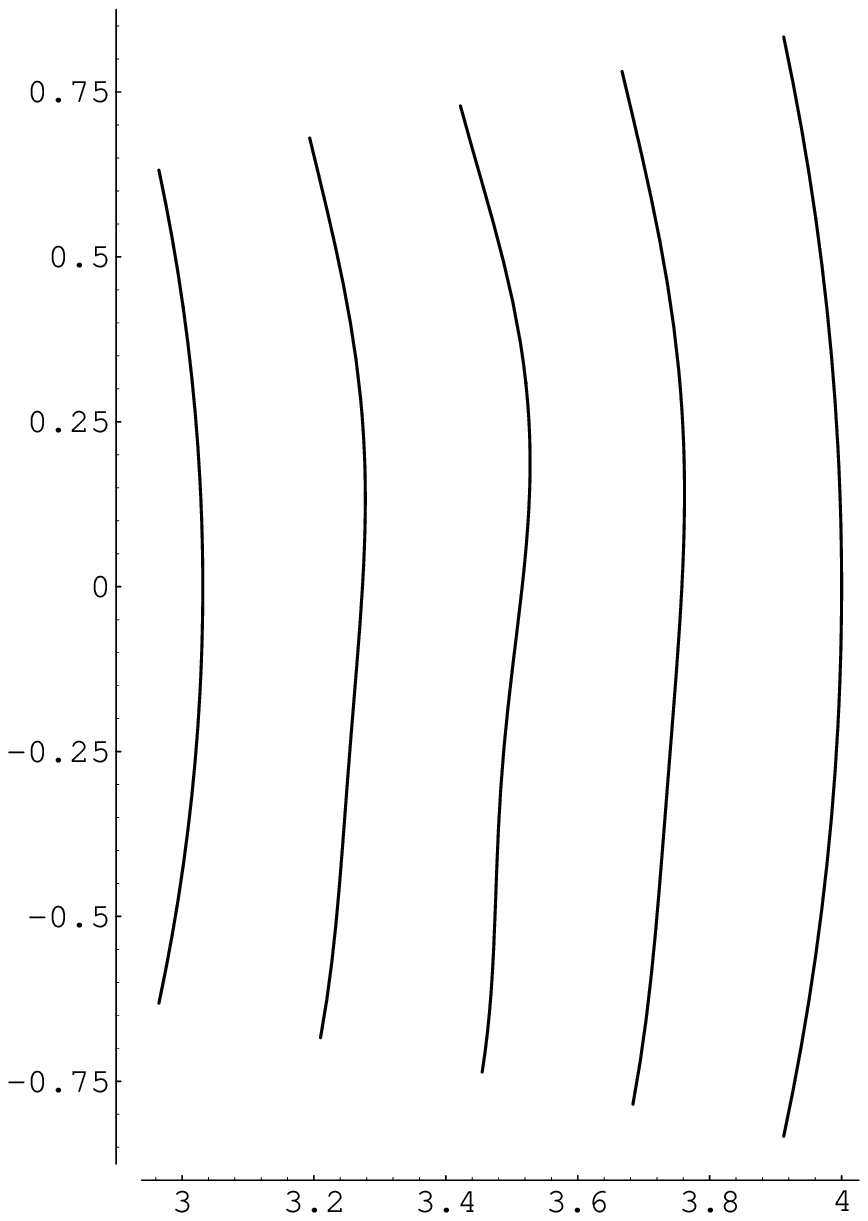,height=1.75in}\hfil
	\psfig{figure=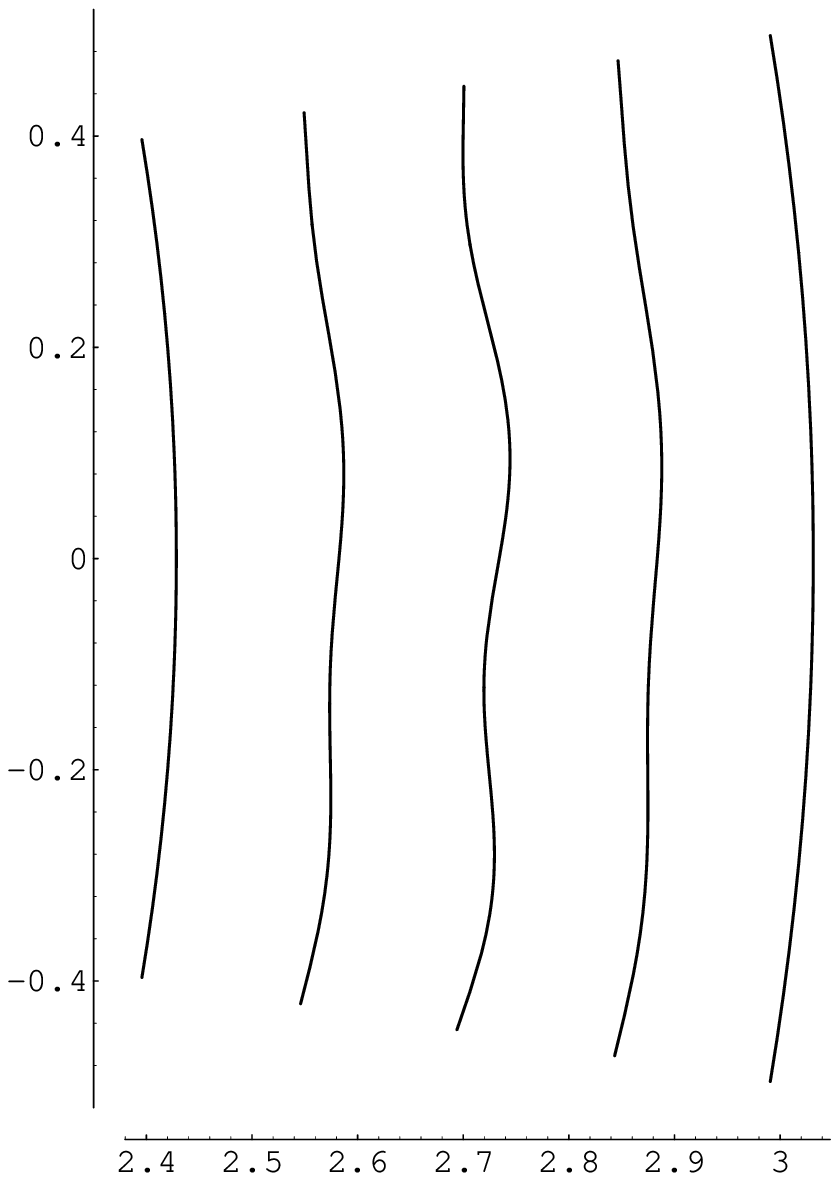,height=1.75in}\hfil
	\psfig{figure=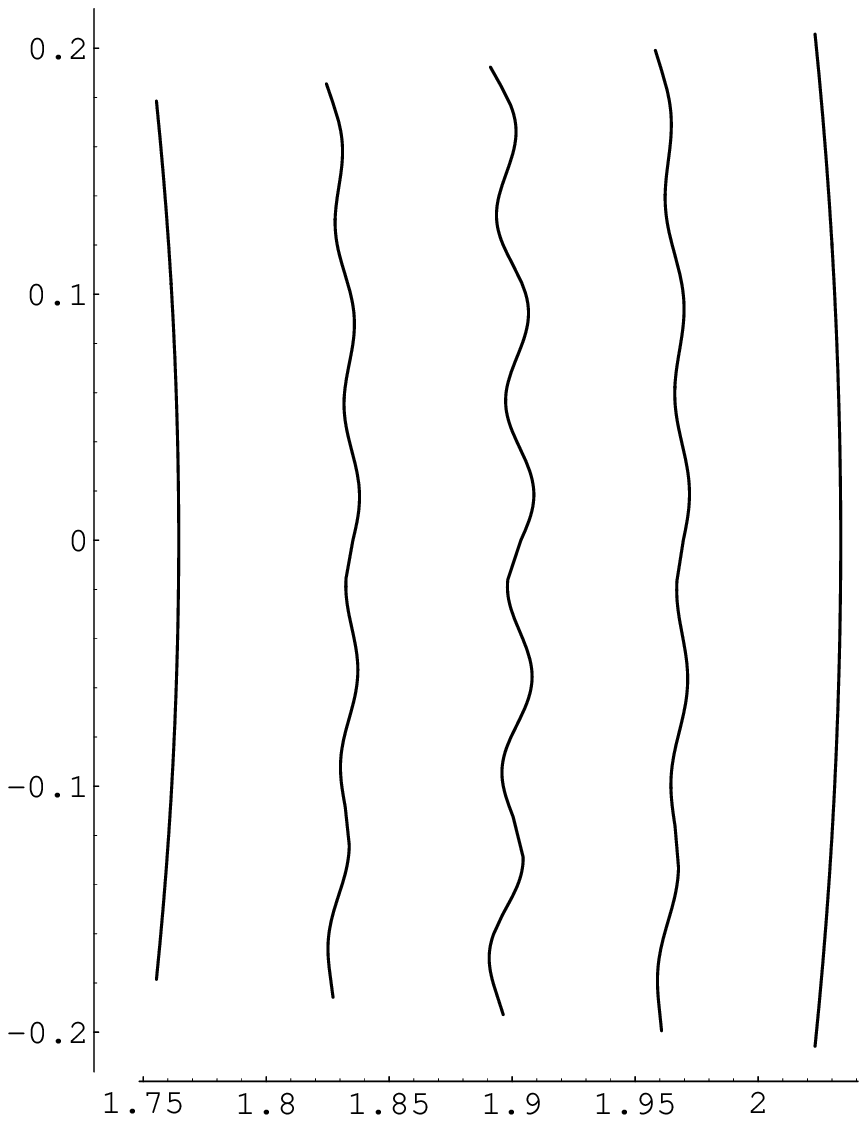,height=1.75in}\hfil
	\psfig{figure=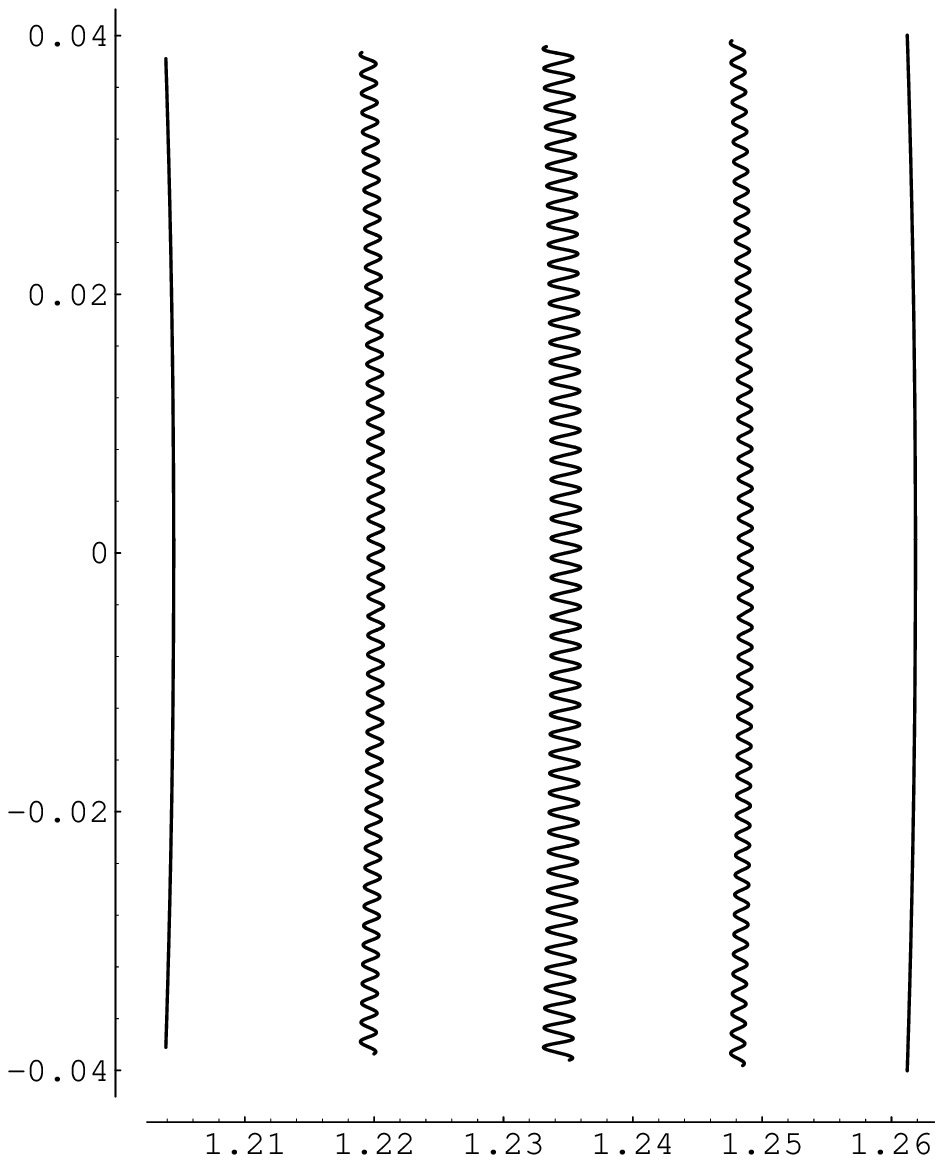,height=1.75in}
 }
 \longcaption{
     Leaves of a circular foliation on the outer component of $A_0$, and
     pullbacks to $A_1$, $A_3$, and $A_8$ when $\alpha=5/8$.  The
     convergence to the circle is only $C^0$.} 
\vfil
 {\vsize=.66\hsize
 \hbox to \hsize{
  \hbox to .46\hsize{\vbox to \vsize{\hsize=.46\hsize
        \vfil
	\psfig{figure=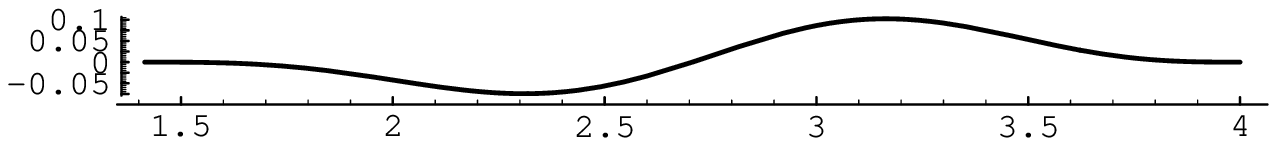,width=\hsize}
	\vfil
        \hbox to \hsize
           {\psfig{figure=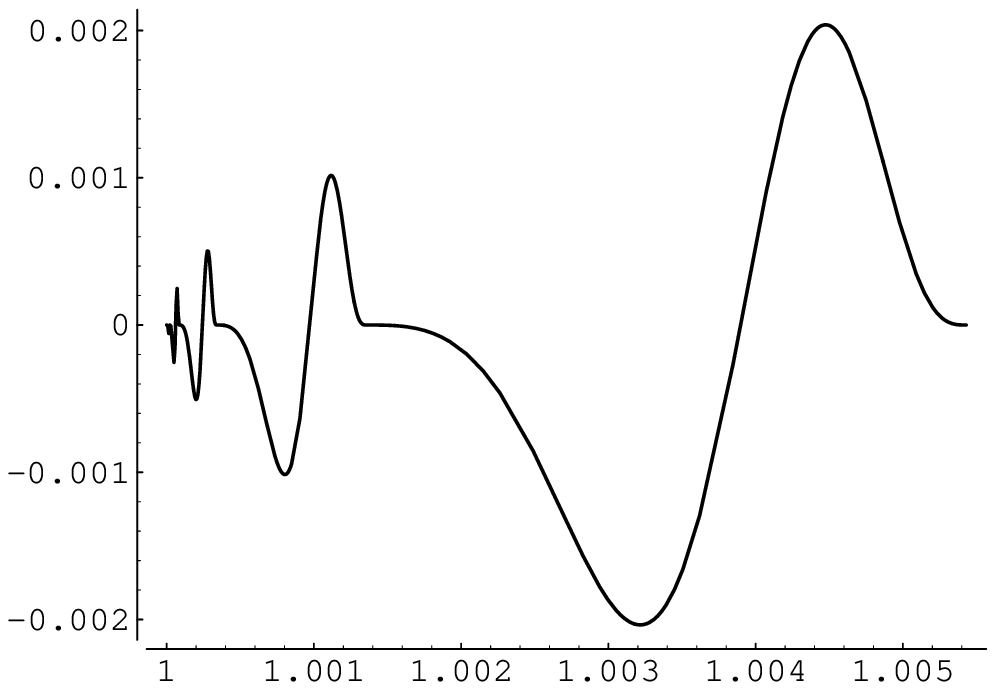,width=.9\hsize}\hfil}
	\vfil
        \centerline{$\alpha=2$}
        \vfil
	\widecaption{
         Above, we show a leaf of a transverse foliation on the outer
         component of $A_0$ which is close to the radial one. Below is  a
         leaf of the corresponding foliation on the outer component of
         $\union_{n=3}^\infty A_n$ when $\alpha=2$; notice the lack of
         smoothness near the limit at $(1,0)$.}  
  }}
  \hfill
  \hbox to .46\hsize{\vbox to \vsize{\hsize=.46\hsize
        \hbox to \hsize
           {\hfil\psfig{figure=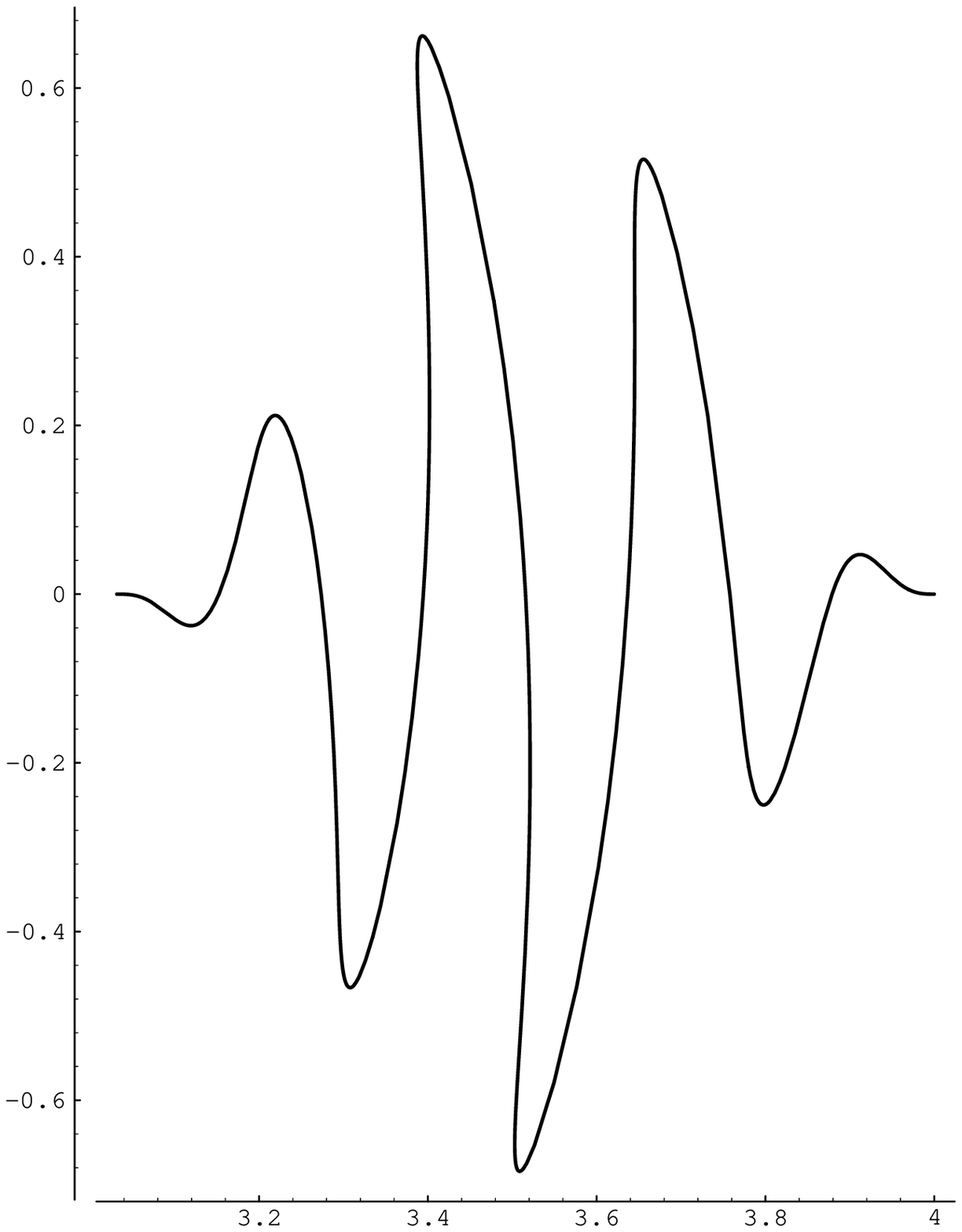,width=.7\hsize}\hfil}
	\vfil
	\psfig{figure=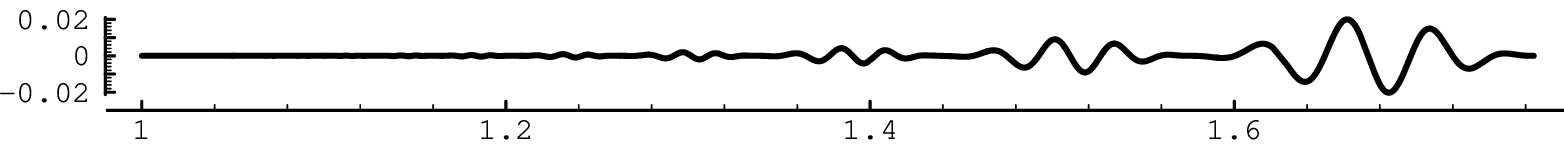,width=\hsize}
	\vfil
        \centerline{$\alpha=5/8$}
        \vfil
	\widecaption{
	As on the left, we have a leaf of the initial foliation above and
        a leaf of the induced foliation below.  Even though the initial
        foliation is quite far from the radial one, as it approaches $S^1$
        the $C^3$ limit of the resulting leaf is exactly radial.} 
  }}
  \hfil
 } 
 } 
}\endinsert

\smallskip
Next, consider a transverse foliation of $A - f_0^{-1}(A)$ by smooth arcs,
running from one boundary component to another in each component annulus,
which is transverse to the foliation by round circles and compatible with the
dynamics.  Pull  back by the dynamics to obtain a foliation of $A - S^1$ by
smooth curves which is transverse to the circular foliation. Since $f_0^{-1}$
is a contraction, each of these curves limits on $S^1$, and at least two
curves land  at each point of $S^1$ (one from the inside and one from the
outside). Moreover, each of these curves is an angular graph
$\(r,\theta(r)\)$ and is uniformly $C^k$ for all
$k\le{{\ln2}\over{\ln(2\alpha)}}$.  
If $\alpha < 1$, the resulting foliation extends to all of $A$ and all
leaves are uniformly $C^k$ for all $k < m$.  However, if $\alpha >1$ and the
initial foliation is not exactly radial, the curves can not meet smoothly at
$S^1$.

\theorem Theorem \StructStability. 
Fix $\alpha \neq 1$ and a concentric annulus $A$ containing the unit
circle in its interior, and let $m={{\ln(2 \alpha)}\over{\ln 2}}$.
\item{\smbull} If $\alpha > 1$, then for all $k < m$ there exists
$\delta_k$ so that when $|c| < \delta_k$, any initial circular $C^k$
foliation on $A - f_c^{-1}(A)$ which is close to the round foliation will
pull back and extend to an $f_c$-invariant $C^k$ foliation on $A$.
\item{\smbull} If $\alpha < 1$, then for all $k < 1/m$ there exists
$\delta_k$ so that when $|c| < \delta_k$, any $C^k$ transverse foliation on
$A - f_c^{-1}(A)$ which is close to the radial foliation and dynamically 
compatible, will pull back and extend to an $f_c$-invariant $C^k$ foliation
on $A$.

\theorem Corollary \WhenJSmooth. \hfil
\item{\smbull} If $\alpha > 1$ and  $|c|<\delta_k$, the Julia set
$J(\alpha,c)$ is a $C^k$ curve.  
\item{\smbull} If $\alpha < 1$ and $|c|<\delta_k$, the Julia set
$J(\alpha,c)$ intersects every leaf of the corresponding radial foliation in
a single point.  

\midinsert{
 \centerline{ \psfig{figure=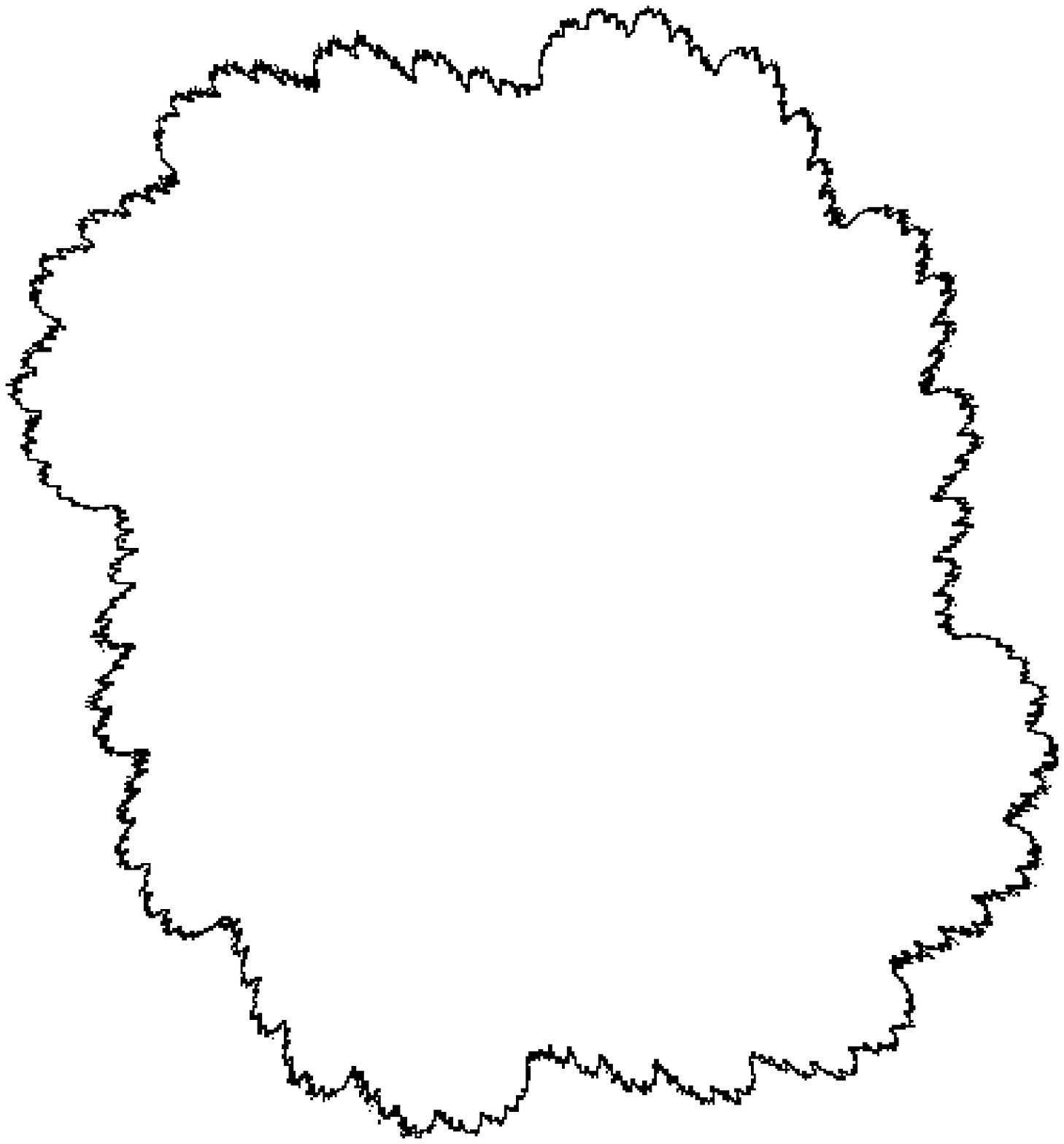,height=1.5in}\hss
  	      \psfig{figure=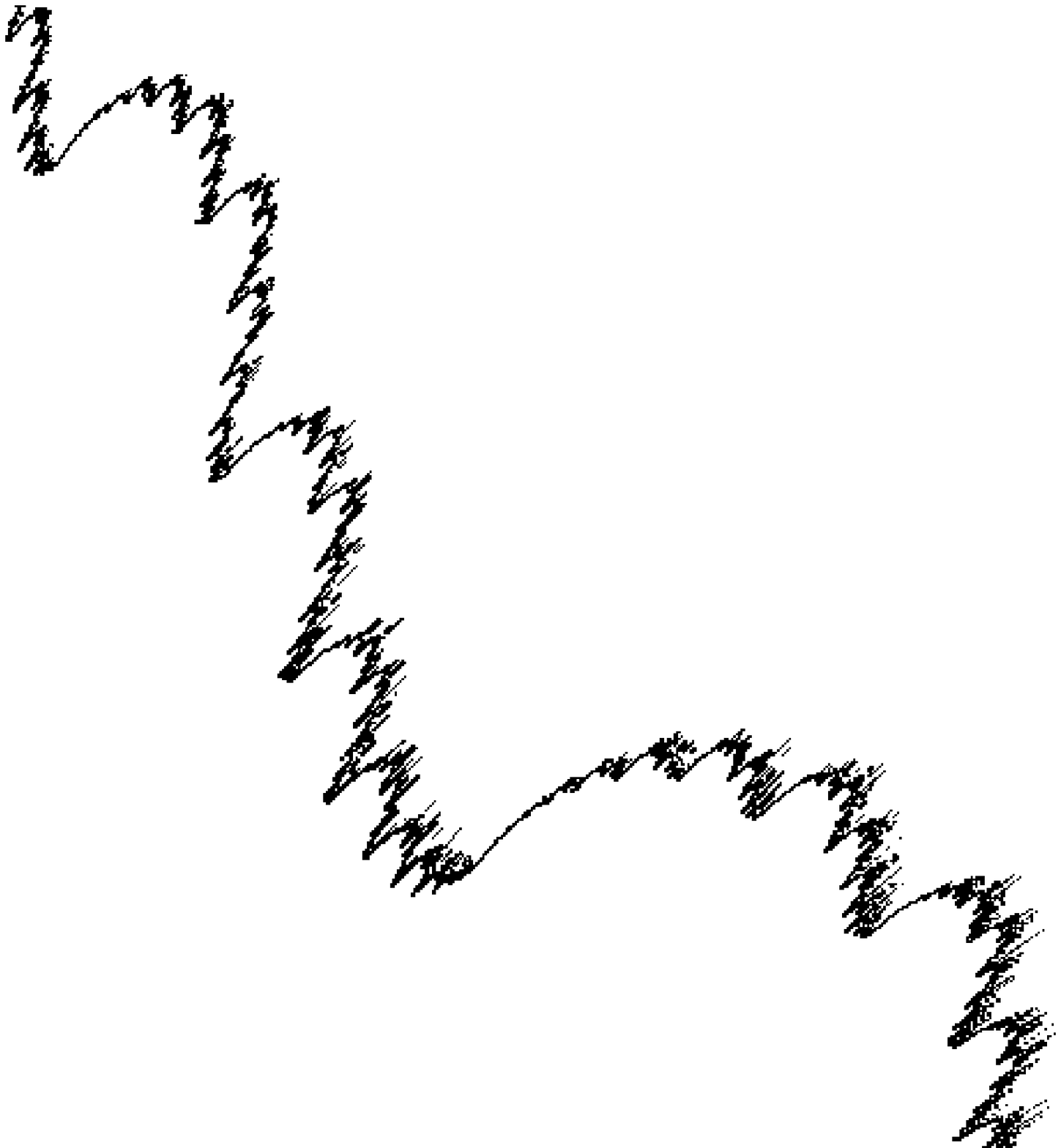,height=1.5in}
 	}
 \longcaption{The Julia set for $\alpha=.75, c=.1+.1i$, and a blowup of it.
  This illustrates \WhenJSmooth: although the Julia set is not at all
  smooth, it looks like the graph of a polar function $r=f(\theta)$.  
}
}\endinsert

\proof
Most of the technical details can be found in \cite{HPS} (diffeomorphisms).
Since the maps we discuss here have degree 2, the initial setup is a little
bit different. Consider the following cone fields on $\C-\{0\}$:
$$\openup1\jot\eqalignno{
 C^+(r, \theta) 
 &= \left\{ R{{\partial}\over{\partial r}} +
	    \Theta{{\partial}\over{\partial\theta}}
         \;\Big|\; |R| \leq |\Theta|  \right\}		\cr
 C^-(r, \theta) 
 &= \left\{ R{{\partial}\over{\partial r}} +
	\Theta{{\partial}\over{\partial\theta}}
	\;\Big|\; |R| \geq |\Theta| \right\}		\cr
}$$
When $\alpha > 1$, $f_0^{-1}$ maps $C^+$ strictly into itself. Consider
the annulus $A$. For $|c|$ small enough, $f_c^{-1}$ maps the cone field $C^+$
on $A$ strictly into itself. Choose an initial circular foliation on
$A_0 = A-f_c^{-1}(A)$ where the tangent vectors at each point are in the
cone $C^+$, and extend this to a foliation ${\cal F}_0$ on all of $A$ which
has the same property. Define a new foliations ${\cal F}_1$ as follows: pull
back ${\cal F}_0$ on $A_0$ by $f_c^{-1}$ to obtain a foliation on $A_1 =
A-f_c^{-2}(A)$. Extend this to all of $A$ as before to obtain ${\cal F}_1$.
We iterate this procedure to obtain a sequence of circular foliations ${\cal
F}_n$ on $A$.  

Now choose $k < m$, and assume that the leaves of ${\cal F}$ are $C^k$. The
techniques in \cite{HPS} show that when $|c|$ is small enough, the sequence
${\cal F}_n$ is $C^k$-compact and therefore has a limit point 
${\cal F}_{\infty}$ which only depends on the choice in $A - f_c^{-1}(A)$.
Since the foliations ${\cal F}_n$ agree on larger and larger domains, 
${\cal F}_{\infty}$ is the only limit point. In particular, the Julia set
$J(\alpha, c)$ is $C^k$.

\smallskip
When $\alpha < 1$, the situation is reversed: the cone field $C^-$ is mapped
into itself by $f_0^{-1}$. When $|c|$ is sufficiently small, $f_c^{-1}$ maps
$C^-$ restricted $A$ into itself. Now consider a transverse foliation on
$A-f_c^{-1}(A)$ which is dynamically compatible and whose tangent-line field
is in the cone field $C^-$. Extend this foliation to all of $A$ so that
the tangent-line field is in the cone field everywhere. We now repeat the
pull-back construction.  When $k < {{1}\over{m}}$ we now have a $C^k$-compact
sequence of transverse foliations on $A$, for $|c|$ sufficiently small. And
again, there is a single limit point which only depends on the initial choice
in $A - f_c^{-1}(A)$.  

We finally argue that the Julia set intersects every leaf in exactly one
point.   The Julia set $J(\alpha, c)$ certainly intersects every leaf in at
least one point. If it intersects in say two points, we can iterate forward 
and conclude that there are points of the Julia set in $A - f_c^{-1}(A)$.
This is a contradiction.  
\QED

Since the construction of the foliations in the \StructStability\ involves
choices, one may wonder whether it is possible to make canonical choices.
This is indeed the case on the unbounded component of the complement of the
Julia set. Construct an invariant foliation near infinity and pull back by
the dynamics. When $\alpha > 1$, one can then make a canonical choice of
$C^k$ circular foliation on the closure of the unbounded component, while
when $\alpha < 1$ one can make a canonical choice of $C^k$ transverse
foliation. An interesting aspect is that this construction is also possible
when $\alpha = 1$, the conformal case. Although \StructStability\ no longer
holds in this case, one still obtains both foliations (radial lines and
equipotential lines), but by quasi-circles, respectively, quasi-arcs (see
\cite{DH}) rather than $C^k$ curves.

\midinsert{
  \centerline{\psfig{figure=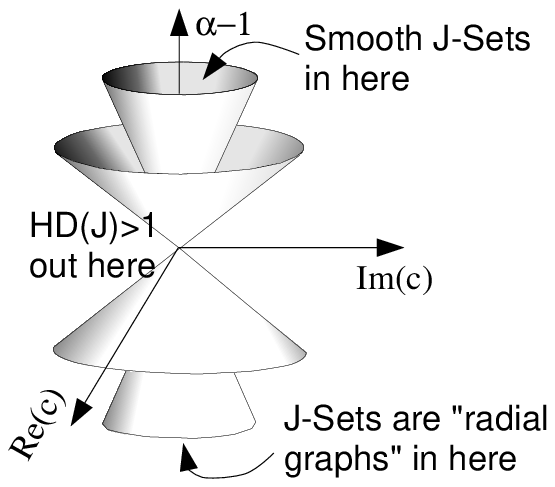,height=2in}}
  \longcaption{A schematic classification of the smoothness of the
        Julia sets near the map $z \mapsto z^2$ in $\alpha$-$c$  parameter
        space. 
	\namefigure{\JiangFig}}
}
\endinsert

\remark
Jiang has shown \cite{J1} that for all $c \ne 0$ with $|c|$ sufficiently
small, there is a $\tau_c > 0$ such that for $\alpha$ with $1-\tau_c \le
\alpha \le 1+\tau_c$, the Julia set $J(\alpha, c)$ has Hausdorff dimension
greater than 1.  See \JiangFig.

	\ifSingleSecs \vfil\eject \fi

	\ifx\macrosLoaded\undefined 
        
	\pageno=12
\fi
\secno=3 
\def\tr{{\rm tr}}

\section{Fixed Points}

In the holomorphic case ($\alpha$ =1) there is one component of the
interior of the connectedness locus which one can understand in all
detail. This is the period one component. For each parameter value in
this component, the corresponding map has a single attracting fixed
point, which moreover attracts the critical point.  This component is
a disk and its boundary is a cardioid.  For every parameter value in
this boundary the corresponding map has a neutral fixed point. When
the eigenvalue of this fixed point is a root of unity, say $\e^{2\pi
\i\,p/q}$, the corresponding parameter value occurs at the
intersection of the closures of two connected components of the
interior of the Mandelbrot set, namely the period one component and a
component where there is a periodic attractor of period $q$.

Part of the key to this picture is the study of the Leau bifurcation
\cite{M}. Here one considers the holomorphic one-parameter family of
holomorphic germs defined near the origin: $$ P_{\lambda}(z) = \lambda
z + z^2 h(z), \quad P_{\lambda}(0) = 0 $$ when $\lambda$ is in the
neighborhood of a root of unity.  This study of the period one
component applies to other hyperbolic components as well. If one
considers a component where one has a periodic attractor of period
$q$, then at each point of the boundary of this component one has a
neutral periodic cycle and, taking the $q^{\rm th}$ iterate, reduces
the study of the bifurcation to that of the Leau bifurcation.  In
particular the boundary of such a component is an algebraic curve.

When $\alpha \neq 1$ our understanding is already incomplete for the
period one component, which we define as the set of parameters $c$ in
the connectedness locus for which $f_{\alpha, c}$ has an attracting
fixed point.  Moreover, the analysis we carry out in the period one
component does not automatically extend to the components
corresponding to periodic attractors of higher period. We show below
that when $\alpha \neq 1$ and an attractor is present, the critical
point is not necessarily attracted to it.

\medskip
Fix $\alpha$. We first analyze the fixed point picture. Note that for
every $z_0$ there is a $c$ such that $z_0$ is a fixed point of
$f_{\alpha, c}$, namely $c=z_0 -
z_0^{\alpha+1}{\conj{z_0}}^{\alpha-1}$.  If $z_0$ is a fixed point of
$f_{\alpha, c}$, the derivative $D(z_0)$ of $f_{\alpha, c}$ at $z_0$ is 
	$$(\alpha + 1) z_0^{\alpha} {\conj z}_0^{\alpha-1}dz +
	  (\alpha - 1) z_0^{\alpha+1} {\conj z_0}^{\alpha-2} d{\conj z}.
$$ The point $z_0$ is an attracting fixed point if the eigenvalues of
$D(z_0)$ are both in the unit disk. In the closure of the set of such
attracting fixed points, there are three important curves:
\item{\smbull} The curve $\delta$ where the determinant of $D(z_0)=1$.
\item{\smbull} The curve $\gamma_+$ where $D(z_0)$ has an eigenvalue of $+1$.
\item{\smbull} The curve $\gamma_-$ where $D(z_0)$ has an eigenvalue of $-1$.

\noindent
A point $z_0=r_0 \e^{\i\theta_0}$ is on $\delta$ if and only if $$
\det D(z_0)= 4 \alpha r_0^{4\alpha - 2} = 1, $$ and therefore $\delta$
is a circle of radius $(4 \alpha)^{1\over{2-4 \alpha}}$.

A point $z_0$ is on $\gamma_+$ if and only if $$ 1 - \tr(D(z_0)) +
\det(D(z_0)) = 0 $$ or equivalently, $$ 1-2(\alpha + 1) r_0^{2\alpha
-1}\cos(\theta_0) + 4\alpha r_0^{4\alpha - 2}=0.  $$ We claim that
$\gamma_+$ is a smooth simple closed curve. This equation has at most two
roots $r_0^{2\alpha - 1}$.  Because we want to only consider the solutions
for $r_0$ positive, we must have $\cos\theta_0 \geq 0$.  Moreover, the
discriminant of this equation is non-negative for $\cos^2(\theta_0) >
4\alpha/(\alpha+1)^2$, that is, in an angular sector about the real axis
(when $\alpha=1$, this sector reduces to a single point) and the discrimant
vanishes at the ends of that angular sector. Consequently, $\gamma_+$ is a
topological circle.  One can check that the curve is $C^1$; we leave that to
the reader.

\smallskip
The curve $\gamma_+$ intersects the curve $\delta$ in two points (one
when $\alpha$ = 1).  Notice that $$\gamma_- = -\gamma_+$$ because
$z_0$ is on $\gamma_-$ if and only if $1 + \tr(D(z_0)) + \det(D(z_0))
= 0$.

\midinsert{
  {\hbox to \hsize{
	\psfig{figure=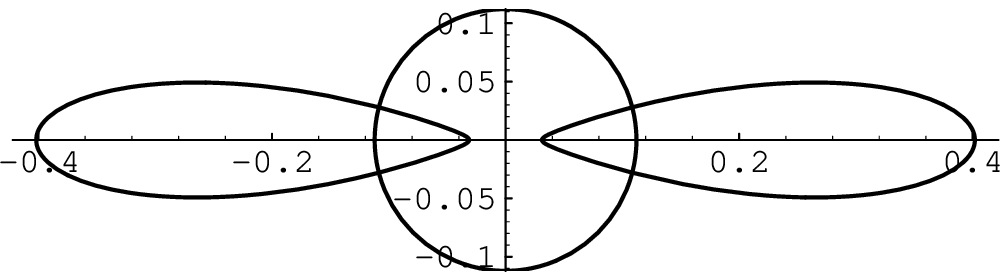,width=.24\hsize}\hfil
	\psfig{figure=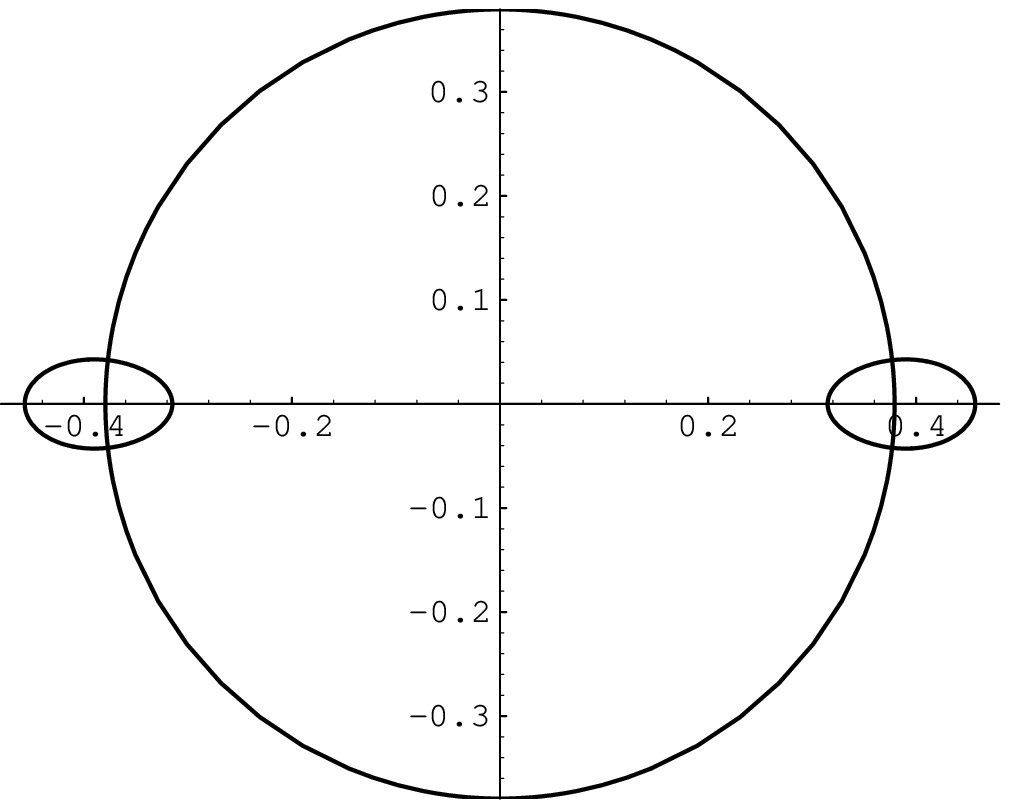,width=.24\hsize}\hfil
	\psfig{figure=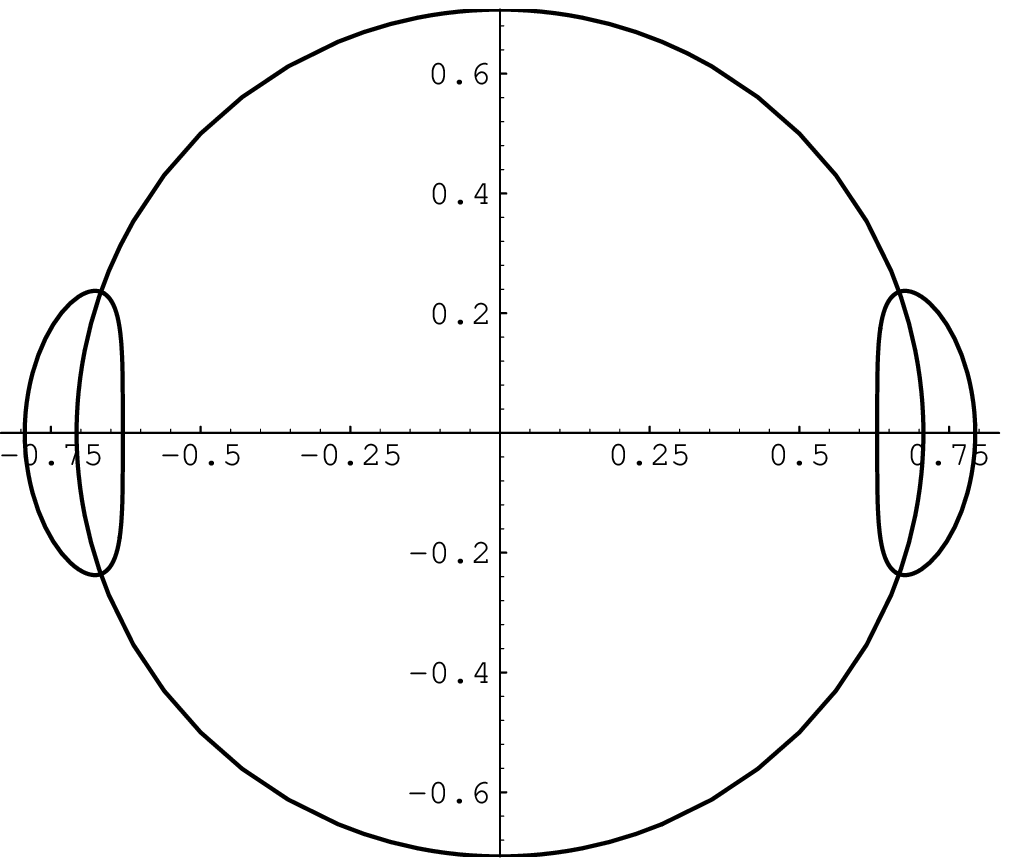,width=.24\hsize}\hfil
	\psfig{figure=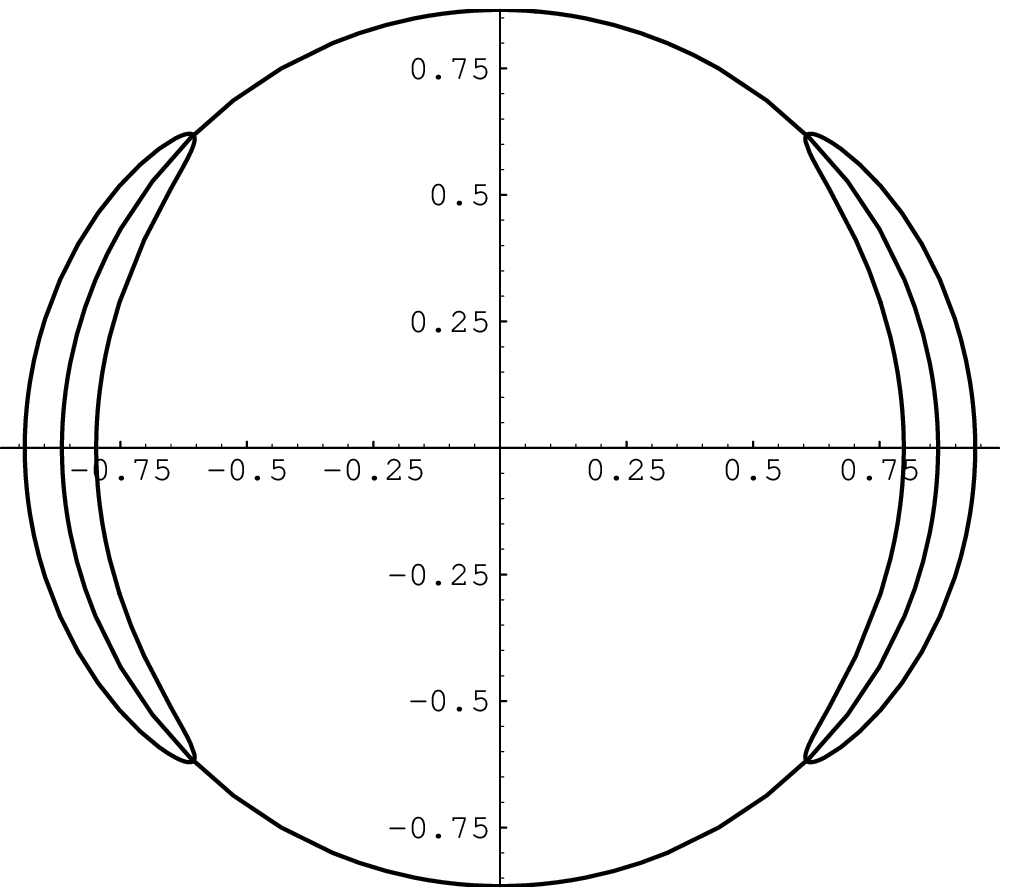,width=.24\hsize} }}
	\shortcaption{The curves $\delta$, $\gamma_+$, and $\gamma_-$
	for $\alpha=.6$, $\alpha=.8$, $\alpha=2$, and $\alpha=6$.}
	}\endinsert

\noindent Define 
$$ P_1(\alpha) = \{ c \st f_{\alpha, c}\ \hbox{has an attracting fixed
point} \}.  $$ The previous analysis immediately provides us with
insight about $P_1(\alpha)$.  Consider the map $p:\C \maps \C$ which
assigns to each $z$ the parameter value $c$ for which $z$ is a fixed
point of $f_c$: $$f_{p(z)}(z) = z$$ One obtains that $$ p(z) = z -
z^{\alpha + 1} {\conj z}^{\alpha - 1}.$$ When $z$ is real, $p(z)$ is
real and $p$ commutes with conjugation. One checks that $p$ is
injective on $\gamma_\pm$, injective on $\delta$ when $\alpha \geq 1$,
and has a single point of multiplicity 2 on $\delta$ when $\alpha <
1$.

\midinsert{
 {\hbox to \hsize{
	\psfig{figure=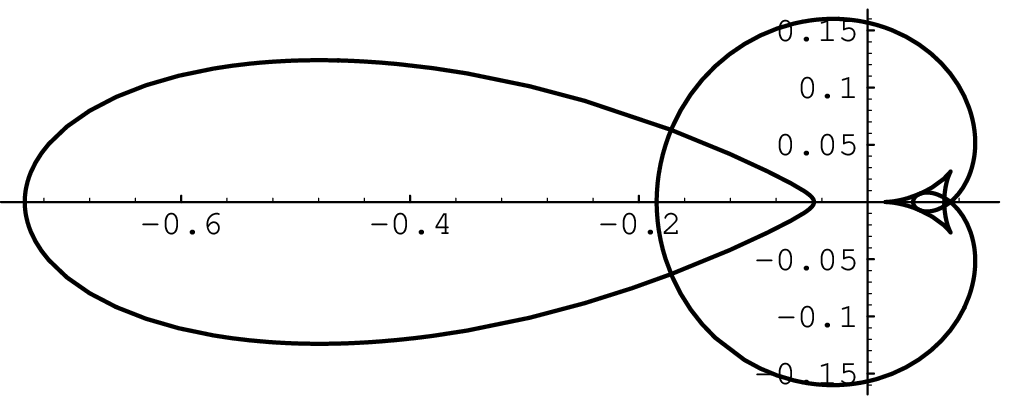,height=.24\hsize}\hfil
	\psfig{figure=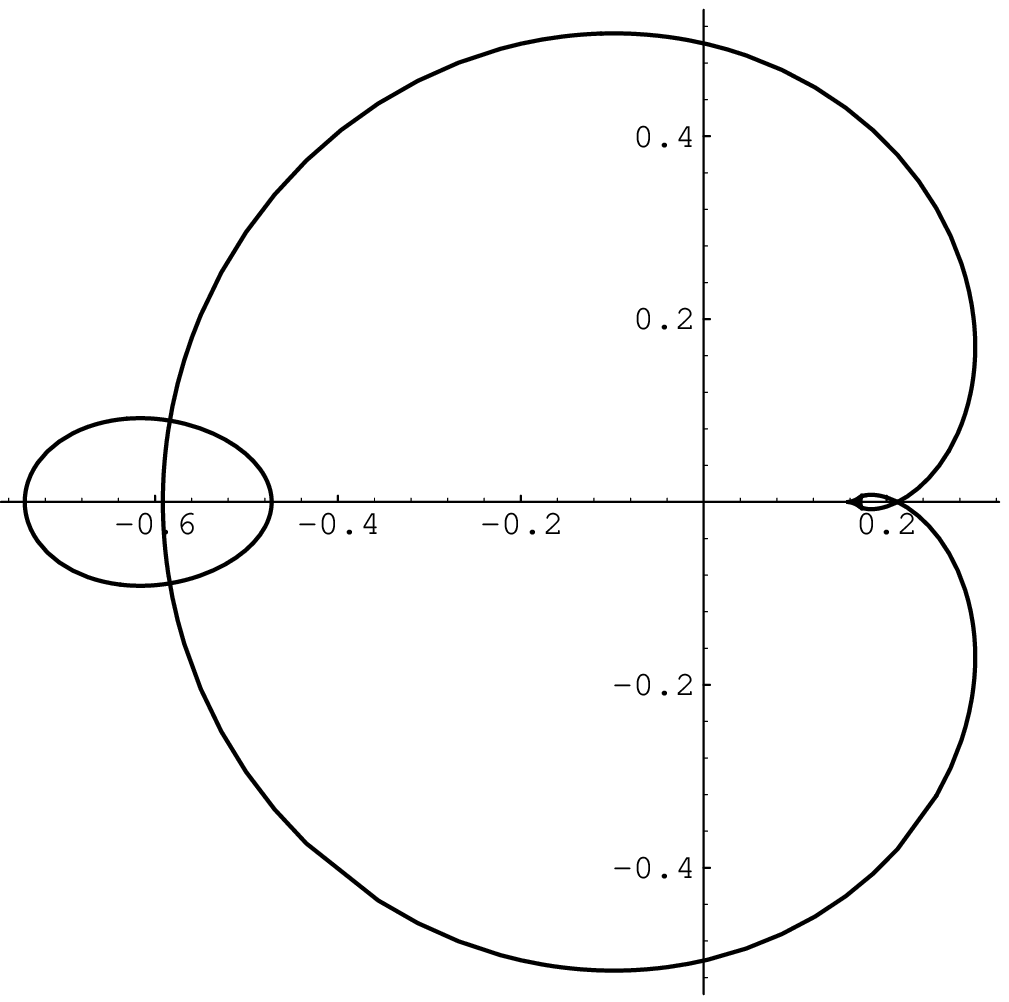,height=.24\hsize}\hfil
	\psfig{figure=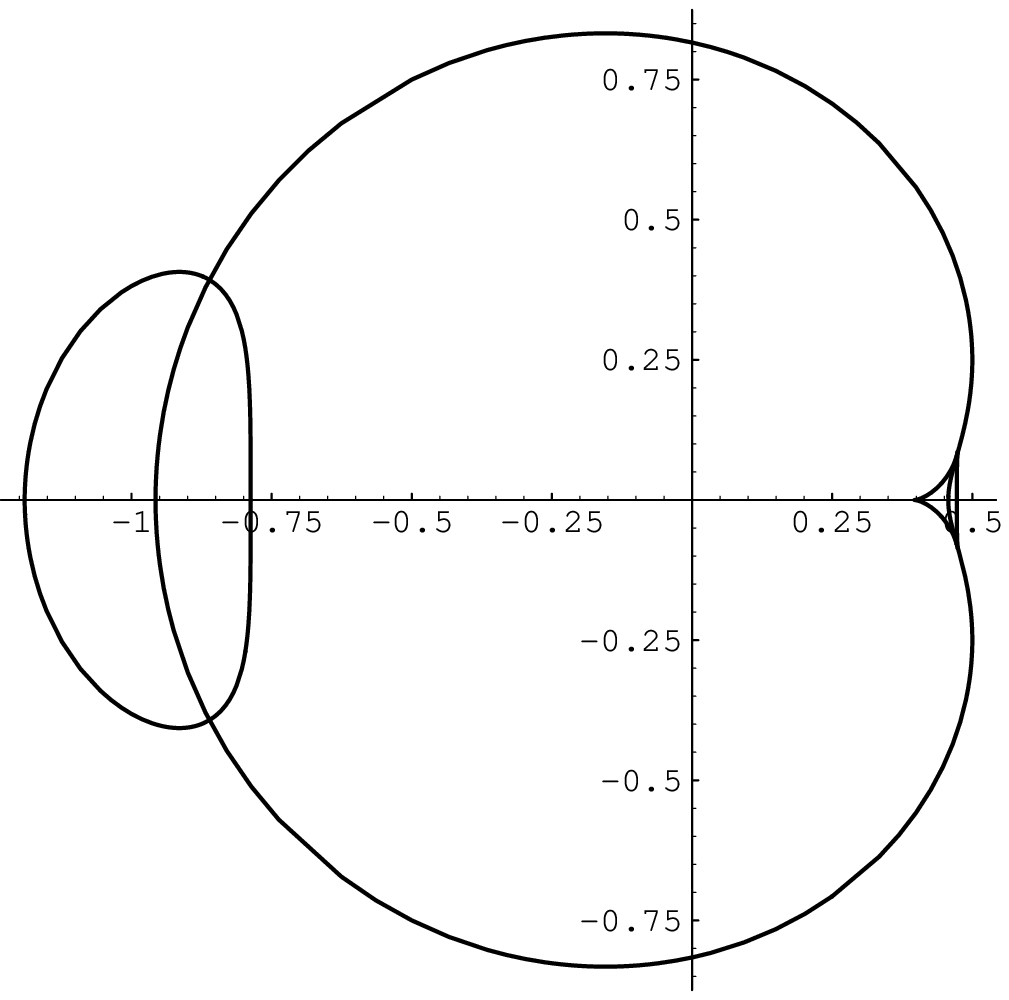,height=.24\hsize}\hfil
	\psfig{figure=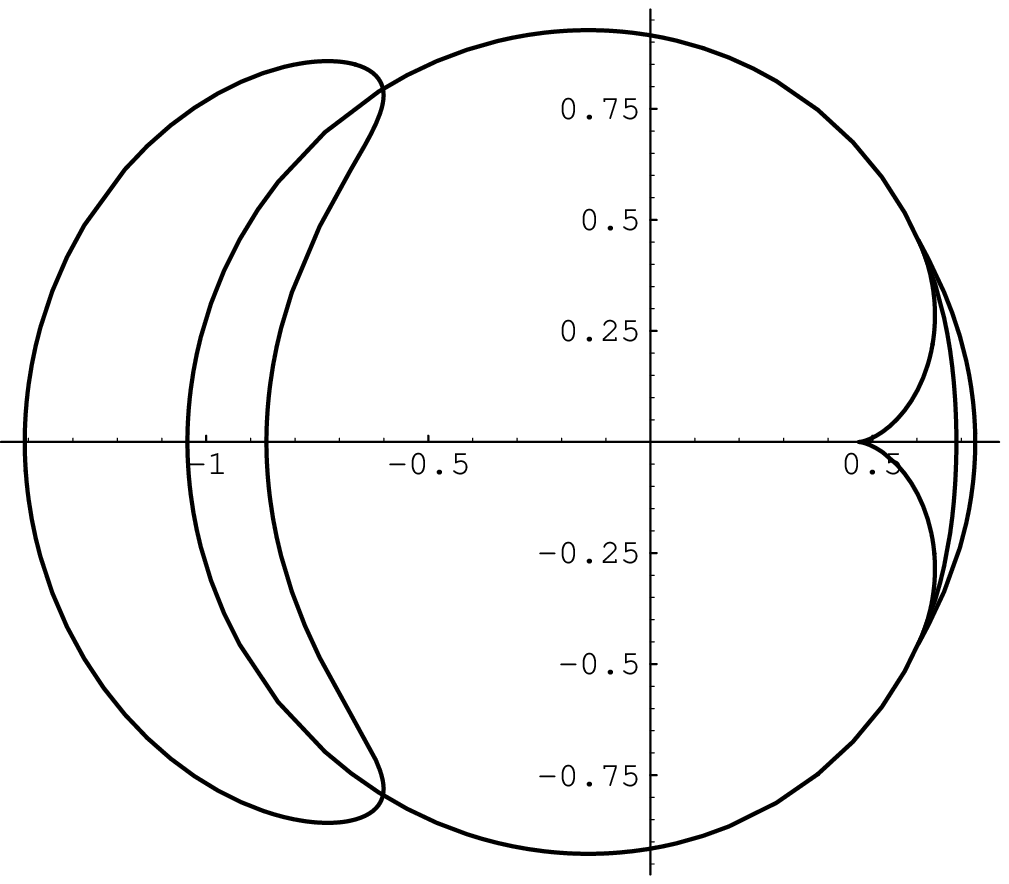,height=.24\hsize} }}
	\shortcaption{The curves $p(\delta)$, $p(\gamma_+)$, and
	     $p(\gamma_-)$ $\alpha=.6$, $\alpha=.8$, $\alpha=2$, and
             $\alpha=6$.\namefigure{\PofCurves} }
  }\endinsert 

\midinsert{
 {\hbox to \hsize{
	\psfig{figure=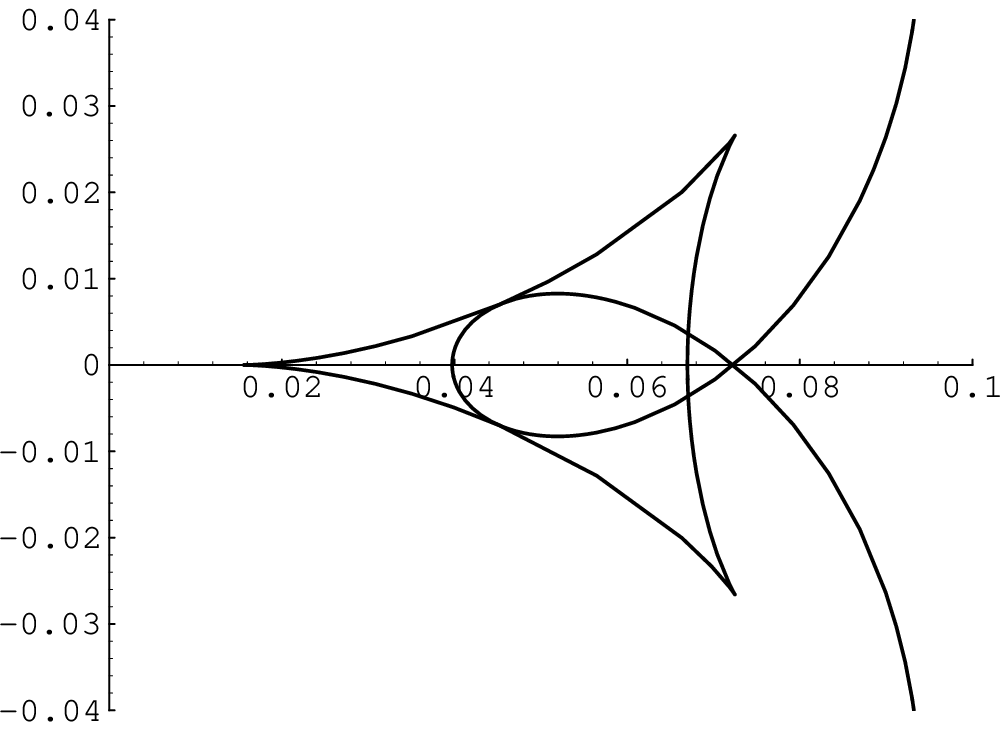,height=.24\hsize}\hfil
	\psfig{figure=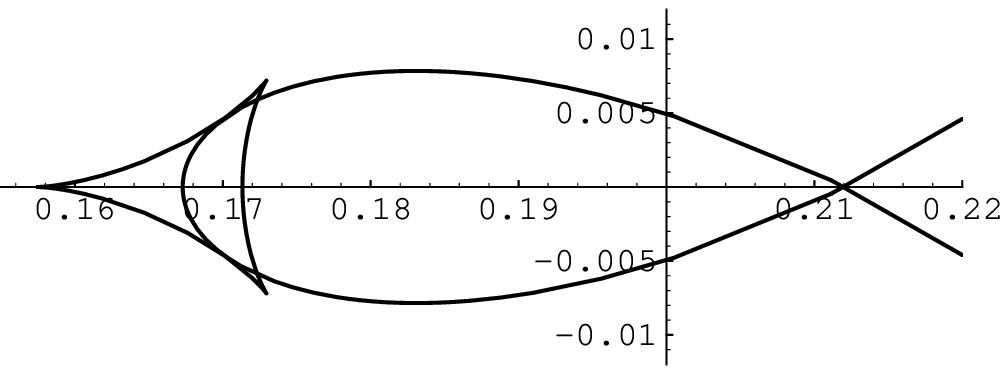,height=.24\hsize}\hfil
	\psfig{figure=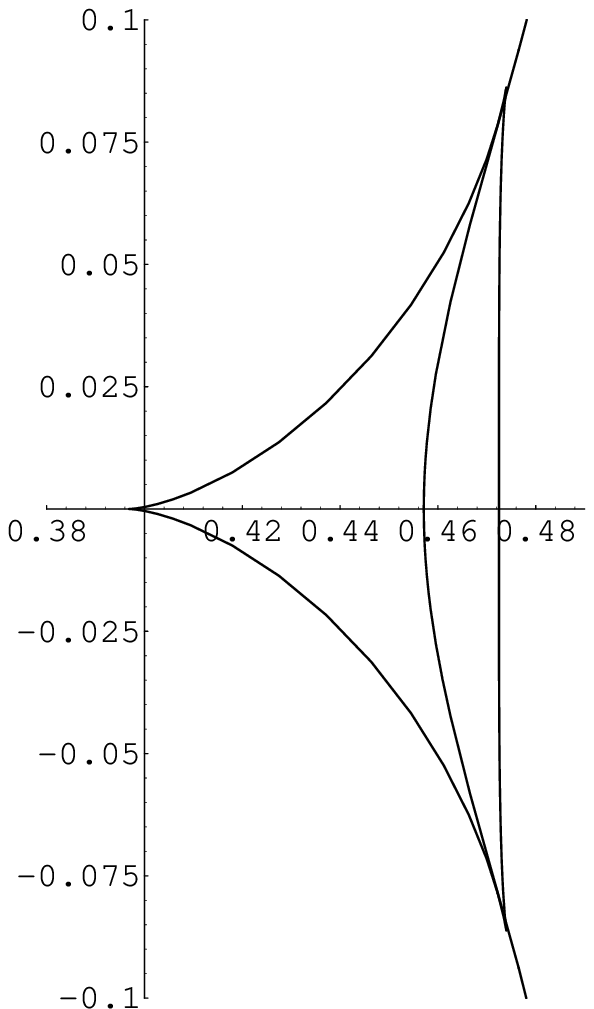,height=.24\hsize}\hfil
	\psfig{figure=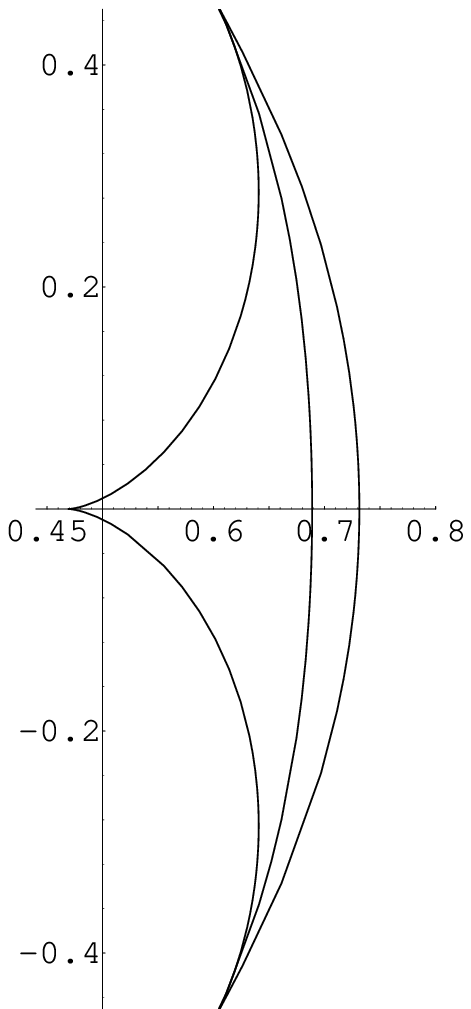,height=.24\hsize} }}
	\shortcaption{A blowup of \PofCurves\ near $p(\gamma_+)$.}
 }\endinsert

\smallskip
We will now discuss the dynamics of the fixed points $z_0\in
p^{-1}(c)$ for $c$ in ${\C}$.  We present the outcome first, followed
by a partial analysis.  The bifurcations occur as you pass through the
curves $p(\gamma_{\pm})$ and $p(\delta)$. One can show that, for
$\alpha\ne 1$, $p(\delta)$ is a lima\c{c}on, $p(\gamma_-)$ is
diffeomorphic to a circle, and that $p(\gamma_+)$ is a simple closed
curve with three cusps.


{\parindent=6em

\bigbreak
\midinsert{
  \medskip \centerline{\psfig{figure=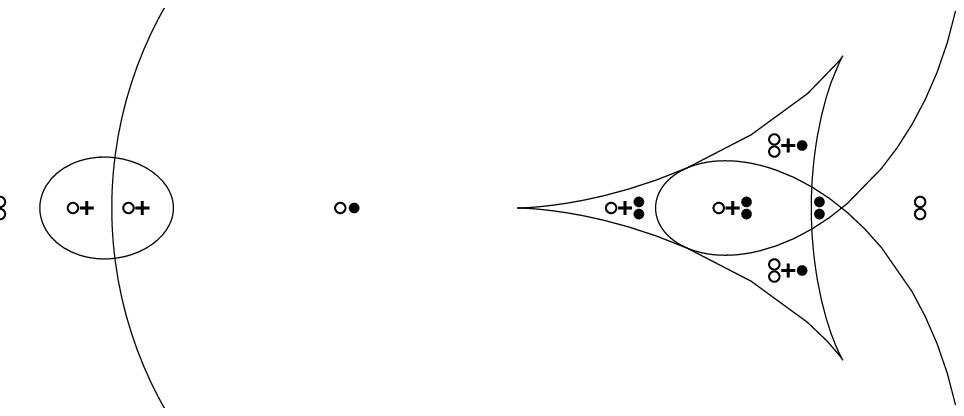,width=.75\hsize}}
  \longcaption{The types of fixed points occuring in the $c$-plane when $1/2
     < \alpha < 1$.  Each region is labelled by symbols indicating the types
     of fixed points which exist there: a symbol $\OneAttr$ is used for each
     attracting fixed point, $\OneRep$ for each repelling fixed point, and
     $\Saddle$ for each saddle. For example, in the region marked
     $\OneRep\Saddle$, there is exactly one repelling fixed point, one
     saddle, and no attracting fixed point.  The upper and lower portions of
     the lima\c{c}on $p(\delta)$ have been omitted, and the region on the right
     side has been greatly magnified.
     \namefigure{\FPtypesLess}}
  }\endinsert

\noindent
{\sc Case $1/2 < \alpha < 1$} (\FPtypesLess):\quad In this case the lima\c{c}on
$p(\delta)$ has an inner loop. In this case, there can be attracting fixed
points which fail to attract the critical point.

\item{Region~$\TwoRep$} Outside the lima\c{c}on~$p(\delta)$ and outside
$p(\gamma_-)$, there are always two repelling fixed points.

\item{Region~$\OneRep\Saddle\TwoAttr$} Inside $p(\gamma_+)$ there are 4
components cut out by the lima\c{c}on.
Region~$\OneRep\Saddle\TwoAttr$ is given by the two pieces which
intersect the real line.  There are two attracting points, one
repelling point, and one saddle.  Part the curve $p(\delta)$ crosses
this region, but crossing this curve only changes the product of the
eigenvalues of the saddle from less than one to greater than one; no
bifurcation occurs.

\item{Region~$\TwoRep\Saddle\OneAttr$} This region consists of the two
components inside $p(\gamma_+)$ which do not intersect the real line;
here there is one attracting fixed point, a saddle, and two repelling
fixed points.  Crossing the curve $p(\delta)$ into
region~$\OneRep\Saddle\TwoAttr$ causes one of the repelling points to
undergo a Hopf bifurcation (see below) and become attracting.  As one
crosses the curve $p(\gamma_+)$ into the $\OneRep\OneAttr$ region, a
saddle and repelling fixed point collide and cancel each other.

\item{Region~$\TwoAttr$} Outside $p(\gamma_+)$ but inside the inner loop of
the lima\c{c}on, we have two attracting fixed points.  When one
crosses $p(\delta)$ into region~$\OneRep\OneAttr$, one of the
attracting points becomes a repelling point, generally with a Hopf
bifurcation.  Entering this region from $\OneRep\Saddle\TwoAttr$
causes the repelling point and the saddle to cancel each other.  We
note that since there are two attracting fixed points, there must be
at least one which doesn't attract the critical point.  In fact, when
$c$ is real the critical point iterates to infinity.

\item{Region~$\OneRep\Saddle$} Inside $p(\gamma_-)$, there is always one
repelling fixed point and one saddle. As above, the part of the
$p(\delta)$ inside this region doesn't cause a bifurcation.  For $c$
real near $p(\gamma_-)$, there is a period two attractor which fails
to attract the critical point; it is attracted to the saddle
instead. As one leaves this region into region~$\TwoRep$, the saddle
splits into a repelling fixed point and a period two saddle.

\item{Region~$\OneRep\OneAttr$}  In this region, which is inside the main
loop of the lima\c{c}on, we have one attractive and one repelling
fixed point. When one enters this region from
region~$\OneRep\Saddle\TwoAttr$, one of the attracting points and the
saddle collide and cancel.  When one enters here from
region~$\OneRep\Saddle$, the saddle point merges with a period~2
attractor and an attracting fixed point is created.  When one crosses
into $\TwoRep$, the attracting fixed point becomes repelling and
typically a Hopf bifurcation occurs.  We will discuss the direction of
the Hopf bifurcation at the end of this section.

\bigbreak

\midinsert{
 \centerline{\psfig{figure=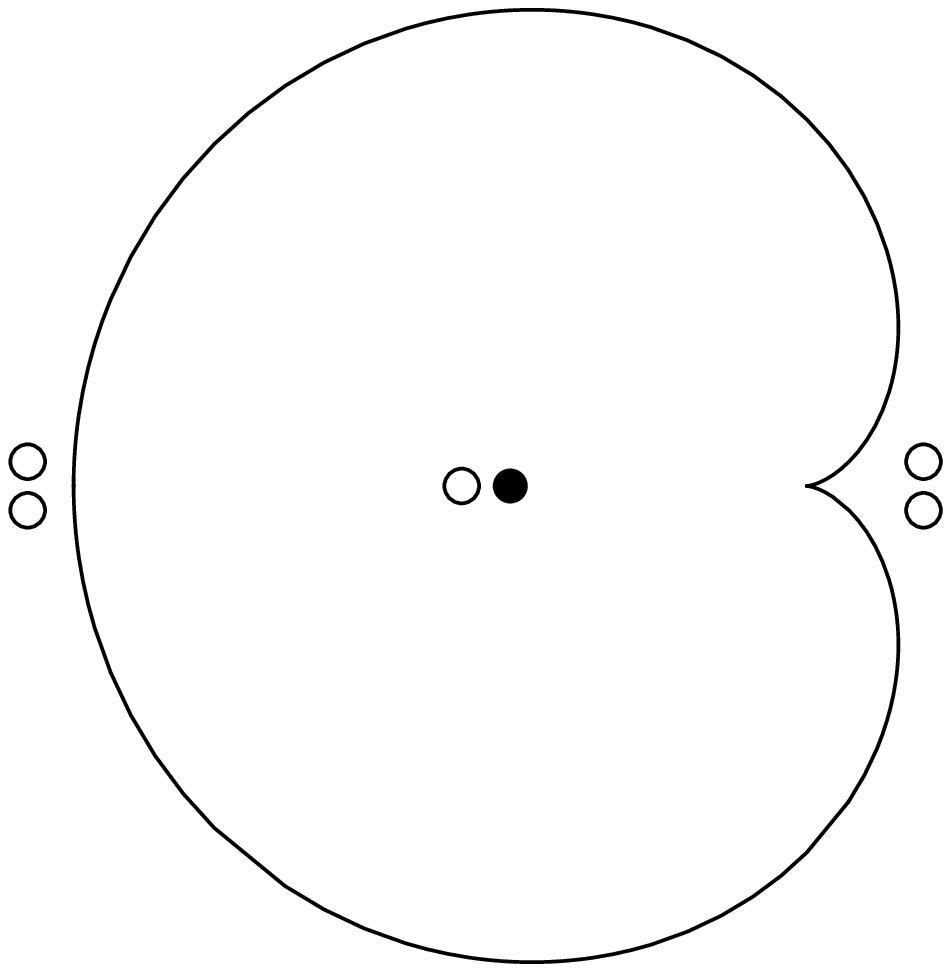,height=1 in}} 
 \longcaption{The types of fixed points occuring in the $c$-plane when
     the map is conformal ($\alpha=1$). The labelling is as in \FPtypesLess. 
 \namefigure{\FPtypesConf}}
  }\endinsert 

\noindent
{\sc Case $\alpha = 1$} (the conformal case, \FPtypesConf):
\quad In this case the lima\c{c}on is a cardioid and $p(\gamma_\pm)$ are
points on the real axis.

\item{Region~$\OneRep\OneAttr$} Inside the cardioid there is one attracting
and one repelling fixed point.  As must occur in a conformal system,
the attracting fixed point always attracts the critical point.

\item{Region~$\TwoRep$} Outside the cardioid there are two repelling fixed
points. In this case, no Hopf bifurcation can occur; when going
through a point on the cardioid for which the derivative at the fixed
point is of the form $\e^{2\pi\i p/q}$, a Leau-Fatou flower
bifurcation occurs. (See
\cite{M}). 
\bigbreak

\midinsert{
 \centerline{\psfig{figure=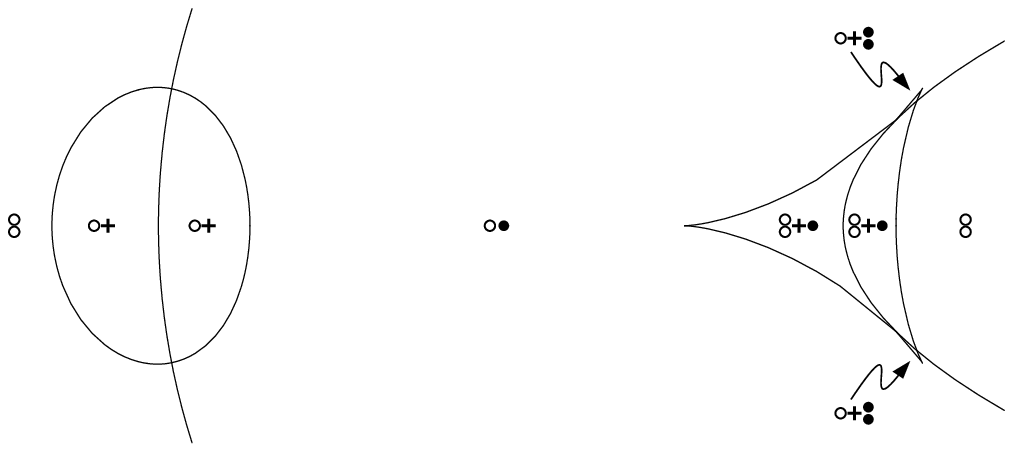,width=.75\hsize}}
 \longcaption{ The types of fixed points occuring in the
     $c$-plane when $\alpha>1$.  As in \FPtypesLess, the upper and lower
     pieces of the lima\c{c}on have been omitted, and the right-hand side
     has been greatly magnified.
     \namefigure{\FPtypesMore}}
  }\endinsert

\noindent
{\sc Case $\alpha > 1$} (\FPtypesMore):\quad In this case the lima\c{c}on is
convex or has a dimple.  In this case we conjecture (based on
numerical evidence) that the critical point is attracted to the
attractive point when it exists.

\item{Region~$\TwoRep$}   Outside the lima\c{c}on and outside
$p(\gamma_\pm)$ there are two repelling fixed points.

\item{Region~$\TwoRep\Saddle\OneAttr$}  The lima\c{c}on cuts $p(\gamma_+)$
into four pieces.  The two pieces intersecting the real line form
region~$\TwoRep\Saddle\OneAttr$, which has two repelling fixed points,
one attracting fixed point, and one saddle point.  As for $\alpha<1$,
crossing $p(\delta)$ inside this region doesn't cause a bifurcation.
When crossing $p(\gamma_+)$ into region~$\TwoRep$, the attracting
fixed point and the saddle collide; crossing into
region~$\OneRep\OneAttr$ causes one of the repelling fixed points and
the saddle to collide.

\item{Region~$\OneRep\Saddle\TwoAttr$}  This very tiny region consists of
the two components inside both $p(\gamma_+)$ and $p(\delta)$ which do
not intersect the real line.  Here there are 2 attracting fixed
points, one repelling, and one saddle point.  When moving from here to
region~$\OneRep\OneAttr$ an attracting fixed point and a saddle cancel
each other.  When moving from here to region~$\TwoRep\Saddle\OneAttr$,
one of the attracting fixed points loses stability and becomes
repelling, generally via a Hopf bifurcation.

\item{Region~$\OneRep\Saddle$} Inside $p(\gamma_-)$, we have one saddle
point and one repeller.  As before, crossing $p(\delta)$ inside this
region causes no bifurcation. When crossing into $\TwoRep$, the saddle
splits into a period two saddle and a repelling fixed point.

\item{Region~$\OneRep\OneAttr$} Inside the lima\c{c}on and outside
$p(\gamma_\pm)$, there is one attracting and one repelling fixed
point.  When one crosses $p(\gamma_+)$ from
region~$\TwoRep\Saddle\OneAttr$, one of the repelling fixed points and
the saddle collide.  When one crosses from region~$\OneRep\Saddle$ to
here an attracting period~2 orbit merges with the saddle to form an
attracting fixed point.

}

\bigskip
We now give the analysis that leads to the bifurcation pictures
above.  We first consider the case where the parameter $c$ is real.

\theorem Proposition \FourFixedPts.
If $c$ is real, the map $f_{\alpha,c}$ has at most four fixed points.

\proof
The point $r\e^{\i\theta}$ is a fixed point if $$
r\e^{i\theta}-r^{2\alpha}\e^{2\i\theta} = c. $$ Looking at the
imaginary part, we see that $$ r\sin(\theta) -
2r^{2\alpha}\sin(\theta)\cos(\theta) = 0,$$ and so either
$\sin(\theta)=0$ or $r=2r^{2\alpha}\cos(\theta)$.  In the first case,
there are two real solutions if $c$ is less than $c_0$, where
$f_{\alpha,c_0}(x)$ is tangent to $y=x$.  In the second case, we
obtain $r^{2\alpha}=c$ after substituting into the real part of the
original equation, and thus we get at most one value for $r$.
Substituting this value for $r$ into the original equation gives a
quadratic equation in $\e^{\i\theta}$ which has a solution for $c$
greater than $c_1<c_0$.
\QED

We can explicitly calculate the type of fixed points which occur on
the real line: the only way the type of fixed points (or total number
of fixed points) can change is by crossing one of the curves
$p(\delta)$ or $p(\gamma_{\pm})$.  Assuming these curves intersect as
discussed earlier, the types of fixed point occuring in each region
can be calculated by considering all possible bifurcations which can
occur.  We know that $p(\delta)$ is a lima\c{c}on.  The difficult part
of the analysis is then to figure out how the curves $p(\gamma_{\pm})$
cross the lima\c{c}on.  First we show that $p(\gamma_-)$ intersects
the lima\c{c}on as shown in the figures by showing that $p$ is
injective on the left half plane.

\theorem Proposition \PInjects.
The function $p$ is injective on the left half plane.

\proof
Observe that $p:\C \maps \C$ is proper and surjective. Consequently, $p$ maps
closed sets to closed sets. Let ${\cal L}$ denote the closed left half-plane:
$$	{\cal L} = \{ z  \st \re{z} \leq 0 \}.$$
Note that $p$ maps the negative real axis onto itself and the imaginary axis
onto a parabola-shaped curve which intersects only at the origin:
$$	p(\i y) = |y|^{2 \alpha} + \i y.$$
Since $p$ has no singularities on ${\cal L}$, $p$ is open and orientation
preserving on ${\cal L}$. In particular, $p(\interior{\cal L})$ is open.
Since $p({\cal L})$ is closed, we conclude that $p(\interior{\cal L})$ is
contained in the component of the complement of the image of the imaginary
axis which contains the negative real axis. Thus $p({\cal L})$ is the
closure of this component, since $p({\cal L})$ is closed. The map $p$ is
proper on ${\cal L}$, and maps ${\cal L}$ onto this component, so the degree
of $p$ is well defined. Since $p^{-1}(0) = \{0\}$, this degree is one. We
conclude that $p$ maps the left half plane diffeomorphically onto the
component described before.   
\QED 

\noindent
We now consider $p(\gamma_+)$:

\item{\smbull} {\it The image of the intersections of $\gamma_+$ and $\delta$:}
At these two points $Dp$ has $0$ as a double eigenvalue.  These points can
be explicitly computed, and one finds the rank there to be one.  This
explains the two points of tangency in the image of the two intersections.

\item{\smbull} {\it The points where the tangent line to $\gamma_+$ is in
the kernel of $Dp$:} One can calculate that there are exactly three such
points, one of which is real (the point $c_1$ described in \FourFixedPts)
and the other two are complex conjugates.  This explains the three cusp
points.

\item{\smbull} {\it The points where the tangent line to $\gamma_+$ is
horizontal or vertical:} One can calculate that there is only one point,
$c_1$, which has a horizontal tangency.  When $\alpha<2$ there is only one
point (the point $c_0$ described in \FourFixedPts) at which $p(\gamma_+)$
has a vertical tangency.

\noindent
Using the above one can show that $p(\gamma_+)$ intersects the 
lima\c{c}on as indicated in the pictures.  It is not hard to show
that $p(\gamma_+)$ and $p(\gamma_-)$ do not intersect.

\subsection{Hopf Bifurcation}
Consider a small disc $D$ with center $c_0\in\delta$ with $Df_{c_0}$ having
complex conjugate eigenvalues of absolute value 1 at one fixed point.  We
wish to discuss the bifurcation picture in this disc. For $c$ in this disc,
we can smoothly parametrize the corresponding fixed point $z(c)$, so that
$z(c_0) = z_0$. When $D$ is small enough this map $z:D\maps\C$ is a
diffeomorphism. In particular, $z(D)$ intersects $\delta$ and $D-\delta$
consists of two regions, one where the fixed point is attracting and one where
the fixed point is repelling.  On the boundary of these regions the fixed
point is neutrally stable. One should in general expect a Hopf bifurcation.
By a Hopf bifurcation, we mean that as $c$ passes through the curve
$\delta$, the fixed point $z_0$ will change stability and an invariant
circle will be created or destroyed (see \cite{MM} or \cite{D}). This
behavior is more precisely described in terms of normal forms.

\smallskip
Assume that we have chosen $z_0$ so that its eigenvalues are non-resonant:
not first, second, third, or fourth roots of unity. Then one can change
coordinates \cite{MM}, depending smoothly on the parameter $c$, so that one
obtains the following normal form for $f_{\alpha, c}$, for $c \in D$:
$$
F_c(z) = \lambda_c z (1 + v_c |z|^2) + O(z^5)
$$
Notice that the eigenvalue $\lambda_c$ and the coefficient $v_c$ depend
also on $\alpha$. The map $c \mapsto \lambda_c$ is a diffeomorphism on $D$
and intersects the unit circle. The bifurcation theory for $c$ near $c_0$
depends on $\re{v_{c_0}}$, provided it does not vanish.

\theorem Claim \HopfDirection.
Assume that the eigenvalue $\lambda_{c_0}$ is non-resonant. Then:
\item{\smbull} When $1/2 < \alpha < 1$, $\re{v_{c_0}}$ is positive.
\item{\smbull} When $\alpha = 1$, $v_{c_0}$ vanishes.
\item{\smbull} When $\alpha > 1$,  $\re{v_{c_0}}$ is negative.

\midinsert{
  \centerline{\psfig{figure=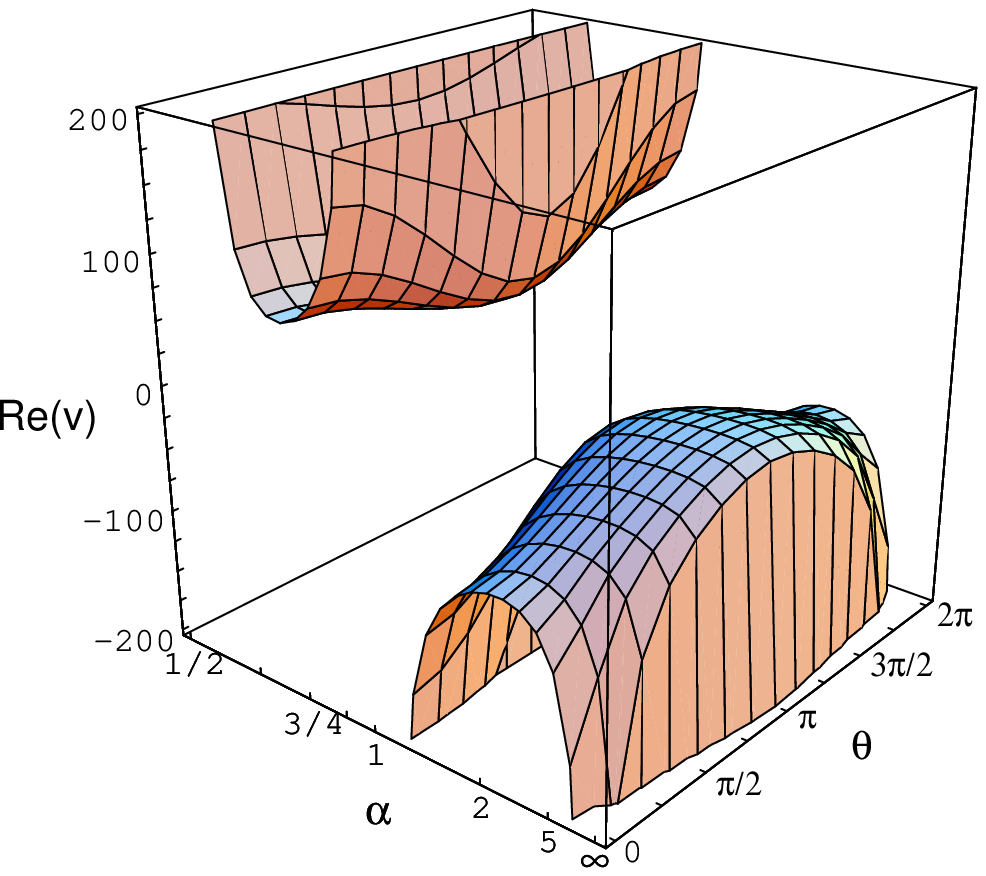,height=.48\hsize}\hfil
	      \psfig{figure=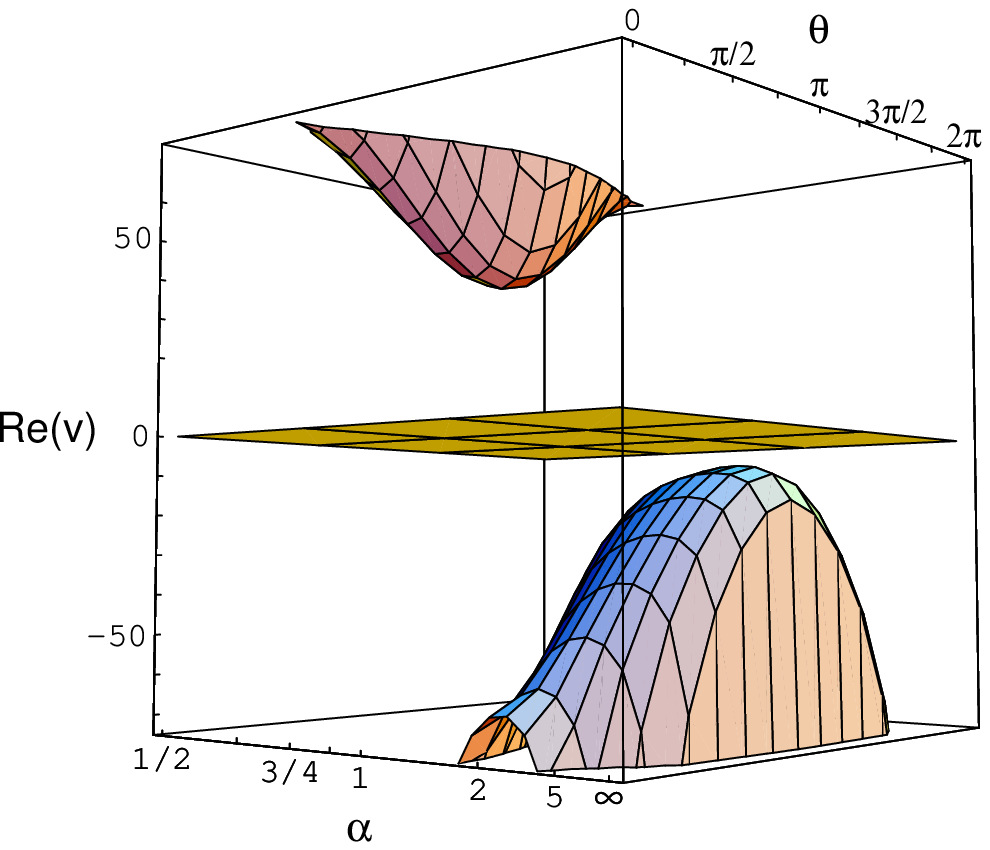,height=.48\hsize}}
   \longcaption{Graphs of $\re{v_{c_0}}$ as a function of $\alpha$ and
$\theta=\arg{\lambda_{c_0}}$.  We have modified the $\alpha$ scale so that
the intervals $({1\over2}, 1)$ and $(1, \infty)$ have the same length.  On
the right is a closeup near the $\alpha$-$\theta$ plane, which we have
shaded to emphasize the claim. \namefigure{\HopfPics}} 
}\endinsert

\smallbreak\noindent{\it Justification.\enspace}
When $\alpha = 1$, this is obvious. For any other values of $\alpha$ there 
seems to be no easy proof.  The only more or less straightforward case is an
infinitesimal computation near the holomorphic case $\alpha = 1$.
Conceivably, a computer-assisted proof of this could be done using interval
arithmetic.   However, we feel this claim does not merit the effort of a
difficult and tedious proof, and have used {\it Mathematica} \cite{W} to
perform the coordinate changes and compute $v_{c_0}$ on a large grid of
parameter values.  See \HopfPics.  For $\alpha < 1$, we obtained
(numerically) that $\re{v_{c_0}} > 41.487$, and for $\alpha >  1$ that 
$\re{v_{c_0}} < -8.594$. 
\QED

By \HopfDirection, the sign of the real part of the coefficient $v_{c_0}$
depends only on $\alpha$,  and hence we know in which direction the Hopf 
bifurcation occurs. Assuming that the disc $D$ is small enough, we have the
following dichotomy:

\item{\smbull} When $\alpha < 1$ and $|\lambda_c| < 1$, there exists an
invariant circle near $z(c)$ which is repelling in the normal direction. 
For $\lambda_c$ outside the closed unit disc, there is no invariant circle 
near the point $z(c)$. 

\item{\smbull} When $\alpha > 1$, we have the opposite situation: for
$|\lambda_c| < 1$, there is no invariant circle close to $z(c)$, and for
$\lambda_c$ outside the closed unit disc, there exists an invariant circle
which is attracting in the normal direction.  We conjecture that the
critical point is still attracted to this circle.

	\ifSingleSecs \vfil\eject \fi

	\ifx\macrosLoaded\undefined 
        
	\pageno=20
\fi
\secno=4 

\section{Remarks on the Topology of the Connectedness Locus}

In the holomorphic case ($\alpha = 1$), the connectedness locus is usually
referred to as the Mandelbrot set, and is known to be connected \cite{DH};
its complement in the Riemann sphere is conformally equivalent to the open
disk. Furthermore, every connected component of its interior is a
topological disk; these disks are referred to as either {\bit hyperbolic
components} or {\bit queer components}.  For any map lying within a
hyperbolic component, there is a periodic attractor which necessarily
attracts the critical point, and each hyperbolic component has a special
point called the {\bit center} (the parameter for which the critical orbit
is periodic).  Within any connected component of the interior of the
Mandelbrot set, all maps except possibly one are topologically (in fact,
quasi-conformally) conjugate; this exception is the center of a hyperbolic
component. 

A long standing conjecture is that there are no ``queer components'' in the
Mandelbrot set; that is, all components of the interior are hyperbolic. This
hyperbolicity conjecture is implied by the local connectivity of the
Mandelbrot set (see \cite{DH}). Yoccoz has shown that local connectivity
holds for a  ``substantial'' part of Mandelbrot set \cite{Hu}, and Swiatek
\cite{S}, McMullen \cite{McM}, and Lyubich \cite{Ly} have established the
hyperbolicity conjecture along the real line by rather different techniques.

The situation when $\alpha \neq 1$ is quite different, as one should
expect, because the iterates of the maps are not uniformly
quasi-conformal.  Douady and Hubbard's proof that the Mandelbrot set
is connected relies on the conformal structure; we see no way to adapt
it to the non-conformal case.  Furthermore, there is no mathematical
relationship between the hyperbolicity conjecture and the local
connectivity in this case.  However, numerical evidence strongly
suggests the following.

\theorem Conjecture \ConnLocusConnected.
For all $\alpha > 1/2$, the connectedness locus $\ConnLocus$ is
connected.  However, $\ConnLocus$ is not locally connected for
all $\alpha \ne 1$. 

\midinsert{
  \hbox to \hsize {
   \vsize=.4\hsize
   \hbox to .468\hsize{\vbox to \vsize{
     \vfil
     \psfig{figure=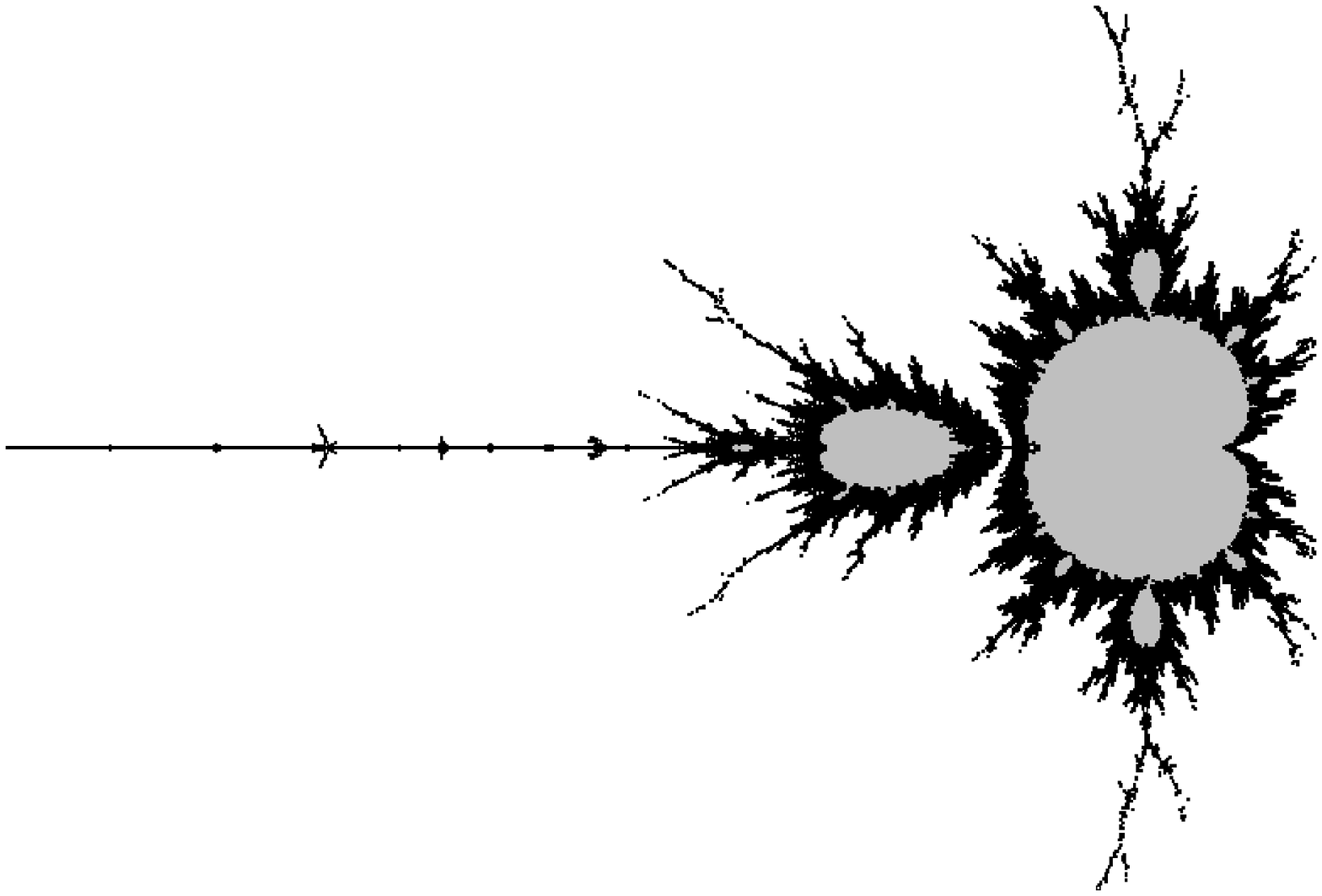,width=.468\hsize}
     \vfil
   }}
   \hfil
   \hbox to .234\hsize{\vbox to \vsize{
     \vfil
     \psfig{figure=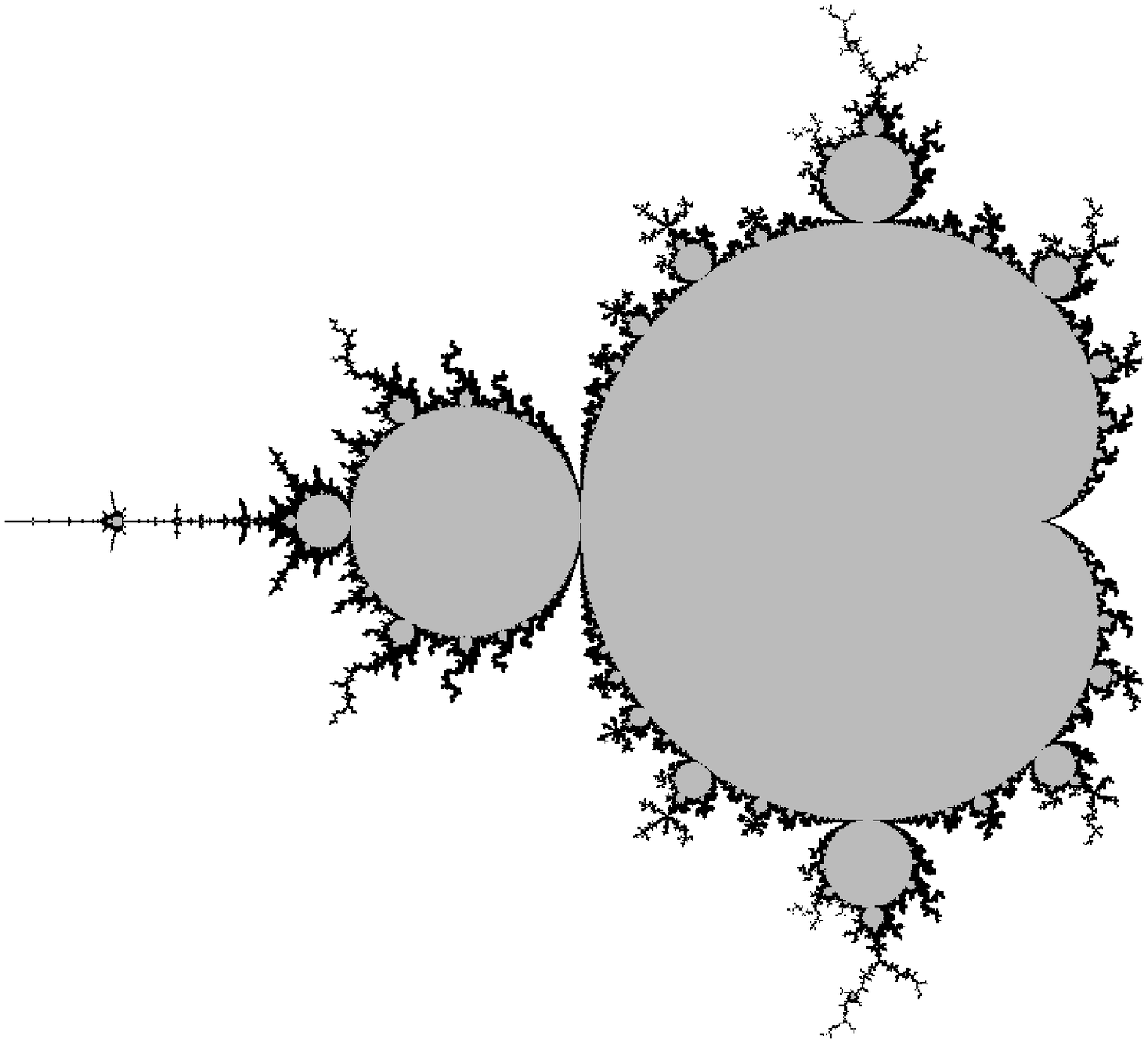,width=.234\hsize}
     \vfil
   }}
   \hfil
   \hbox to .162\hsize{\vbox to \vsize{
     \vfil
     \psfig{figure=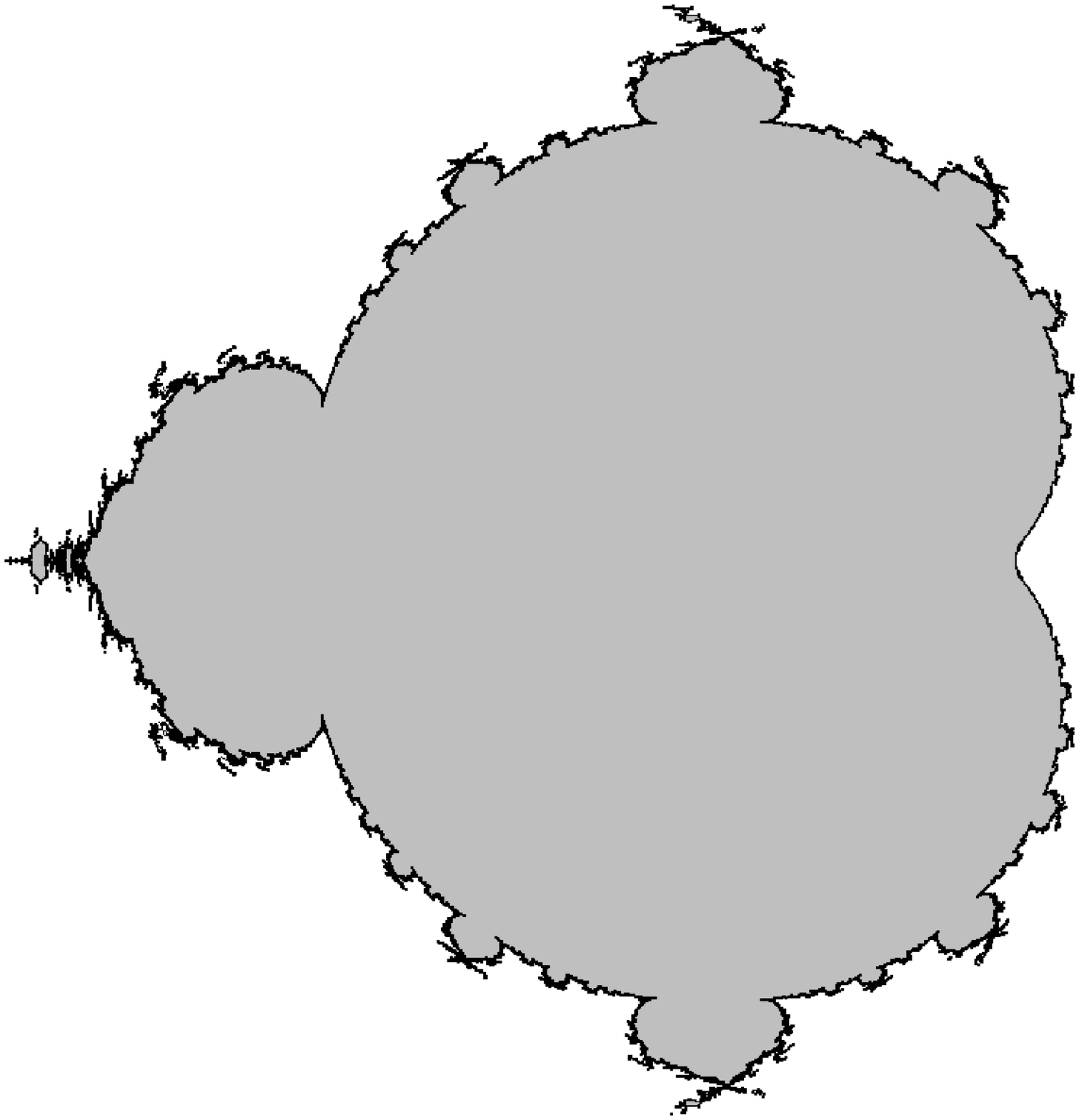,width=.162\hsize}
     \vfil
   }}
  }
  \shortcaption{The connectedness loci $\ConnLocus$ for $\alpha=.75$,
                $\alpha=1$, and $\alpha=1.5$.} 
}\endinsert

The apparent lack of local connectivity of $\ConnLocus$ in the
non-conformal case is at least partially related to the presence of
saddle points and their stable and unstable manifolds. Such invariant
saddles make a qualitative understanding of the dynamics difficult and
a quantitative understanding nearly impossible. In particular, when
$\alpha < 1$, one can readily see the difficulty caused by the
saddles.

\namenextfigure{\CLocusBlowup}  
Specifically, consider the interval of real parameters for which
$f_{\alpha,c}|_\R$ has an attracting (in $\R$) fixed point which attracts
the critical point. For some interval of parameters, this fixed point is
not an attracting fixed point on $\C$, but is repelling in the imaginary
direction (for example, in the region~\OneRep\Saddle\ discussed in
section~4).  Denote this fixed point by $z_c$ and consider its global
stable manifold $W^s(z_c)$. Because the dynamics is non-invertible, this
global stable manifold is topologically more complicated than for a
diffeomorphism.  The critical point is in this stable manifold and one
might hope that the filled-in Julia set is the closure of this stable
manifold. Now consider a parameter value $c'$ which is nearby, but not
real. Consider the corresponding fixed point $z_{c'}$ and the corresponding
global stable manifold. The critical point is not necessarily contained in
this global stable manifold. In fact, the global stable manifold changes
with the parameter; sometimes the critical point escapes to infinity and
sometimes it is in $W^s(z_c')$. The detailed structure of the
connectedness locus is unclear, but it has the topological appearance of a
stable manifold. The rough structure of the $\ConnLocus$ for these
parameters is that the main lobe (which contains those values of the $c$
for which there is an attracting fixed point) is connected to the period
two lobe, (containing those values for which there is an attracting cycle
of period two) are connected only by a segment in the real line. A very
complicated comb-like structure limits on part of this segment.
See \CLocusBlowup.

\medskip
\midinsert{
  \hbox to \hsize {
   \vsize=.55\hsize
   \hbox to 1.05\vsize{\vbox to \vsize{
     \vfil
     \psfig{figure=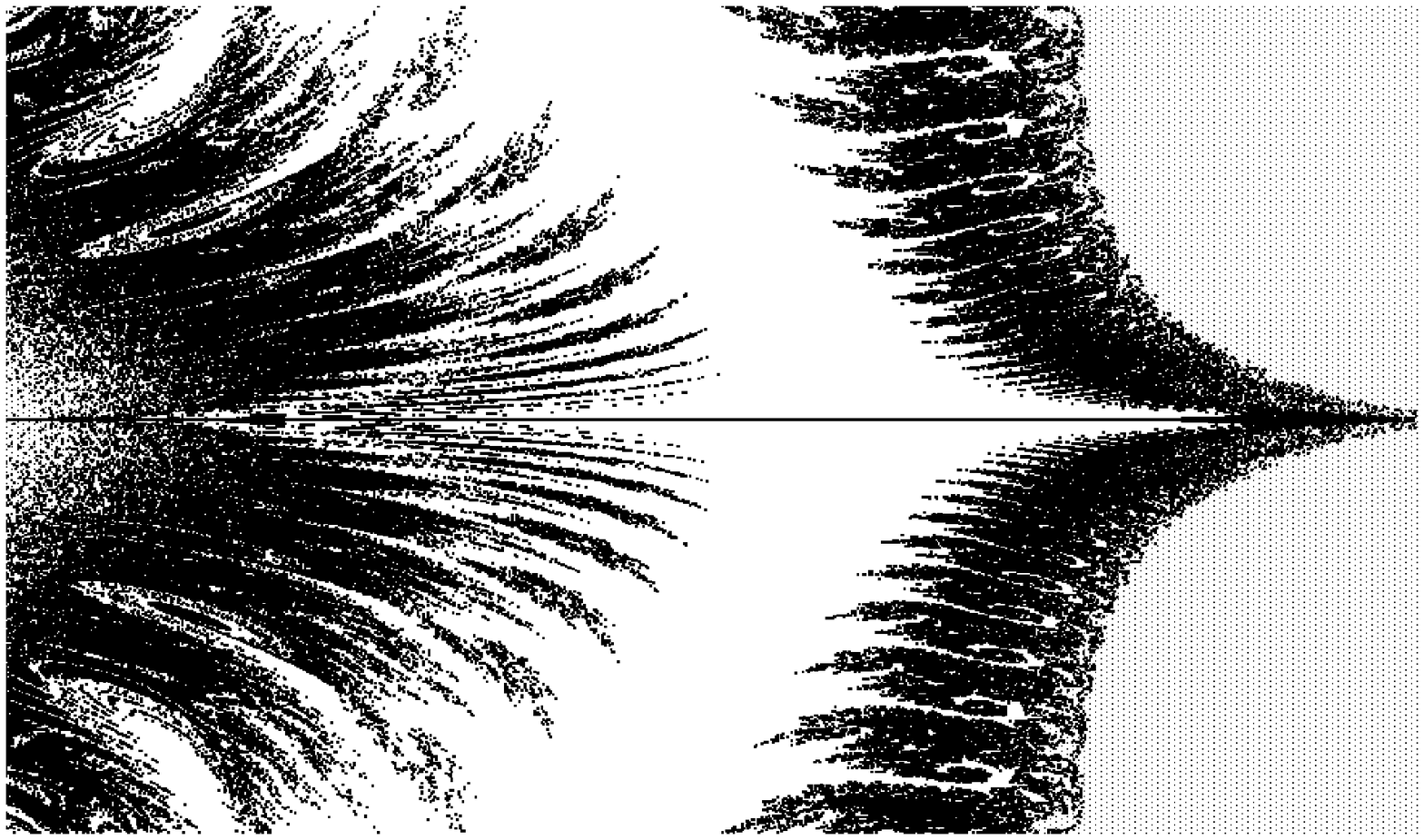,width=\vsize}
     \vfil
   }}
   \hfil
   \psfig{figure=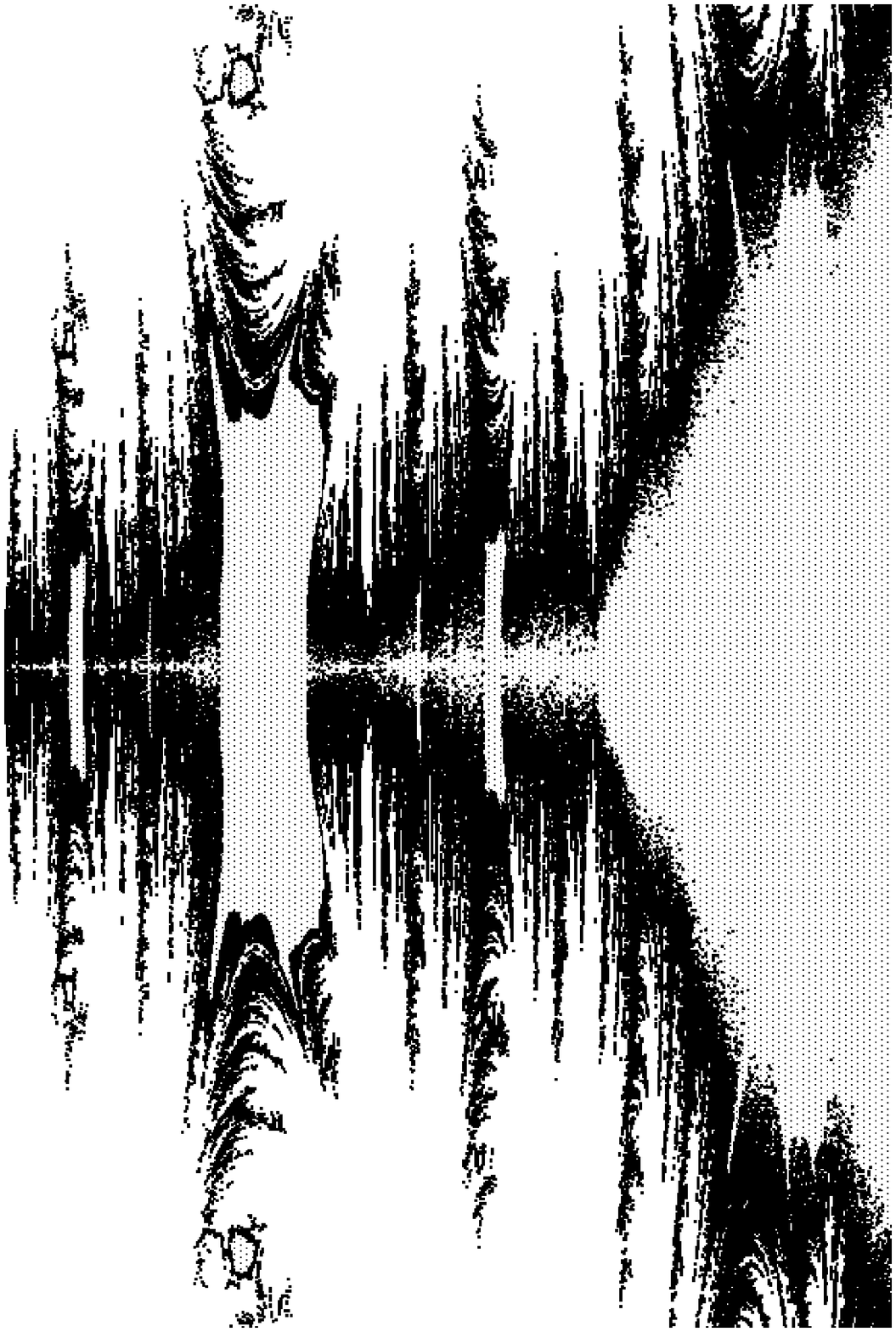,height=\vsize}
  }
  \longcaption{Blowups of the $\alpha=.75$ and $\alpha=1.5$ connectedness
               loci showing the apparent lack of local connectivity.}   
}\endinsert

When $\alpha > 1$,  there is also an apparent lack of local connectivity
near the real line, but in a dynamically different part of $\ConnLocus$
(for example, between the limit of period doubling and the creation of
period~5 orbit; see \CLocusBlowup).  At this point, we have no real
understanding of what causes this, however.

\smallskip
Some insight into the topology of the boundary of the main lobe of
$\ConnLocus$ can be gained by looking again at the Hopf bifurcation near
the conformal case.
We shall analyze the type of bifurcations which occur near parameters for
which there is a fixed point of multiplier $\lambda$, where $\lambda$ is a
$q^{\rm th}$ root of unity ($q \ge 5$).  We shall omit most of the tedious
calculations here; the interested reader should refer to section~5 of
\cite{BSTV}. 

In this situation, one can change coordinates \cite{MM} so that we have a
two complex-parameter family of maps defined near the origin in the complex
plane: 
	$$F_{\mu,a}(z) =  \lambda z (e^\mu + a |z|^2 + z^q) + 
	  O( |a z^4|, |a\mu z^3|, |\mu z^{q+1}|, |z^{q+2}|).$$
(Here $a \sim (1 - \alpha)/\lambda$.) We are interested in the $q$-periodic
points of the map $F_{\mu,a}$; one easily sees that for such a point we have 
	$$z  = z (e^{q\mu} + q a |z|^2 + q z^q) + 
	  O( |a z^4|, |a \mu z^3|, |\mu z^{q+1}|, |z^{q+2}|).$$

We first consider the parameter $a$ to be real, negative, small and fixed.
(this corresponds to $\alpha > 1$.)  When $\re \mu \leq 0$, the fixed point
$z=0$ is an attractor; the product of its eigenvalues is less than 1. When
$\re\mu > 0$, the fixed point is repelling, but for $\re \mu$ sufficiently
small, there exists an attracting, invariant (Hopf) circle whose diameter is
of order $\sqrt{{{\re \mu}\over{|a|}}}$. One can show that there is also a
$q$-periodic orbit located approximately on the circle of radius
$|a|^{{{1}\over{q - 2}}}$. When $|\mu|$ is small, it is easily seen that this
orbit is repelling.

In a horn shaped domain in the $\mu$-plane, the rotation number on the
invariant circle is $p/q$ (recall $\lambda = \e^{2\pi\i p/q}$);
this horn is in fact the $p/q$-resonance horn or Arnol'd tongue
\cite{Ar, ACHM, Ha}. Within this horn, there are two additional
${q}$-periodic orbits: one is a saddle and the other is an attractor; the
invariant circle around the repelling fixed point is the closure of the
unstable manifold of the $p/q$-saddle, which contains the attracting
orbit.  If $\mu$ leaves the horn ``through the side'', i.e., fix $\re\mu$
and vary $\im \mu$, the saddle and the attractor collide, and although there
is still an invariant circle, the rotation number is no longer $\arg
\lambda$.  If instead we allow  $\re \mu$ to increase sufficiently, the
saddle and the repeller collide, leaving a single attracting orbit of period
$q$.   However, before this collision occurs, the invariant circle looses
smoothness and becomes only a topological circle.  This loss of smoothness
occurs when the eigenvalues of the $p/q$-sink become complex.  See
\cite{ACHM, \S~8}. 

\midinsert{
 \hbox to \hsize{
    \hfil
    \vtop{\hsize=.4\hsize
    \hbox to \hsize{\psfig{figure=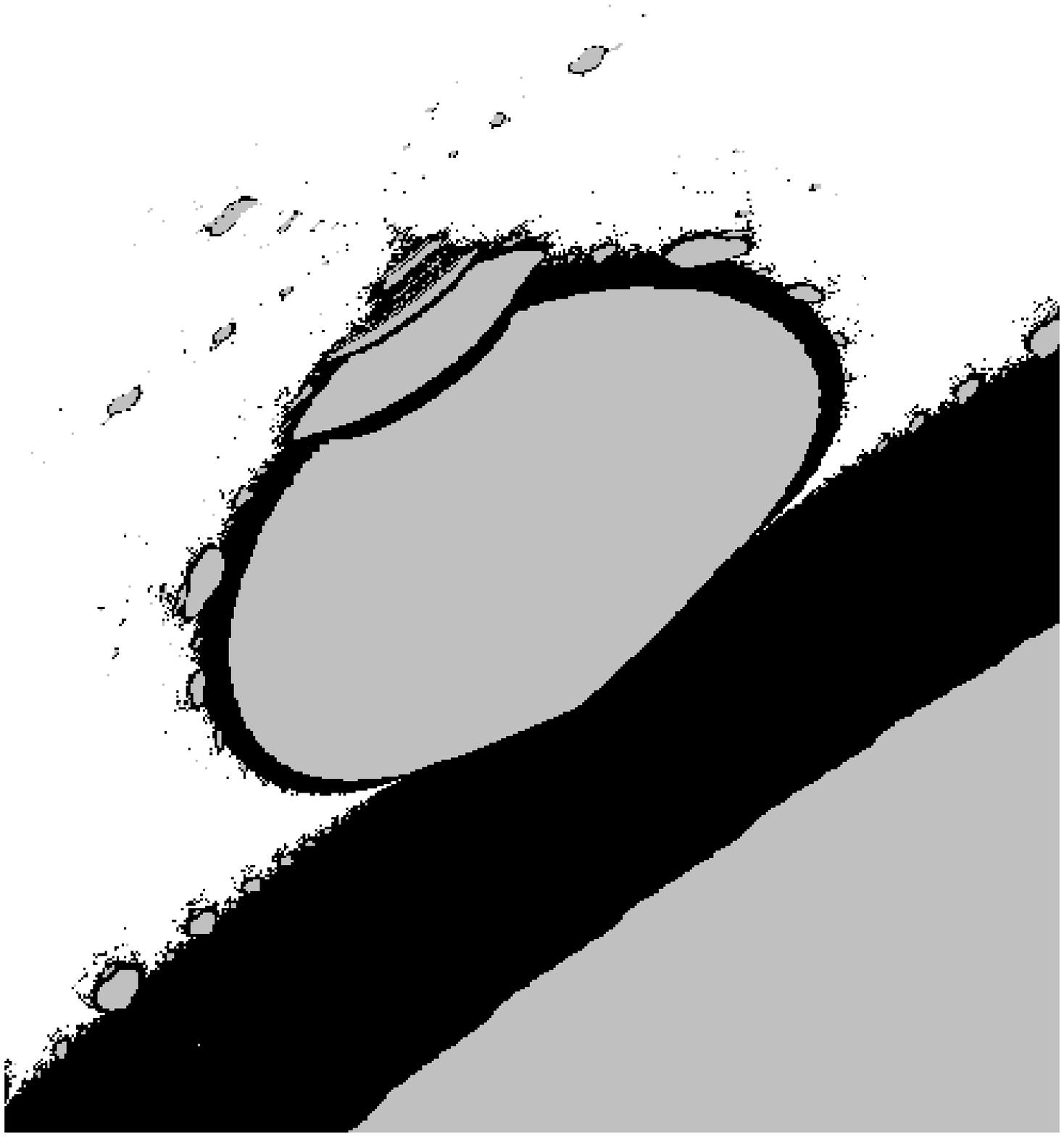,width=\hsize}\hfil}
    \widecaption{The 2/5 limb for $\alpha = 1.5$.  The regions in grey
indicate that the critical point converged to a periodic attractor of
moderate ($<100$) period within a few hundred iterations.  Note the grey
horn-like region at the base. The boundary of the figure appears
disconnected due to the algorithm used to produce the picture; a different
algorithm gives a much ``thicker'' boundary.
\namefigure{\TwoFifthsMore}} 
    } 
 \hfil
 \vtop{\hsize=.4\hsize
    \hbox to\hsize{\hfil\psfig{figure=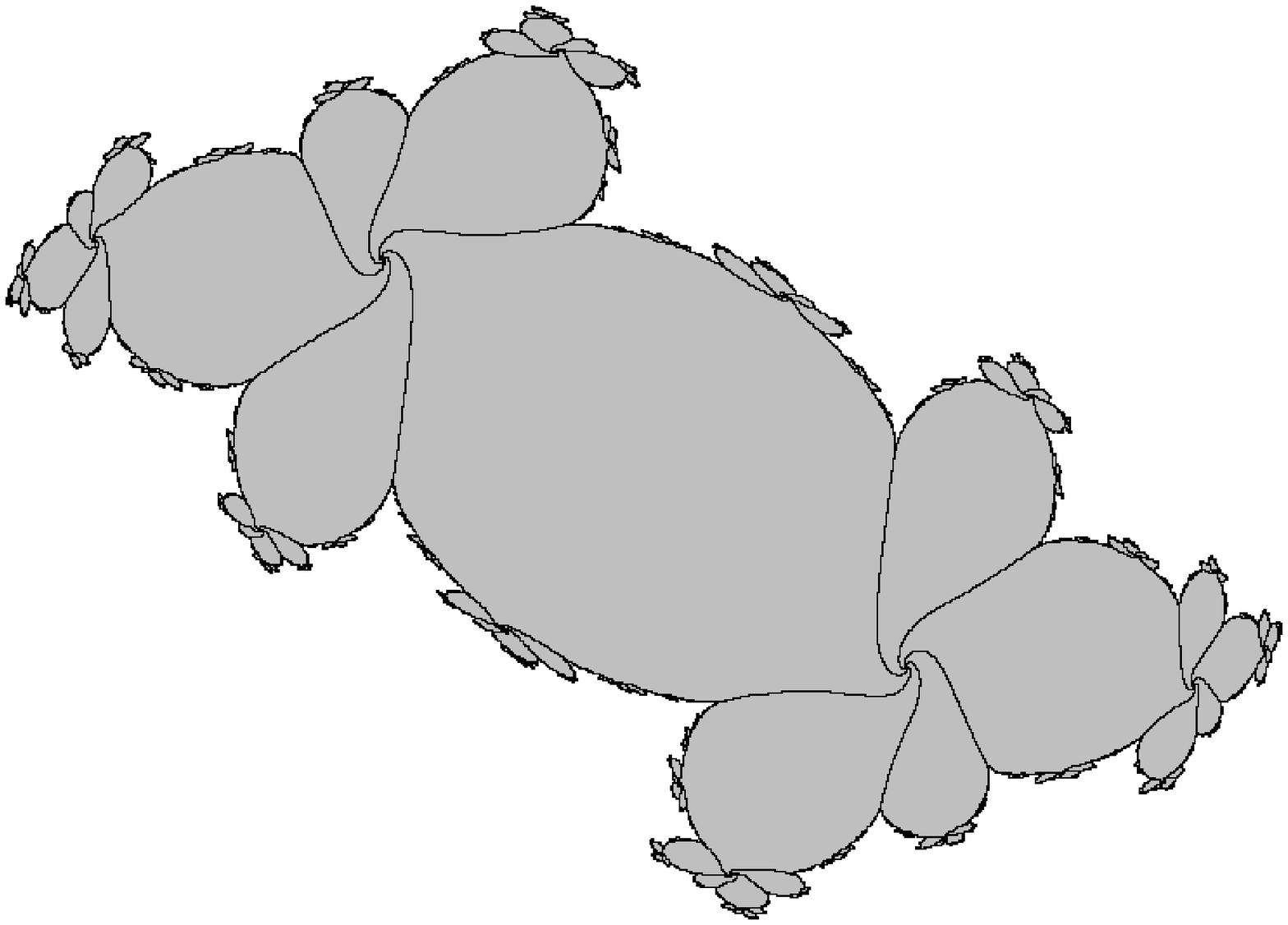,width=\hsize}}
    \widecaption{A Julia set for $c$ within the 2/5 resonance horn for
    $\alpha=1.5$.  The large grey areas are the basin of attraction of the
    $2/5$ attracting periodic orbit.  This attracting orbit has complex
    eigenvalues, and so the unstable manifold of the $2/5$ periodic saddle
    does not form a smooth invariant circle.  This saddle orbit lies on the
    five smooth curves which divide the grey regions; these curves are the
    stable manifold of the saddle.}}
 \hfil
 }
 }
\endinsert

Experimental evidence indicates that the critical orbit remains bounded
for all parameter values discussed above:  it is attracted to either the
attracting fixed point ($\re\mu<0$), the invariant circle, or the
$p/q$-periodic attractor.  This gives some explanation for the appearance of
the connectedness locus near the main lobe: the $p/q$-lobe sits at the end
of a resonance horn, and is hence attached along an arc of values. See
\TwoFifthsMore. 

\smallskip
When $a$ is positive (corresponding to $\alpha<1$), the picture is the other
way around. As above, the fixed point is attracting for $\re\mu<0$ and
repelling for $\re\mu>0$, but the invariant circle is repelling and exists
only when $\re\mu<0$.  Within a horn of $\mu$ values, the invariant circle
contains a $p/q$-saddle and a $p/q$ repeller, and is the closure of the
stable manifold of the saddle. For all $\re\mu$ negative and $\mu$
sufficiently small, there is another $q$-periodic orbit nearby, which is
attracting. 

\midinsert{
 \centerline{\psfig{figure=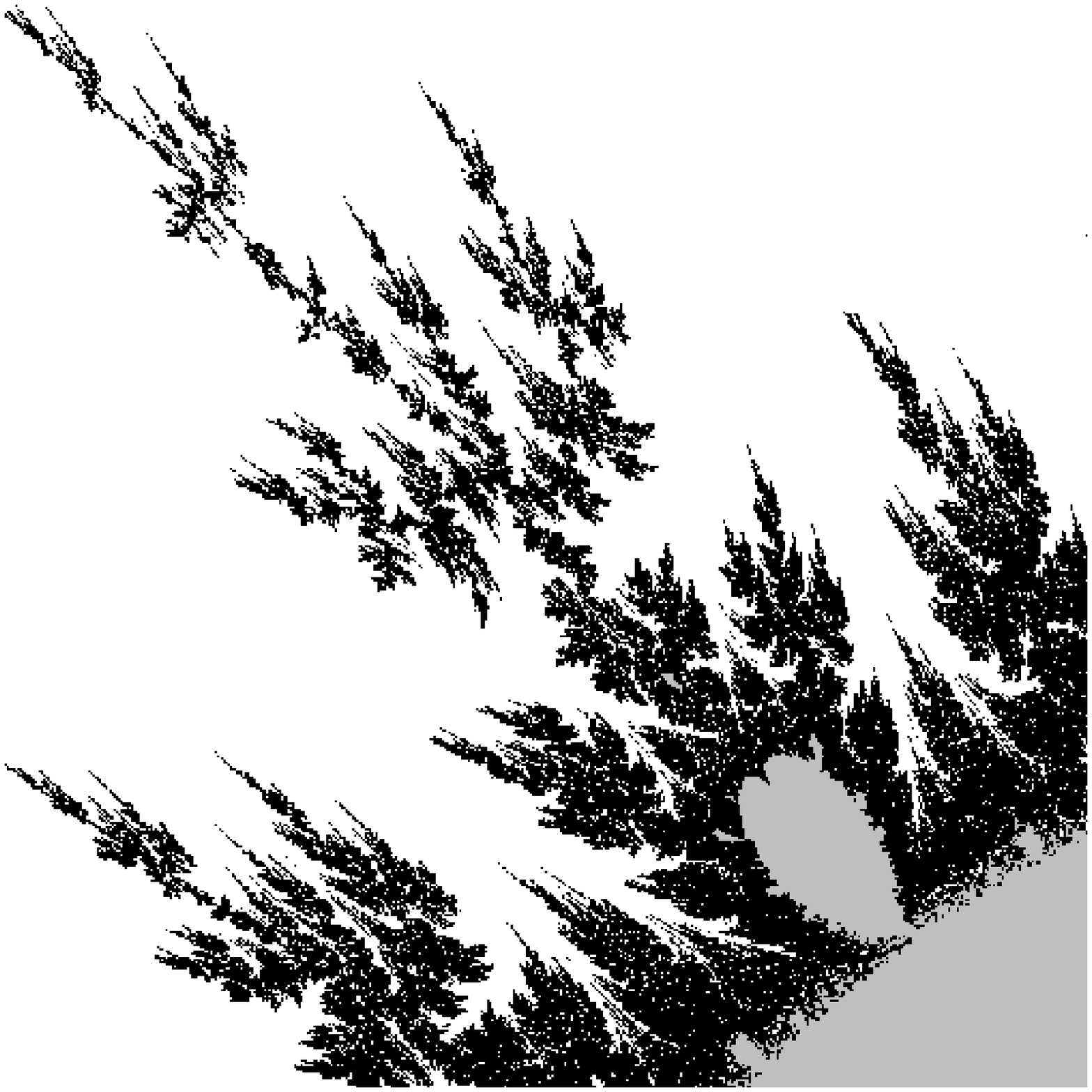,width=.45\hsize} \hfil
	     \psfig{figure=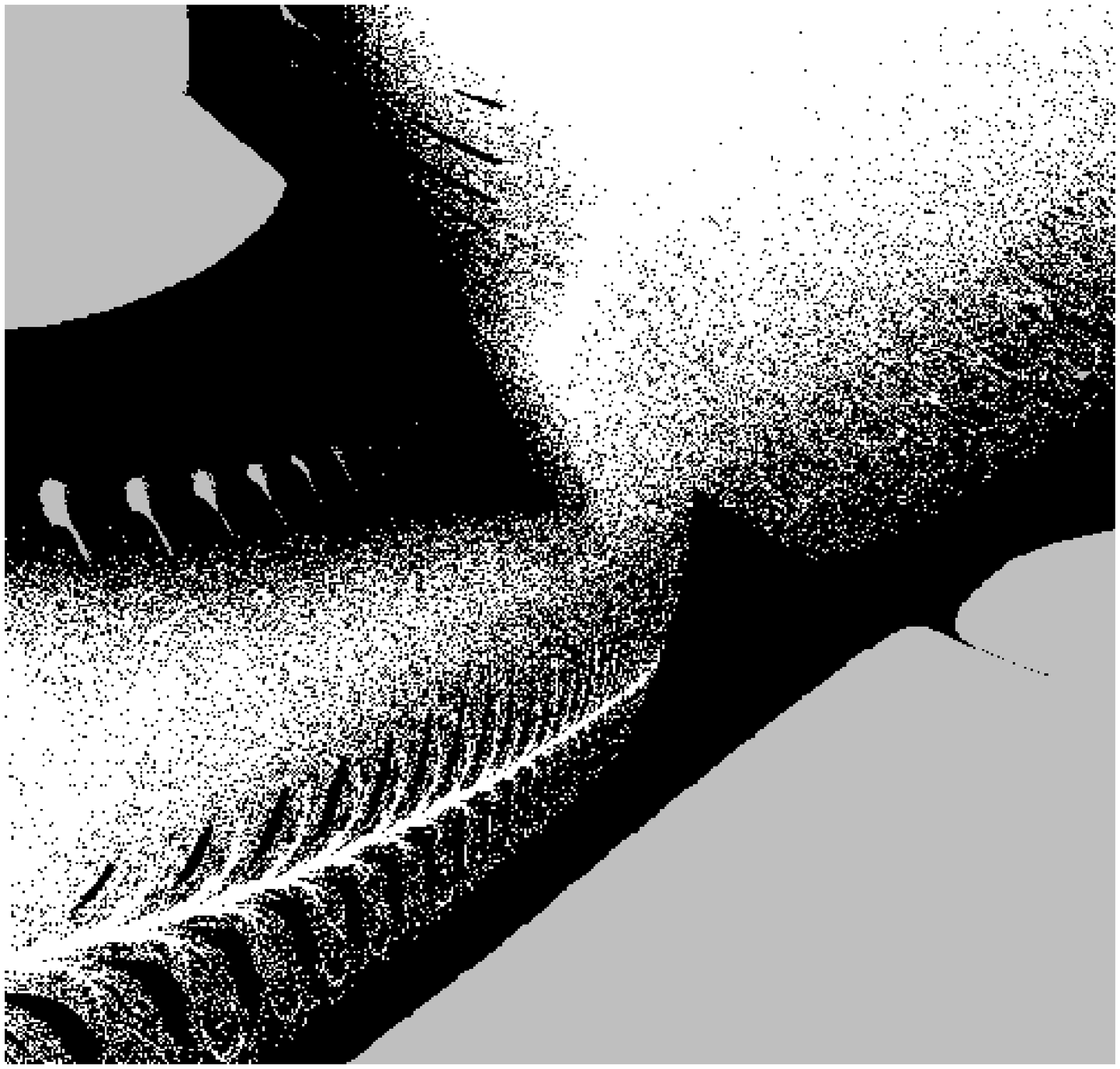,width=.45\hsize}}
 \longcaption{The 2/5 limb for $\alpha =0.75$, and a blowup (at right).
The two figures were produced with different algorithms; in the left-hand
picture, grey denotes parameters for which the critical point failed to
escape within 256 iterations, while on the right such parameter values are
colored black, and grey is used to indicate convergence of 0 to an
attracting orbit of moderate period. \namefigure{\TwoFifthsLess}} 
 }
\endinsert

However, in this case the relationship between the Arnol'd tongue and the
connectedness locus is quite different.  Since the circle is repelling, for
many parameter values in the horn, the critical orbit does not limit on the
attracting fixed point; it can escape to $\infty$, and hence the filled-in
Julia set will be disconnected.  Thus, one cannot readily detect the presence
of the Arnol'd tongues from the connectedness locus alone, as in the case of
$\alpha>1$.  

\namenextfigure{\JulInHorn}  
There is, however, a horn-like structure which is readily
apparent along the boundary of $\ConnLocus$.  See \TwoFifthsLess.  
This is related to the presence of $q$-periodic saddles.  Near this horn,
there are three period $q$ orbits, as well as the attracting fixed point.  Two
of the periodic orbits are repelling, and the the other is a saddle.  The
horn corresponds to the parameter values for which the critical point lies
between the one side of the stable manifold of a saddle point $z$ and the
other side of the stable manifold for its image $f_c(z)$.  See \JulInHorn.
At the point of the horn, a saddle connection occurs: one side of the local
stable manifold of the periodic saddle $z$ is the local unstable manifold
for its image $f_c(z)$.

\midinsert{
 \hbox to \hsize {
  \hfil
  \vtop{\hsize=.4\hsize
 \hbox to \hsize{\psfig{figure=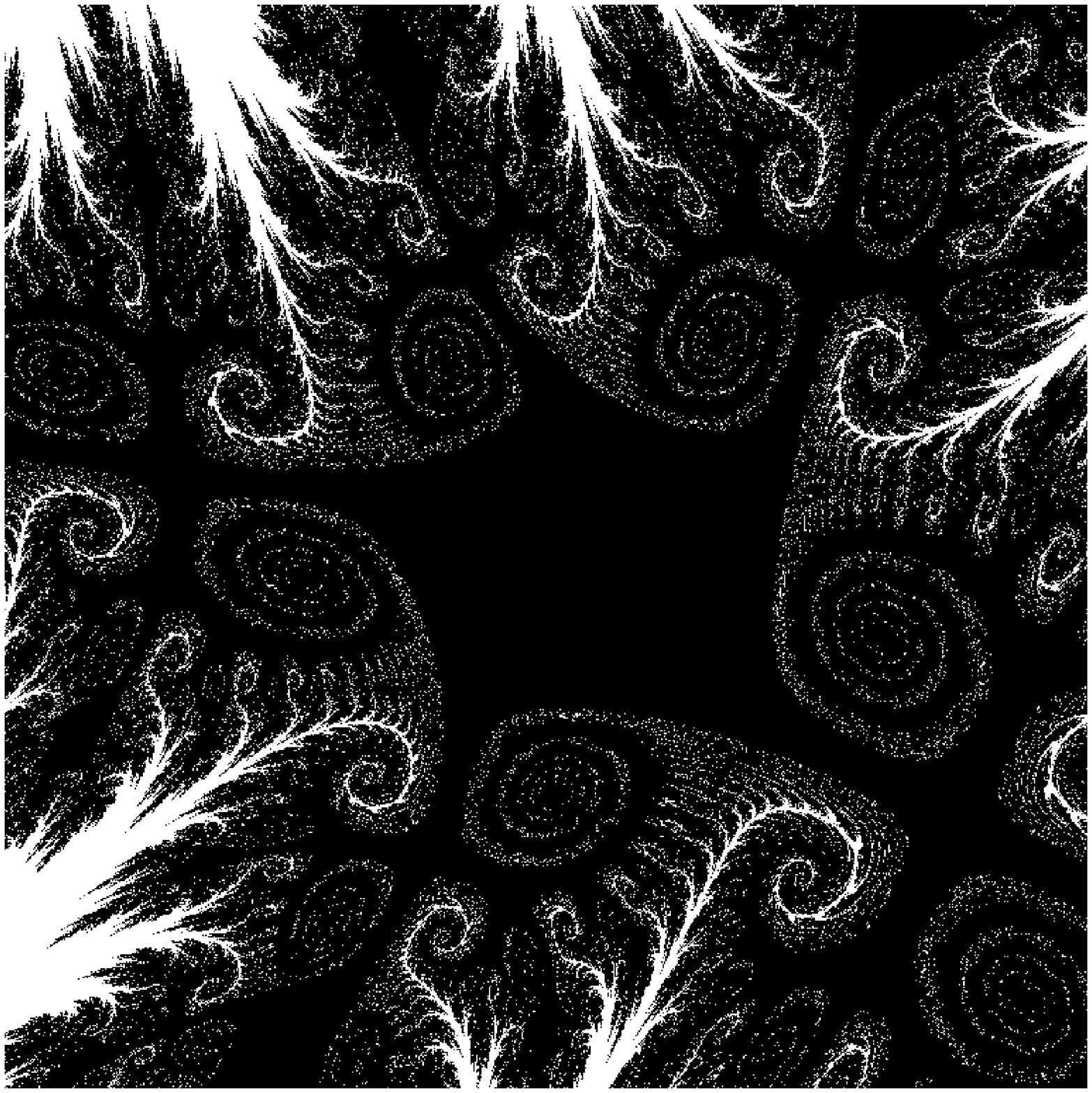,width=\hsize}\hfil}
 \widecaption{A close-up of the filled Julia set for $\alpha=.75$ and $c$
within the horn in \TwoFifthsLess.  There is an attracting fixed point at
the center of the picture, with a pair of period~5 repellers surrounding it
(at the ends of the black and white spirals nearest the center), and a
period~5 saddle at the edge of the large black region.  The critical point,
$0$, lies in the large black cross shaped region near the lower right.} 
}
 \hfil
 \vtop{\hsize=.4\hsize
 \hbox to \hsize{\hfil\psfig{figure=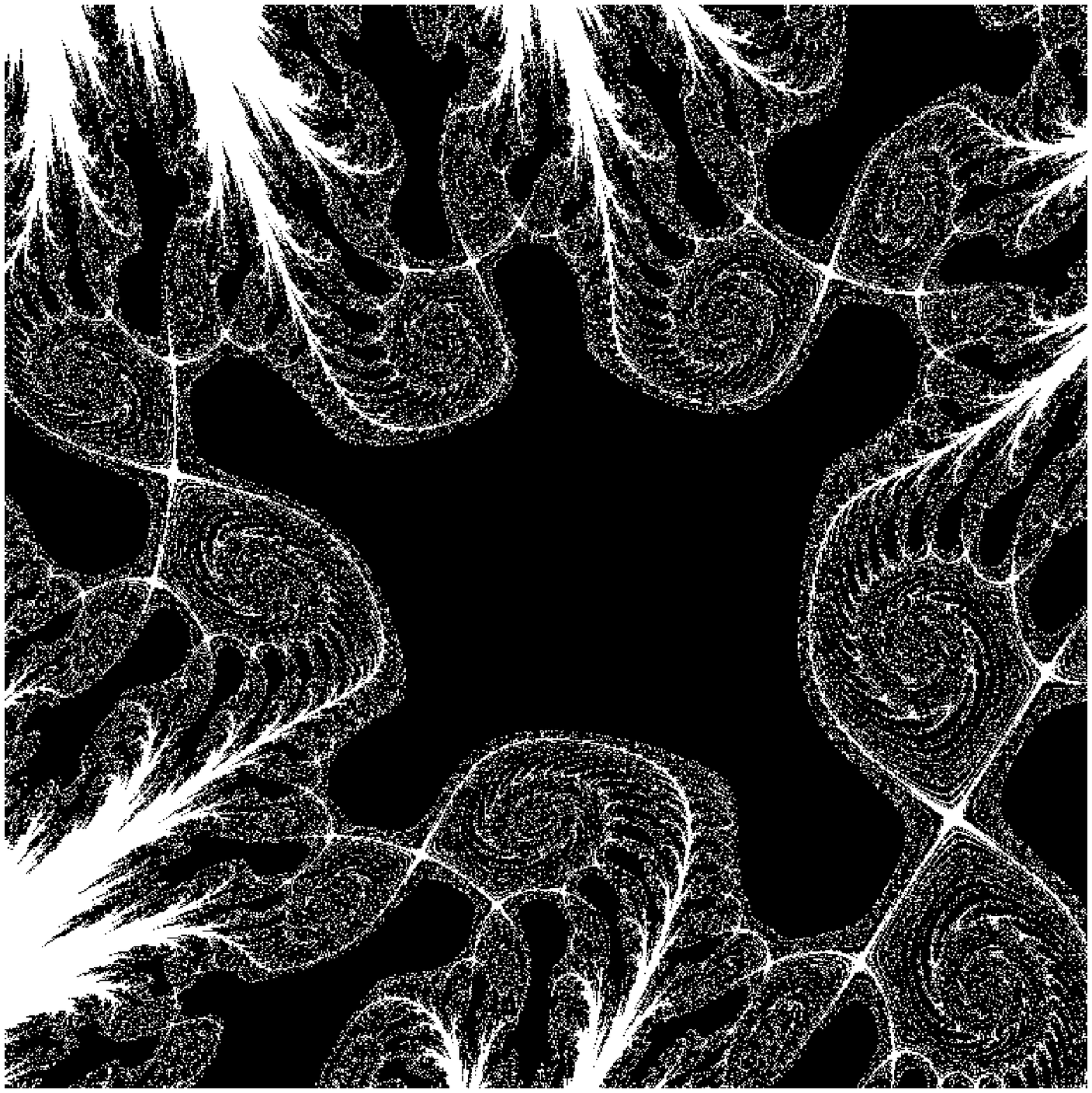,width=\hsize}}
 \widecaption{A close-up of the filled Julia set for $\alpha=.75$ and $c$
	just outside the horn in \TwoFifthsLess.  As in \JulInHorn, there is
 	period~5 saddle, two period~5 repellers, and an attracting fixed
        point, in approximately the same locations.  However, in this case
        the attractor fails to attract the critical point; 0 iterates to
        infinity, and the filled Julia set is disconnected, although not
        totally disconnected.} 
  }
  \hfil
}
}
\endinsert

Attached to the top of the horn is a curve for which the critical orbit
remains bounded, although it is not attracted to an attractor.  For these
parameters, the critical point lies on the stable manifold of one of the
points of the period-$q$ saddle orbit discussed above.  This orbit appears
to persist long enough to attach the period-$q$ lobe (within which there is
an attracting period $q$ orbit) to the main lobe.  Thus, $\ConnLocus$ is
not disconnected as it appears in \TwoFifthsLess.

	\ifSingleSecs \vfil\eject \fi

%
	\ifx\macrosLoaded\undefined 
        
\fi
\secno=-1 

\section{References}

{
\parskip=.75\baselineskip
\frenchspacing

\item{[A]} L. Ahlfors:
	{\sl Lectures on Quasiconformal Mappings}.
	Wadsworth, 1987.

\item{[ACHM]} D. Aronson, M. Chory, G. R. Hall, R. McGehee:
	Bifurcations from an Invariant Circle for Two-Parameter Families of
	Maps of the Plane: A Computer-Assisted Study. 
	{\sl Comm. Math. Phys.} {\bf 83} (1982), 303--354.

\item{[Ar]} V. I. Arnol'd:
	Small Denominators I. On the Mappings of the Circumference of the
        Circle onto Itself. {\sl Translations Am. Math. Soc., 2nd series}
	{\bf 46} (1965), 213--284.
	
\item{[BSTV]} B. Bielefeld, S. Sutherland, F. Tangerman, and J.J.P. Veerman:
	{\sl Dynamics of Certain Non-Conformal Degree Two Maps of the Plane}.
	Stony Brook Institute for Math. Sciences preprint \#1991/18.
	
\item{[Bl]} P. Blanchard:
	Complex Analytic Dynamics on the Riemann Sphere, 
	{\sl Bull. Amer. Math. Soc. {\bf 11}} (1984), 85-141.

\item{[D]} R. Devaney:
	{\sl An Introduction to Chaotic Dynamical Systems}, $2^{\rm nd}$ ed.
        Addison-Wesley, 1989.

\item{[DH]} A. Douady, and J. Hubbard:
	{\sl Etude Dynamique des Polyn\^omes Complexes I, II}, 
	Publications Math\'ematiques d'Orsay {\bf 84--02} (1984), and
        {\bf 85--04} (1985).

\item{[Ha]} G. R. Hall:
	Resonance Zones in Two-Parameter Families of Circle Homeomorphisms,
	{\sl SIAM J. Math. Anal.} {\bf 15} (1984), 1075--1081.

\item{[HPS]} M. Hirsch, C. Pugh, and M. Shub,
	{\sl Invariant Manifolds,}
	Lecture Notes in Mathematics, Vol. 583. 
	Springer-Verlag, 1977.

\item{[Hu]} J.H. Hubbard:  
	 Local Connectivity of Julia Sets and Bifurcation Loci: Three
        Theorems of J.-C. Yoccoz, pp. 467--511 in
	{\sl Topological Methods in Modern Mathematics, A Symposium in
	 Honor of John Milnor's 60th Birthday},
	 Publish or Perish, 1993.

\item{[J1]} Y. Jiang,
	Generalized Ulam-von Neumann Transformations,
	Thesis, CUNY Graduate Center, 1990.

\item{[J2]} Y. Jiang,
	Dynamics of Certain Nonconformal Semigroups,
	{\sl Complex Variables} {\bf 22} (1993), 27--34.

\item{[L]} O. Lehto,
	{\sl Univalent Functions and Teichm\"uller Spaces}.
	Springer-Verlag, 1987.

\item{[Ly]} M. Lyubich, 
	{\sl  Geometry of Quadratic Polynomials:  Moduli, Rigidity and Local
        Connectivity},
	Stony Brook Institute for Math. Sciences preprint \#1993/9.

\item{[McM]} C. McMullen,
	{\sl Complex Dynamics and Renormalization},
	preprint, U.C. Berkeley Math. Dept., 1993.

\item{[MSS]} R. Ma\~n\'e, P. Sad, and D. Sullivan,
	On the Dynamics of Rational Maps,
	{\sl Ann. Scient. \'Ec. Norm. Sup., $4^e$ s\'erie}, {\bf 16} (1983),
	193--217.
        
\item{[MM]} J. Marsden and M. McCracken,
	{\sl The Hopf Bifurcation and its Applications}.
	Springer-Verlag, 1976.

\item{[M]} J. Milnor,
	{\sl Dynamics in One Complex Variable: Introductory Lectures},
	Stony Brook Institute for Math. Sciences preprint \#1990/5.

\item{[S]} G. \'Swi\c{a}tek,
	{\sl Hyperbolicity is Dense in the Real Quadratic Family},
	Stony Brook Institute for Math. Sciences preprint \#1992/10.
 
\item{[W]} S. Wolfram,
	{\sl Mathematica: A System for Doing Mathematics by Computer},
	Addison-Wesley, 1988.
}

	\ifSingleSecs \vfil\eject \fi

\end